\documentclass{amsart}
\usepackage{amssymb}
\usepackage{tikz-cd}
\usepackage{CJKutf8}
\usepackage{amssymb, latexsym, MnSymbol, amscd ,graphicx,psfrag}
\usepackage[mathscr]{euscript}
\usepackage{chngcntr}

\usepackage[letterpaper, margin=1.4in]{geometry}
\usepackage{hyperref}
\hypersetup{
    colorlinks=true,
    linkcolor=black,
    citecolor=black,
    urlcolor=black
}

\usepackage{color, xcolor, enumitem, stmaryrd, url, verbatim}
\usepackage[all]{xy}

\def\orange { \textcolor{orange} }

\newcommand{\Ga}{\Gamma}
\newcommand{\Si}{ {\Sigma} }

\newcommand{\si}{ {\sigma} }

\newcommand{\bC}{ {\mathbb{C}} }

\newcommand{\bF}{\mathbb{F}}

\newcommand{\bh}{\mathbf{h}}

\newcommand{\bL}{\mathbb{L}}
\newcommand{\bP}{\mathbb{P}}
\newcommand{\bQ}{\mathbb{Q}}
\newcommand{\bR}{\mathbb{R}}

\newcommand{\bT}{\mathbb{T}}

\newcommand{\bZ}{\mathbb{Z}}

\newcommand{\cB}{\mathcal{B}}
\newcommand{\cC}{\mathcal{C}}
\newcommand{\cD}{\mathcal{D}}

\newcommand{\cF}{\mathcal{F}}
\newcommand{\cH}{\mathcal{H}}
\newcommand{\cI}{\mathcal{I}}
\newcommand{\cL}{\mathcal{L}}
\newcommand{\cM}{\mathcal{M}}
\newcommand{\cO}{\mathcal{O}}
\newcommand{\cP}{\mathcal{P}}
\newcommand{\cQ}{\mathcal{Q}}

\newcommand{\cX}{\mathcal{X}}
\newcommand{\cY}{\mathcal{Y}}
\newcommand{\cV}{\mathcal{V}}

\newcommand{\cZ}{\mathcal{Z}}

\newcommand{\inner}[1]{\langle  #1 \rangle}
\newcommand{\floor}[1]{\lfloor  #1 \rfloor}
\newcommand{\ceil}[1]{\lceil  #1 \rceil}

\newcommand{\age}{\mathrm{age}}
\newcommand{\Aut}{\mathrm{Aut}}
\newcommand{\CR}{ {\mathrm{CR}} }

\newcommand{\Hom}{\mathrm{Hom}}

\newcommand{\Spec}{\mathrm{Spec}}

\newcommand{\ev}{\mathrm{ev}}
\newcommand{\val}{ {\mathrm{val}} }
\newcommand{\vir}{ {\mathrm{vir}} }

\newcommand{\inv}{\mathrm{inv}}

\newcommand{\Eff}{{\mathrm{Eff}}}
\newcommand{\Nef}{{\mathrm{Nef}}}

\newcommand{\pt}{\Pic^\mathrm{pt}}

\newcommand{\bb}{\mathbf{b}}
\newcommand{\bd}{\mathbf{d}}

\newcommand{\bj}{\mathbf{j}}
\newcommand{\bk}{\mathbf{k}}
\newcommand{\bm}{\mathbf{m}}

\newcommand{\bw}{\mathbf{w}}
\newcommand{\one}{\mathbf{1}}

\newcommand{\bsi}{{\boldsymbol{\si}}}

\newcommand{\bSi}{\mathbf{\Si}}

\newcommand{\fa}{\mathfrak{a}}
\newcommand{\fb}{\mathfrak{b}}

\newcommand{\fg}{\mathfrak{g}}
\newcommand{\fl}{\mathfrak{l}}
\newcommand{\fm}{\mathfrak{m}}

\newcommand{\fn}{\mathfrak{n}}
\newcommand{\fp}{\mathfrak{p}}
\newcommand{\fq}{\mathfrak{q}}
\newcommand{\fr}{\mathfrak{r}}
\newcommand{\fs}{\mathfrak{s}}

\newcommand{\fx}{\mathfrak{x}}
\newcommand{\fy}{\mathfrak{y}}
\newcommand{\fo}{\mathfrak{o}}
\newcommand{\fu}{\mathfrak{u}}

\newcommand{\su}{\mathsf{u}}
\newcommand{\sv}{\mathsf{v}}
\newcommand{\sw}{\mathsf{w}}
\newcommand{\sfp}{\mathsf{p}}
\newcommand{\sff}{\mathsf{f}}

\newcommand{\sGa}{\mathsf{\Gamma}}

\newcommand{\hY}{\hat{Y}}

\newcommand{\hGa}{\hat{\Gamma}}

\newcommand{\hbeta}{\hat{\beta}}

\newcommand{\hD}{\hat{D}}

\renewcommand{\hY}{\hat{Y}}
\newcommand{\hZ}{\hat{Z}}

\newcommand{\hcY}{\hat{\cY}}
\newcommand{\hcZ}{\hat{\cZ}}
\newcommand{\hcD}{\hat{\cD}}

\newcommand{\tGa}{\widetilde{\Gamma}}

\newcommand{\tSi}{\widetilde{\Sigma}}
\newcommand{\tsi}{\widetilde{\sigma}}
\newcommand{\trho}{\widetilde{\rho}}

\newcommand{\tcX}{\widetilde{\cX}}

\newcommand{\tbh}{\widetilde{\bh}}

\newcommand{\tsu}{\widetilde{\su}}
\newcommand{\tsw}{\widetilde{\sw}}

\newcommand{\tM}{\widetilde{M}}
\newcommand{\tN}{\widetilde{N}}

\newcommand{\tT}{\widetilde{T}}

\newcommand{\tv}{\widetilde{v}}
\newcommand{\tX}{\widetilde{X}}

\newcommand{\tb}{\widetilde{b}}

\newcommand{\tk}{\widetilde{k}}

\newcommand{\tu}{\widetilde{u}}

\newcommand{\tbeta}{\widetilde{\beta}}

\newcommand{\tgamma}{\widetilde{\gamma}}
\newcommand{\ttau}{\widetilde{\tau}}

\newcommand{\tbT}{\widetilde{\bT}}
\newcommand{\tbw}{\widetilde{\bw}}

\newcommand{\tNef}{\widetilde{\Nef}}

\newcommand{\tbSi}{\widetilde{\bSi}}
\newcommand{\tcD}{\widetilde{\cD}}
\newcommand{\talpha}{\widetilde{\alpha}}

\newcommand{\vGa}{\vec{\Gamma}}
\newcommand{\vga}{\vec{\gamma}}
\newcommand{\vj}{\vec{j}}
\newcommand{\vd}{\vec{d}}
\newcommand{\vf}{\vec{f}}
\newcommand{\vk}{\vec{k}}
\newcommand{\vs}{\vec{s}}
\newcommand{\vmu}{\vec{\mu}}

\newcommand{\Mbar}{\overline{\cM}}

\newcommand{\mubar}{\overline{\mu}}

\newcommand{\ebar}{\overline{e}}

\newcommand{\bpsi}{\overline{\psi}}

\renewcommand{\pt}{\mathrm{pt}}
\renewcommand{\Box}{\mathrm{Box}}

\newcommand{\sm}{\mathrm{sm}}
\newcommand{\ord}{\mathrm{ord}}

\newtheorem{dummy}{dummy}[section]
\newtheorem{lemma}[dummy]{Lemma}
\newtheorem{theorem}[dummy]{Theorem}

\newtheorem{corollary}[dummy]{Corollary}
\newtheorem{proposition}[dummy]{Proposition}
\newtheorem{assumption}[dummy]{Assumption}

\theoremstyle{definition}
\newtheorem{definition}[dummy]{Definition}
\newtheorem{example}[dummy]{Example}

\theoremstyle{remark}
\newtheorem{remark}[dummy]{Remark}

\numberwithin{equation}{section}

\begin{document}
\counterwithout{figure}{section}

\counterwithout{equation}{section}
\counterwithout{equation}{subsection}
\title{Multi-component open/relative/local correspondence}

\author{Song Yu}
\address{Song Yu, Yau Mathematical Sciences Center, Tsinghua University, Haidian District, Beijing 100084, China}
\email{song-yu@tsinghua.edu.cn}

\author{Ke Zhang}
\address{Ke Zhang, Department of Mathematical Sciences, Tsinghua University, Haidian District, Beijing 100084, China}
\email{ke-zhang23@mails.tsinghua.edu.cn}

\author{Zhengyu Zong}
\address{Zhengyu Zong, Department of Mathematical Sciences, Tsinghua University, Haidian District, Beijing 100084, China}
\email{zyzong@mail.tsinghua.edu.cn}

\date{}

\begin{abstract}
For a toric Calabi-Yau 3-orbifold relative to $s$ Aganagic-Vafa outer branes, we prove a correspondence among the genus-zero open Gromov-Witten invariants with maximal winding at each brane and: (i) closed invariants of a toric Calabi-Yau ($3+s$)-orbifold; (ii) formal relative invariants of a formal toric Calabi-Yau (FTCY) 3-orbifold with maximal tangency to $s$ divisors; (iii) formal relative invariants of a sequence of FTCY intermediate geometries interpolating dimensions 3 and $3+s$. The correspondence provides examples of the log/local principle of van Garrel-Graber-Ruddat in the multi-component setting and the refined conjecture of Brini-Bousseau-van Garrel via intermediate geometries. It also establishes the multi-component case of the open/closed correspondence proposed by Lerche-Mayr and studied by Liu-Yu. As an application, we obtain examples of the conjecture of Klemm-Pandharipande on the integrality of BPS invariants of higher-dimensional toric Calabi-Yau manifolds. Along the way, we set the basic stages of the relative Gromov-Witten theory of higher-dimensional FTCY orbifolds, generalizing the case of smooth 3-folds by Li-Liu-Liu-Zhou.

\end{abstract}


\maketitle

\setcounter{tocdepth}{1}

\tableofcontents


\section{Introduction}

\subsection{Historical background and motivation}

\subsubsection{The log/local correspondence}\label{sec:log-local}
Let $\cY$ be a smooth projective variety and $\cD = \cD_1 + \cdots + \cD_s$ be a normal crossing divisor whose irreducible components $\cD_i$ are smooth and nef. The \emph{log/local correspondence}, proposed by van Garrel-Graber-Ruddat \cite{vGGR19}, conjectures that the genus-zero log Gromov-Witten invariants of $(\cY,\cD)$ with maximal tangency at each component $\cD_i$ should coincide with the local Gromov-Witten invariants of the total space of the vector bundle $\cO_{\cY}(-\cD_1)\oplus\cdots\oplus\cO_{\cY}(-\cD_s)$. When $\cD$ is smooth, i.e. $s=1$, \cite{vGGR19} provided a proof at the level of virtual fundamental cycles. The log/local correspondence has since been studied in different setups and generalities; see e.g. \cite{BBvG19,BBvG20,BBvG20b,BFGW21,BS23,CvGKT21,vGNS23,NR19,Schuler25,TY20}. Notably, Nabijou-Ranganathan \cite{NR19} gave a counterexample to the conjecture in the general setting where the pair $(\cY,\cD)$ is not log Calabi-Yau. On the other hand, the log Calabi-Yau case is supported by various studies, including the case of toric orbifolds relative to the toric boundary \cite{BBvG19}.

When $(\cY, \cD)$ is log Calabi-Yau, Brini-Bousseau-van Garrel \cite{BBvG20} proposed a refinement of the log/local correspondence in terms of \emph{intermediate geometries}. Specifically, for each $\ell = 0, \dots, s$, consider the log Calabi-Yau pair $(\cY^{(\ell)}, \cD^{(\ell)})$ with
$$
	\cY^{(\ell)} = \bigoplus_{i = 1}^\ell \cO_{\cY}(-\cD_i), \qquad \cD^{(\ell)} = \cD^{(\ell)}_{\ell+1} \cup \cdots \cup \cD^{(\ell)}_s,
$$
where $\cD^{(\ell)}_i$ is the preimage of $\cD_i$ under the projection $\cY^{(\ell)} \to \cY$. This sequence of geometries interpolate the initial pair $(\cY, \cD) = (\cY^{(0)}, \cD^{(0)})$ and the eventual local geometry $\cY^{(s)}$, and the genus-zero maximally-tangent log Gromov-Witten invariants of these geometries are expected to all be identified \cite[Conjecture 1.1]{BBvG20}. In particular, each intermediate step may be viewed as a non-compact instance of \cite{vGGR19}, where a log tangency condition is exchanged by a local condition.


\subsubsection{The open/closed correspondence}\label{sec:open-closed correspondence}
The open/closed correspondence, proposed by Mayr \cite{Mayr01} and Lerche-Mayr \cite{LM01}, conjectures that the disk invariants of a toric Calabi-Yau 3-fold $\cX$ relative to a Lagrangian $\cL$ of Aganagic-Vafa type should coincide with the genus-zero closed Gromov-Witten invariant of a dual toric Calabi-Yau 4-fold $\tcX$. When $\cL$ is an outer brane, the correspondence has been established mathematically by Liu and the first author \cite{LY21,LY22}, generalizing examples that arise from Looijenga pairs \cite{BBvG20,BBvG20b}. It has lead to new structural results for the open and closed Gromov-Witten theories, including open Witten-Dijkgraaf-Verlinde-Verlinde equations for $(\cX, \cL)$ \cite{YZ23} and BPS integrality structures \cite{Yu24}. The correspondence also admits a B-model, Hodge-theoretic counterpart that is compatible under open and closed mirror symmetry \cite{LY22,Yu25}.



As indicated in the proof when $\cX$ is smooth \cite{LY21}, the open/closed correspondence factors through the log/local correspondence in the following way. By the open/relative correspondence studied in \cite{FL13,FLT12}, the open Gromov-Witten invariants of $(\cX,\cL)$ coincide with the relative Gromov-Witten invariants of a pair $(\cY,\cD)$, where $\cY$ is a partial compactification of $\cX$ obtained by adding an irreducible toric divisor $\cD$ whose position depends on the position of the Aganagic-Vafa brane $\cL$. The corresponding toric Calabi-Yau 4-fold $\tcX$ is then given by the total space of $\cO_{\cY}(-\cD)$, and the genus-zero maximally-tangent relative invariants of $(\cY, \cD)$ are shown to coincide with the local Gromov-Witten invariants of $\tcX$. This provides a class of examples for the log/local correspondence in the extended setting where the base $\cY$ can be non-compact, in view of the identification of log and relative invariants \cite{AMW14}. The relative perspective also generalizes to the case where $\cX$ can be an orbifold, although not explicitly discussed in the proof in \cite{LY22}. We note that $\tcX$ may be replaced by its natural toric semi-projective partial compactification \cite{LY22}, which is convenient for the B-model correspondence and mirror symmetry considerations.



We note that the open/closed correspondence has been studied for other target geometries, including the quintic Calabi-Yau 3-fold \cite{AL23} and the projective line \cite{Zong25,YZ25}.

\subsubsection{The multi-component case}

\label{sec:open/relative/local}

In this paper, we consider the aforementioned chain of correspondences in the multi-component case, where the open geometry is a toric Calabi-Yau 3-orbifold $\cX$ relative to $s$ Aganagic-Vafa outer branes $\cL_1, \dots, \cL_s$. We will construct a toric Calabi-Yau ($3+s$)-orbifold $\tcX$ from $(\cX, \cL := \cL_1\sqcup \dots \sqcup\cL_s)$ as the corresponding closed geometry, and prove that the open Gromov-Witten invariants of $(\cX, \cL)$ coincide with the closed Gromov-Witten invariants of $\tcX$. This generalizes the single-component case ($s=1$) in \cite{LY21,LY22}.


As in \cite{FL13,FLT12}, the open Gromov-Witten invariants of $(\cX, \cL)$ can be explained as relative Gromov-Witten invariants. In the multi-component case, the natural relative geometry is a relative \emph{formal} toric Calabi-Yau 3-orbifold $(\hcY, \hcD = \hcD_1\sqcup \dots \sqcup\hcD_s)$, where $\hcD_1, \dots, \hcD_s$ are disjoint irreducible formal toric divisors. From the viewpoint of the open/relative correspondence, which identifies the open invariants of $(\cX, \cL)$ and the (formal) relative invariants of $(\hcY, \hcD)$, the Lagrangian boundary condition at each brane $\cL_i$ can be interpreted as tangency condition to the formal divisor $\hcD_i$. Moreover, based on $(\hcY, \hcD)$, we will construct a sequence of intermediate geometries
$(\hcY^{(\ell)}, \hcD^{(\ell)})$, $\ell = 0, \dots, s$, which is a relative formal toric Calabi-Yau ($3+\ell$)-orbifold. In particular, $(\hcY, \hcD) = (\hcY^{(0)}, \hcD^{(0)})$, and $\hcY^{(s)}$ is the formal completion of closed geometry $\tcX$ along its toric 1-skeleton. We give an identification of the (formal) relative Gromov-Witten invariants of all the intermediate geometries, verifying the multi-component log/local correspondence and the refined \cite[Conjecture 1.1]{BBvG20} in this setting. In particular, this gives a factorization of the multi-component open/closed correspondence into a sequence of instances of the log/local correspondence. In the single-component ($s=1$) orbifold case, this provides the missing relative perspective in \cite{LY22}.

We want to emphasize that our target geometries are \emph{non-compact} and the Gromov-Witten invariants are defined  \emph{equivariantly}. 
We also emphasize that working with orbifolds is essential in the multi-component case: even if $\cX$ is a smooth toric Calabi-Yau 3-fold, the corresponding geometries $\tcX$ and $\hcY^{(\ell)}$ are orbifolds in general when $s>1$ (see e.g. Remark \ref{rem:OrbifoldExample}). We note in addition that non-formal constructions for the relative and intermediate geometries are possible but require additional technicalities; we give an illustration in the example $\cX = \bC^3$ in Section \ref{sect:example}.

In the multi-component case, we may derive structural results similar to \cite{Yu24,YZ23} using the open/closed correspondence. In addition, we may replace $\tcX$ with its natural toric semi-projective partial compactification and develop the B-model correspondence in a way similar to \cite{LY22,Yu25}.

\subsection{Relative FTCY orbifolds and their Gromov-Witten theory}


The natural relative and intermediate geometries involved in the multi-component correspondences are relative \emph{formal} toric Calabi-Yau (FTCY) orbifolds. Smooth relative FTCY 3-folds and their formal relative Gromov-Witten theory are introduced by Li-Liu-Liu-Zhou \cite{LLLZ09} in the establishment of the mathematical theory of the topological vertex. The construction is motivated by the observation that the equivariant Gromov-Witten theory of a toric orbifold is determined by the equivariant information of the formal neighborhood of its toric 1-skeleton. The advantage of this perspective is that one may construct and study a formal geometry from a decorated graph, called the FTCY graph, which locally models the toric 1-skeleton in some toric orbifold, even if a global model for the graph does not exist. In this paper, we set up the details of relative FTCY orbifolds in a general dimension $r \ge 3$ and their formal relative Gromov-Witten theory, which may be of independent interest. We remark that concurrent work of Hu \cite{Hu} also establishes the theory for relative FTCY 3-orbifolds, as part of the development of the \emph{orbifold} topological vertex.
   
An FTCY graph will be a certain decorated graph with a local embedding into $\bR^{r-1}$, which assigns directions and rational weights to the edges in a locally consistent way. The construction in this paper compares to that in \cite[Section 3]{LLLZ09} in the following ways:
\begin{itemize}
    \item We extend the construction from dimension 3 to a general dimension $r$.
    
    \item We extend the construction to orbifolds.
    
    \item We place relative conditions only at \emph{univalent} vertices in the graph, which corresponds to the case of \emph{regular} graphs in \cite{LLLZ09} (see Definition 3.2).
\end{itemize}
We note that Hu \cite{Hu} treats the general 3-dimensional orbifold case where the graphs can admit \emph{bivalent} vertices. The general non-regular case in higher dimensions is left to future work. When the relative condition (or equivalently the set of univalent vertices) is empty, the FTCY graphs to be defined are generalizations of the toric graphs of genuine toric orbifolds considered in \cite[Section 8]{Liu13} and specializations of abstract stacky GKM graphs considered in \cite[Section 4]{LS20}. We note that in the more general setting of \cite{LS20}, one needs to assign additional data to the vertices and edges including isotropy groups and their representations at the flags. In the toric Calabi-Yau case, however, the additional data are redundant as they can be retrieved from the positioning of the edges. 

We also discuss the equivariant formal relative Gromov-Witten invariants of relative FTCY orbifolds, inspired by and generalizing the case of smooth relative FTCY 3-folds in \cite{LLLZ09}. To contain the scope of the exposition, we will only treat the case required for the application in this paper, where:
\begin{itemize}
    \item the stable maps have genus-zero domains;
    \item the stable maps have maximal tangency with each component of the formal divisor;
    \item each component of the formal divisor has cyclic generic stabilizer group.
\end{itemize}
We exemplify the computation of the formal relative invariants via virtual localization \cite{GP99, GV05, Liu13} by the invariants involved in our correspondences. The general case is again left to future work.

\subsection{Main result: multi-component Gromov-Witten correspondences}\label{sec:main result}
Let $\cX$ be a toric Calabi-Yau 3-orbifold (Section \ref{sect:TCY3}) and let $\cL = \cL_1 \sqcup \cdots \sqcup \cL_s$ be a disjoint union of $s$ Aganagic-Vafa outer branes (Section \ref{sect:AVBranes}) in $\cX$ satisfying Assumption \ref{assump:Branes}. We specify a choice of \emph{parallel framings} of the $s$ branes, determined by a generic parameter $f \in \bQ$. From this data, we construct a relative FTCY 3-orbifold $(\hcY, \hcD)$ from $(\cX, \cL)$ (Section \ref{sect:RelativeGeometry}). The formal toric divisor $\hcD$ has $s$ connected components
$$
	\hcD = \hcD_1 \sqcup \cdots \sqcup \hcD_s
$$
where the position of $\hcD_i$ is determined by that of $\cL_i$. We also construct a toric Calabi-Yau ($3+s$)-orbifold $\tcX$ (Section \ref{sect:LocalGeometry}). Moreover, we construct a sequence of intermediate geometries 
$$
	(\hcY^{(0)}, \hcD^{(0)}) = (\hcY, \hcD), (\hcY^{(1)}, \hcD^{(1)}), \dots, (\hcY^{(s)}, \hcD^{(s)})
$$
that interpolate the relative and local geometries (Section \ref{sect:IntermediateGeometry}). For $\ell = 1, \dots, s$, $(\hcY^{(\ell)}, \hcD^{(\ell)})$ is a relative FTCY ($3+\ell$)-orbifold and can be viewed as the geometric formal vector bundle associated to the rank-$\ell$ locally free sheaf
$$
	\bigoplus_{i = 1}^\ell \cO_{\hcY}(-\hcD_i)
$$
on the formal scheme $\hcY$. The formal divisor
$$
	\hcD^{(\ell)} = \hcD^{(\ell)}_{\ell+1} \sqcup \cdots \sqcup \hcD^{(\ell)}_{s}
$$
has $s-\ell$ connected components, where $\hcD^{(\ell)}_i$ is the preimage of the divisor $\hcD_i$ under the projection $\hcY^{(\ell)} \to \hcY$. At each intermediate step, $\hcY^{(\ell)}$ can be viewed as the geometric formal line bundle associated to the invertible sheaf
$$
	\cO_{\hcY^{(\ell-1)}}(-\hcD^{(\ell-1)}_\ell)
$$
on $\hcY^{(\ell-1)}$. Eventually, when $\ell = s$, $\hcD^{(s)}$ is empty and $\hcY^{(s)}$ is the formal completion of $\tcX$ along its toric 1-skeleton.



Let $X, L, \tX ,\hY, \hY^{(\ell)}$ denote the coarse moduli spaces of $\cX, \cL,\tcX,  \hcY, \hcY^{(\ell)}$ respectively. We take an effective curve class $\beta \in \Eff(X)$ of $X$ and winding data $\bd = ((d_1, \lambda_1), \dots, (d_s, \lambda_s))$, where $d_i \in \bZ_{\ge 1}$ records the winding number at the brane $\cL_i$ and $\lambda_i$ is an element of the generic stablizer group of $\cL_i$ recording the monodromy. The data determines an effective relative curve class $\beta' \in \Eff(X,L)$ and an effective curve class $\tbeta \in \Eff(\tX)$ of $\tX$. For the formal relative geometries, we consider effective curve classes $\hbeta \in \Eff(\hY)$ (which is identified with $\Eff(\hY^{(\ell)})$ for all $\ell$) that project to $\beta'$ under the natural projection $\pi: \Eff(\hY) \to \Eff(X,L)$ induced by the inclusion of the toric 1-skeleton. The winding data $\bd$ further determines the twisting data $\bk = (\tk_1, \dots, \tk_s)$ where $\tk_i$ is an element of the generic stabilizer group $\hcD_i$ recording the twisting. With the above topological data, we consider the following three classes of Gromov-Witten invariants:
\begin{itemize}
	\item The \emph{open} invariant $\langle \gamma_1,\dots,\gamma_n  \rangle^{\cX, \cL, T_f}_{\beta',\bd}$ of $(\cX, \cL, f)$ with maximal winding at each component $\cL_i$, where $\gamma_1, \dots, \gamma_n$ are (equivariant lifts of) insertions from the Chen-Ruan cohomology $H^2_{\CR}(\cX;\bQ)$.
	
	\item The \emph{formal relative} invariant $\langle \gamma_1^{(0)},\dots,\gamma_n^{(0)} \mid \gamma_{n+1}^{(0)},\dots,\gamma_{n+s}^{(0)} \rangle^{\hcY/\hcD, T_f}_{\hbeta,\bk}$ of $(\hcY, \hcD)$ with maximal tangency at each component $\hcD_i$, where for $i = 1, \dots, n$, $\gamma_i^{(0)}$ is a lift of $\gamma_i$, and for $i = n+1, \dots, n+s$, $\gamma_i^{(0)}$ is a new, fixed insertion to $\hcD_{i-n}$ in the twisting class $\tk_i$.
	
	\item The \emph{local} invariant $\langle \tgamma_1, \dots, \tgamma_{n+s} \rangle^{\tcX, T_f}_{\tbeta}$ of $\tcX$, where for $i = 1, \dots, n$, $\tgamma_i$ is a lift of $\gamma_i$, and for $i = n+1, \dots, n+s$, $\tgamma_i$ is a new, fixed insertion transformed from the relative insertion $\gamma_i^{(0)}$.
	
\end{itemize}
See Section \ref{sect:CorrData} for details of the corresponding topological data and insertions. Our first main result is an identification of the above classes of invariants.

\begin{theorem}[See Theorem \ref{thm:OpenClosedStatement}]\label{thm:main1}
We have the open/relative/local correspondence
$$
	\langle \gamma_1,\dots,\gamma_n  \rangle^{\cX, \cL, T_f}_{\beta',\bd} 
	= (-1)^{\sum_{i = 1}^s (\ceil{\frac{d_i}{\fa_i}} -1)} \sum_{\pi(\hbeta) = \beta'}\langle \gamma_1^{(0)},\dots,\gamma_n^{(0)} \mid \gamma_{n+1}^{(0)},\dots,\gamma_{n+s}^{(0)} \rangle^{\hcY/\hcD, T_f}_{\hbeta,\bk} 
	= \langle \tgamma_1, \dots, \tgamma_{n+s} \rangle^{\tcX, T_f}_{\tbeta}.
$$
\end{theorem}

The sign in front of the relative invariant is consistent with the prediction of the log/local correspondence (generalized to the orbifold setting), where the parameters $\fa_i \in \bZ_{\ge 1}$ are determined by the framing $f$. Note that the open/closed correspondence holds without a sign discrepancy.

We explain Theorem \ref{thm:main1} in terms of correspondences of Gromov-Witten invariants of the intermediate geometries, in alignment with \cite[Conjecture 1.1]{BBvG20}. For $\ell = 0, \dots, s$, we consider the \emph{formal relative} invariant $\langle \gamma_1^{(\ell)},\dots,\gamma_{n+\ell}^{(\ell)} \mid \gamma_{n+\ell+1}^{(\ell)},\dots,\gamma_{n+s}^{(\ell)} \rangle^{\hcY^{(\ell)}/\hcD^{(\ell)}, T_f}_{\hbeta,\bk^{(\ell)}}$ of $(\hcY^{(\ell)}, \hcD^{(\ell)})$ with maximal tangency at each component $\hcD^{(\ell)}_i$, where the first $\ell$ relative insertions have been transformed into insertions to the total space $\hcY^{(\ell)}$, and $\bk^{(\ell)} = (\tk_{\ell+1}, \dots, \tk_s)$ records the twisting at the remaining divisor components. We again refer to Section \ref{sect:CorrData} for details.

\begin{theorem}[See Theorem \ref{thm:RelativeStatement}]\label{thm:main2}
The correspondence among the formal relative invariants of the intermediate geometries hold. Inductively, we have
$$
	\langle \gamma_1^{(\ell)},\dots,\gamma_{n+\ell}^{(\ell)} \mid \gamma_{n+\ell+1}^{(\ell)},\dots,\gamma_{n+s}^{(\ell)} \rangle^{\hcY^{(\ell)}/\hcD^{(\ell)}, T_f}_{\hbeta,\bk^{(\ell)}} 
	= (-1)^{\ceil{\frac{d_{\ell+1}}{\fa_{\ell+1}}} -1} \langle \gamma_1^{(\ell+1)},\dots,\gamma_{n+\ell+1}^{(\ell+1)} \mid \gamma_{n+\ell+2}^{(\ell+1)},\dots,\gamma_{n+s}^{(\ell+1)} \rangle^{\hcY^{(\ell+1)}/\hcD^{(\ell+1)}, T_f}_{\hbeta,\bk^{(\ell+1)}}
$$
for all $\ell = 0, \dots, s-1$.
\end{theorem}

Each intermediate step in Theorem \ref{thm:main2} can be viewed as an instance of the log/local correspondence, where the relative condition at the divisor $\hcD^{(\ell)}_{\ell+1}$ is exchanged for a local condition in $\hcY^{(\ell+1)}$. Theorem \ref{thm:main2} is proved by a detailed analysis and comparison of the virtual localization \cite{GP99, GV05, Liu13} computation of the invariants in the orbifold setting with a mixture of relative and local insertions.


\subsection{Application to BPS integrality}
We apply the multi-component open/closed correspondence in Theorem \ref{thm:main1} to study BPS integrality, which generalizes the single-component ($s=1$) case in \cite{Yu24}. We give a brief summary here and leave the details to Section \ref{sect:BPS}.

For smooth Calabi-Yau 3-folds, \emph{BPS invariants} may be defined from the Gromov-Witten invariants via a famous resummation formula of Gopakumar-Vafa \cite{GV98a,GV98b} and are expected to have integrality and finiteness properties. Such properties have been extensively studied \cite{BP01, Peng07, Konishi06a,Konishi06b, IP18, DIW21}; in particular, Ionel-Parker \cite{IP18} provided a proof for integrality for general compact Calabi-Yau 3-folds. In genus zero, the Gopakumar-Vafa formula recovers the Aspinwall-Morrison multiple-cover formula \cite{AM93} and is generalized to dimensions four and above by Klemm-Pandharipande \cite{KP08}, who conjectured that the resulting BPS invariants are integers. This is verified in the compact case by Ionel-Parker \cite{IP18} and known in examples in the non-compact case in dimensions four \cite{KP08,CMT18,Cao20,CMT22,CKM22,COT22,COT24,BBvG20,BS23,Yu24}, five \cite{PZ10,BBvG20}, and six \cite{BBvG20}.

Another generalization of the Gopakumar-Vafa conjecture concerns open Gromov-Witten theory of Calabi-Yau 3-folds. In the toric case relative to Aganagic-Vafa branes, as considered in this paper, open BPS invariants may be defined from the Gromov-Witten invariants via a resummation formula of Labastida, Mari\~no, Ooguri, and Vafa (LMOV) \cite{OV00,LM00,LMV00,MV02} discovered in the context of Chern-Simons theory and large $N$ duality. Integrality and finiteness properties of these open BPS invariants are verified in examples by \cite{LZ16,LZ19,Zhu19,BBvG20,Zhu22} and in full in \cite{Yu24}. It is proposed that the LMOV formula may be applied to study the integrality structures of Gromov-Witten invariants of more general open geometries \cite[Section 1.3]{Yu24}.

As observed by Brini-Bousseau-van Garrel \cite{BBvG20}, via a comparison of the resummation formulas, the open/closed correspondence of Gromov-Witten invariants may be used to deduce a correspondence of the BPS invariants. BPS integrality may then be carried from the open side to the closed side. With this approach, \cite{BBvG20} obtained examples for the Klemm-Pandharipande integrality conjecture \cite[Conjecture 0]{KP08} when the (non-compact) open geometry $(\cX, \cL, f)$ and closed geometry $\tcX$ arise from Looijenga pairs. Applying the approach to a general smooth single-brane open geometry, \cite{Yu24} obtained a general family of non-compact examples for \cite[Conjecture 0]{KP08} in dimension four; see \cite[Section 4.3]{Yu24} for a detailed comparison of these examples with the previously-known ones.

The same approach may be applied to our present study of open geometries with multiple branes, in the smooth case. As a result, we verify \cite[Conjecture 0]{KP08} in a general class of non-compact toric Calabi-Yau manifolds in each dimension that is at least five. Such examples may arise from an arbitrary toric Calabi-Yau 3-fold and may not necessarily be a local curve or local surface.

\subsection{Outline of the paper}
In Sections \ref{sect:FTCY} and \ref{sec:GW-FTCY}, we set up the definitions of FTCY graphs, their associated relative FTCY orbifolds, and their Gromov-Witten theory. In Section \ref{sect:Geometry}, we study the open geometry $(\cX, \cL, f)$ and construct the corresponding closed geometry $\tcX$ and the formal intermediate geometries $(\hcY^{(\ell)}, \hcD^{(\ell)})$, $\ell = 0, \dots, s$. In Section \ref{sec:GW}, we provide the localization computations of the open, formal relative, and closed Gromov-Witten invariants involved in the correspondences. In Section \ref{sect:Correspondence}, we prove the Gromov-Witten correspondences in Theorems \ref{thm:main1}, \ref{thm:main2}. We discuss the application to BPS integrality in Section \ref{sect:BPS}. Finally, we discuss the example $\cX = \bC^3$ in Section \ref{sect:example}.

\subsection*{Acknowledgments}
We thank Andrea Brini, Bohan Fang, Hechen Hu, and Chiu-Chu Melissa Liu for valuable discussions and constructive feedback. The third author is partially supported by the Natural Science Foundation of Beijing, China (Grant No. 1252008) and NSFC (Grant No. 12571067).


\section{Relative formal toric Calabi-Yau orbifolds}\label{sect:FTCY}
In this section, we include the preliminaries of toric orbifolds with a focus on the Calabi-Yau case. Then, we define formal toric Calabi-Yau (FTCY) graphs and their associated relative FTCY orbifolds. We work over $\bC$ throughout.

\subsection{Toric orbifolds}\label{sect:ToricOrbifolds}
We start by reviewing the basics of \emph{toric orbifolds}, or \emph{smooth toric Deligne-Mumford stacks} with trivial generic stabilizer, and introducing some notation. We refer to \cite{CLS11,Fulton93} for the general theory of toric varieties and to \cite{BCS05, FMN10} for the general theory of smooth toric Deligne-Mumford stacks.

\subsubsection{Extended stacky fan}
Let $\cZ$ be an $r$-dimensional toric orbifold specified by an \emph{extended stacky fan} $\mathbf{\Xi} = (N, \Xi, \alpha)$ in the sense of \cite{Jiang08}, where $N \cong \bZ^r$ is a rank-$r$ lattice, $\Xi$ is a finite simplicial fan in $N_{\bR} \cong \bR^r$, and $\alpha:\bZ^R \to N$ is a group homomorphism determined by a list of vectors $(b_1, \dots, b_R)$ in $N$.\footnote{In this paper, for a lattice $\Lambda$ and $\bF = \bQ, \bR$, or $\bC$, we set $\Lambda_{\bF} := \Lambda \otimes \bF$.} The coarse moduli space $Z$ of $\cZ$ is the simplicial toric variety defined by the fan $\Xi$. Since $\cZ$ is an orbifold, the Deligne-Mumford torus $T$ acting on $\cZ$ is also the dense algebraic torus acting on $Z$, which is $N \otimes \bC^* \cong (\bC^*)^r$.

For each $d = 0, \dots, r$, let $\Xi(d)$ denote the set of $d$-dimensional cones in $\Xi$. In particular, there exists $1 \leq R' \leq R$ such that
$$
   \Xi(1) = \{\bR_{\ge 0}b_1, \dots, \bR_{\ge 0}b_{R'}\}.
$$
If $\alpha':\bZ^{R'} \to N$ is the group homomorphism determined by $(b_1, \dots, b_{R'})$, then the triple $(N, \Xi, \alpha')$ is the \emph{stacky fan} of $\cZ$ in the sense of \cite{BCS05}.

For each $\sigma \in \Xi(d)$, let $\cV(\sigma) \subseteq \cZ$ denote the codimension-$d$ $T$-invariant closed substack of $\cZ$ corresponding to $\sigma$. Let $V(\sigma) \subseteq Z$ denote the codimension-$d$ $T$-orbit closure in $Z$ corresponding to $\sigma$, which is the coarse moduli space of $\cV(\sigma)$. Let
$$
    \iota_{\sigma}: \cV(\sigma) \to \cZ, \quad V(\sigma) \to Z,
$$
denote the inclusion maps.\footnote{By an abuse of notation, in this paper, the letter $\iota$ will be used to denote various natural inclusion maps of substacks/subvarieties, fixed loci of moduli spaces, or cones. The precise meaning and usage will be made clear in the context.}

\subsubsection{stabilizers}
Let $\sigma \in \Xi(d)$. We set index sets
$$
    I_\sigma' := \{ i \in \{1, \dots, R'\} : \rho_i \subseteq \sigma \}, \qquad I_\sigma := \{1, \dots, R\} \setminus I_\sigma'.
$$
Note that $|I_\sigma'|=d$. The generic stabilizer group of the substack $\cV(\sigma)$, denoted $G_\sigma$, is a finite abelian group and can be identified as
\begin{equation}\label{eqn:ConeStab}
    G_\sigma \cong \left. \left( N \cap \sum_{i \in I_\sigma'} \bR b_i \right) \middle/ \sum_{i \in I_\sigma'} \bZ b_i \right. .
\end{equation}
We define
$$
    \Box(\sigma) := \left\{ v \in N: v = \sum_{i \in I_\sigma'} c_ib_i \text{ for some } 0 \le c_i < 1 \right\},
$$
which gives a set of representatives for $G_\sigma$. Given $v = \sum_{i \in I_\sigma'} c_i(v)b_i \in \Box(\sigma)$, we define
$$
    \age(v) := \sum_{i \in I_\sigma'} c_i(v).
$$
Given cones $\tau \subseteq \sigma$ in $\Xi$, we have natural inclusions $G_\tau \subseteq G_\sigma$, $\Box(\tau) \subseteq \Box(\sigma)$.
\subsubsection{Fixed points, torus-invariant lines, fundamental groups, flags}
For each $\sigma \in \Xi(r)$, let $\fp_\sigma := \cV(\sigma)$ denote the corresponding $T$-fixed point in $\cZ$, and $p_\sigma:= V(\sigma)$ denote the corresponding $T$-fixed point in $Z$. For each $\tau \in \Xi(r-1)$, let $\fo_\tau \cong \bC^* \times \cB G_\tau$ denote the corresponding (open) $T$-orbit in $\cZ$, and $o_\tau \cong \bC^*$ denote the corresponding $T$-orbit in $Z$. Let $\fl_\tau:= \cV(\tau)$ denote the corresponding closed $T$-invariant line in $\cZ$, which is the closure of $\fo_\tau$, and $l_\tau := V(\tau)$ denote the corresponding closed $T$-invariant line in $Z$, which is the closure of $o_\tau$. We set
$$
    \Xi(r-1)_c:= \{\tau \in \Xi(r-1) : l_\tau \text{ is compact} \},
$$
and define
$$
    \cZ^1 := \bigcup_{\tau \in \Xi(r-1)} \fl_{\tau}, \quad Z^1:= \bigcup_{\tau \in \Xi(r-1)} l_{\tau}, \quad \cZ^1_c := \bigcup_{\tau \in \Xi(r-1)_c} \fl_{\tau}, \quad Z^1_c:= \bigcup_{\tau \in \Xi(r-1)_c} l_{\tau}.
$$

For each $\tau \in \Xi(r-1)$, let $H_\tau := \pi_1(\fo_\tau)$ be the fundamental group of $\fo_\tau$. The projection $\fo_\tau \to o_\tau$ to the coarse moduli space induces a map
$$
    \pi_\tau: H_\tau = \pi_1(\fo_\tau) \to \pi_1(o_\tau) \cong \bZ
$$
on the fundamental groups, which fits into a split short exact sequence
$$
    \xymatrix{
        1 \ar[r] & G_{\tau} \ar[r] & H_{\tau} \ar[r]^{\pi_\tau} \ar[r] & \bZ \ar[r] &  1.
    }
$$

Let
$$
    F(\Xi):= \{ (\tau, \sigma) \in \Xi(r-1) \times \Xi(r) : \tau \subset \sigma\}
$$
be the set of \emph{flags} in $\Xi$. Given a flag $(\tau, \sigma) \in F(\Xi)$, we have $\fp_\sigma \subset \fl_\tau$ and $p_\sigma \in l_\tau$. Let
$$
\chi_{(\tau, \sigma)}: G_\sigma \to \bC^*
$$
be the representation of $G_\sigma$ on the tangent line $T_{\fp_\sigma}\fl_\tau$. The image of $\chi_{(\tau, \sigma)}$ is $\mu_{\fr(\tau, \sigma)}$, where
$$
    \fr(\tau, \sigma):= \frac{|G_\sigma|}{|G_\tau|},
$$
and for $r \in \bZ_{>0}$, $\mu_r \subset \bC^*$ is the cyclic group of $r$-th roots of unity. The map $\chi_{(\tau, \sigma)}$ and the inclusion $G_\tau \to G_\sigma$ fit into a short exact sequence
\begin{equation}\label{eqn:FlagSeq}
    \xymatrix{
        1 \ar[r] & G_{\tau} \ar[r] & G_{\sigma} \ar[r]^{\chi_{(\tau,\sigma)}} & \mu_{\mathfrak{r}(\tau,\sigma)} \ar[r] & 1.
    }    
\end{equation}

 Let $\fu_{(\tau, \sigma)}:= \fo_\tau \cup \fp_\sigma$, which is an open substack of $\fl_\tau$. The inclusion $\fo_\tau \to \fu_{(\tau, \sigma)}$ induces a surjective map
\begin{equation}\label{eqn:FlagFundGroup}
    \pi_{(\tau, \sigma)}: H_\tau = \pi_1(\fo_\tau) \to \pi_1(\fu_{(\tau, \sigma)}) \cong G_\sigma
\end{equation}
on the fundamental groups.


\subsubsection{Chen-Ruan orbifold cohomology}
Let
$$
    \Box(\cZ) := \bigcup_{\text{cone $\sigma$ in $\Xi$}} \Box(\sigma) = \bigcup_{\text{maximal cone $\sigma$ in $\Xi$}} \Box(\sigma),
$$
which indexes the inertia components of $\cZ$. The inertia stack of $\cZ$ is
$$
    \cI \cZ = \bigsqcup_{j \in \Box(\cZ)} \cZ_j.
$$
In particular, $\cZ_{\vec{0}} = \cZ$ is the untwisted sector. Let
$$
    \inv^*: \cI \cZ \to \cI \cZ
$$
denote the involution on $\cI \cZ$ which $(z,g)$ with $z \in \cZ, g \in \Aut(z)$ to $(z, g^{-1})$.

As a graded vector space over $\bQ$, the \emph{Chen-Ruan cohomology group} \cite{CR04} of $\cZ$ is defined as
$$
    H_{\CR}^*(\cZ; \bQ) := \bigoplus_{j \in \Box(\cZ)} H^*(\cZ_j; \bQ)[2\age(j)],
$$
where $[2\age(j)]$ denotes a degree shift by $2\age(j)$. We write $\one_j$ for the unit of $H^*(\cZ_j; \bQ)$, viewed as an element of $H_{\CR}^{2\age(j)}(\cZ; \bQ)$. In addition, for any subtorus $Q \subseteq T$, the \emph{$Q$-equivariant} Chen-Ruan cohomology group of $\cZ$ is
$$
    H_{\CR, Q}^*(\cZ; \bQ) := \bigoplus_{j \in \Box(\cZ)} H^*_Q(\cZ_j; \bQ)[2\age(j)],
$$
which is a module over $H^*_{Q}(\pt;\bQ)$. The above definitions can be extended to $\bC$-coefficients.

\subsection{The Calabi-Yau case}\label{sect:CYCase}
Consider the case $\cZ$ is \emph{Calabi-Yau}, i.e. the canonical bundle $K_Z$ of the coarse moduli space $Z$ is trivial. Let $M:= \Hom(N, \bZ)$, which is canonically identified with the character lattice $\Hom(T, \bC^*)$ of $T$. The Calabi-Yau condition implies the existence of a primitive vector (or character) 
$$
    \su \in M
$$
such that $\inner{\su, b_i} = 1$ for all $i = 1, \dots, R'$.\footnote{In this paper, we use $\inner{-,-}$ to denote the natural pairings between dual lattices or vector spaces.} That is, $b_1, \dots, b_{R'}$ belong to $N' \times \{1\}$ where $N':= \ker(\su) \subset N$.

Let $T':= \ker(\su) = N' \otimes \bC^* \cong (\bC^*)^{r-1}$ be the Calabi-Yau subtorus of $T$. The action of $T'$ on $\cZ$ has the same fixed points and invariant lines as the action of $T$. The character lattice $\Hom(T', \bC^*)$ or $T'$ is canonically identified with $M' := \Hom(N', \bZ)$, which comes with a natural projection $M \to M'$.

\subsection{Local data at a vertex of an FTCY graph}\label{sect:VertexData}
Now we turn to defining FTCY graphs. Consider the case where the dimension is $r \ge 3$. Similar to Sections \ref{sect:ToricOrbifolds}, \ref{sect:CYCase}, we set up the following notation to be used through the rest of Section \ref{sect:FTCY}. Let $N \cong \bZ^r$, $T := N \otimes \bC^*$, and $M := \Hom(N, \bZ)$. Fix a primitive vector
$$
    \su \in M.
$$
Let $N':= \ker(\su)$, $T' := N' \otimes \bC^*$, and $M' := \Hom(N', \bZ)$, as in Section \ref{sect:CYCase}. We further fix a splitting
$$
    N \cong N' \oplus \bZ \sv
$$
for some $\sv \in N$ satisfying $\inner{\su, \sv} = 1$.


The local data that will be associated to a vertex in an FTCY graph consist of a collection
$$
    \sw_1, \dots, \sw_{r-1} \in M'_{\bQ} \setminus \{0\}
$$
of linearly independent vectors such that 
$$
    M'_v := \sum_{i = 1}^{r-1} \bZ \sw_i
$$
is a rank-$(r-1)$ lattice in $M'_{\bQ}$ that contains $M'$ as a sublattice. We will consider two cases corresponding to two choices of an $r$-th vector as follows:
$$
    \sw_r := \begin{cases}
        -\sw_1 - \sw_2 - \cdots - \sw_{r-1} & \text{in the $r$-valent case,}\\
        -\sw_2 - \cdots - \sw_{r-1} & \text{in the univalent case.}
    \end{cases}
$$
To uniformize the description of the two cases, we introduce the indicator symbol
\begin{equation}\label{eqn:IndicatorSymbol}
    \spadesuit :=  \begin{cases}
        1 & \text{in the $r$-valent case,}\\
        0 & \text{in the univalent case.}
    \end{cases}
\end{equation}
Then
\begin{equation}\label{eqn:WRelation}
    \sw_r = -\spadesuit \sw_1 - \sw_2 - \cdots - \sw_{r-1}.
\end{equation}
The set of vectors
$$
    \{\sw_1, \dots, \sw_{r-1}, \sw_r + \su\}
$$
generates the rank-$r$ lattice $M_v = M'_v \oplus \bZ \su$ in $M_{\bQ}$ which contains $M$ as a sublattice. The set also forms a $\bQ$-basis of $M_{\bQ}$. Equation \eqref{eqn:WRelation} then implies that the dual $\bQ$-basis of $N_{\bQ}$ takes form
$$
    \{b_1 = b'_1 + \spadesuit \sv, b_2 = b'_2 + \sv, \dots, b_{r-1} = b'_{r-1} + \sv, b_r = \sv\}
$$
for some $b'_1, \dots, b'_{r-1} \in N'_{\bQ}$. Since $M_v \supseteq M$, the vectors $b_1, \dots, b_r$ lie in $N$ and generate a sublattice $N_v \subseteq N$ that is dual to $M_v$.

The convex hull of $b_1, \dots, b_r$ in $N_\bR$ is a simplicial $r$-cone, which we denote by $\sigma_v$. Let
$$
    \cZ_v \cong [\bC^r/G_v]
$$
be the canonical $r$-dimensional affine toric orbifold of the cone $\sigma_v$, where $G_v := N/N_v$ is a finite abelian group (cf. \eqref{eqn:ConeStab}). The coarse moduli space $Z_v$ of $\cZ_v$ is the corresponding simplicial toric variety. For $i = 1, \dots, r$, the ray generated by $b_i$ corresponds to a toric divisor $\cD_i \subset \cZ_v$ whose coarse moduli space is a toric divisor $D_i \subset Z_v$.  The Deligne-Mumford torus of $\cZ_v$ is $T$ and contains $T'$ as a codimension-1 subtorus.


The group $G_v$ is the generic stabilizer group of the unique $T$-fixed point $\fp_v := \cV(\sigma_v)$ in $\cZ_v$. For $i = 1, \dots, r$, let $\tau_i$ be the facet of $\sigma_v$ generated by the vectors $b_1, \dots, \widehat{b_i}, \dots b_r$, which corresponds to a $T$-invariant line $\fl_i := \cV(\tau_i)$ in $\cZ_v$. The rigidification of $\fl_i$ is isomorphic to $[\bC/\mu_{\fr_i}]$ and the coarse moduli is isomorphic to $\bC$. The generic stabilizer group $G_i$ of $\fl_i$ fits into a short exact sequence
\begin{equation}\label{eqn:FTCYFlagSeq}
    \xymatrix{
        0 \ar[r] & G_i \ar[r] & G_v \ar[r]^{\chi_i} & \mu_{\fr_i} \ar[r] & 1
    }
\end{equation}
where $\fr_i := \frac{|G_v|}{|G_i|}$ (cf. \eqref{eqn:FlagSeq}). The map $\chi_i: G_v \to \mu_{\fr_i} \subset \bC^*$ is the 1-dimensional representation of $G_v$ on the tangent line $T_{\fp_v}\fl_i$. We note that $\fr_i \sw_i \in M'$. Moreover, each $\sw_i$ is the weight of the $T'$-action on $T_{\fp_v}\fl_i$.

Finally, let
$$
    \hcZ_v
$$
denote the formal completion of $\cZ_v$ along the toric 1-skeleton $\cZ_v^1 = \bigcup_{i = 1}^r \fl_i$. In the univalent (or $\spadesuit = 0$) case, let
$$
    \hcD_v
$$
denote the formal completion of the divisor $\cD_{1}$ along the toric 1-skeleton $\bigcup_{i =2}^r \fl_{i}$, which is a divisor in $\hcZ_v$.

\begin{remark}\rm{
As we will see in the definition of FTCY graphs below, the data of vectors $\sw_1, \dots, \sw_r$ corresponds to \emph{position map} (or \emph{framing map} in the univalent case) in the sense of \cite{LLLZ09} or the \emph{axial function} in the sense of \cite{LS20} at the vertex. 
}
\end{remark}

\subsection{Local data at a compact edge}\label{sect:EdgeData}
The local data that will be associated to a compact edge in an FTCY graph consist of two collections of vectors
$$
    \sw_{1,1}, \dots, \sw_{1,r}, \quad \sw_{2,1}, \dots, \sw_{2,r} \in M'_{\bQ} \setminus \{0\}
$$
arising in the vertex situation in Section \ref{sect:VertexData}. Specifically, they satisfy that
\begin{equation}\label{eqn:DoubleWRelation}
    \sw_{1,1} + \cdots + \sw_{1,r} = 0, \qquad \sw_{2,1} + \cdots + \sw_{2,r} = (1 - \spadesuit)\sw_{2,1}
\end{equation}
where $\spadesuit$ takes the same meaning as in \eqref{eqn:IndicatorSymbol}, and that each collection generates a rank-$(r-1)$ lattice in $M'_{\bQ}$ that contains $M'$. Let $\fr_{1,1}, \fr_{2,1} \in \bZ_{\ge 1}$ be constructed as in Section \ref{sect:VertexData} and associated to $\sw_{1,1}, \sw_{2,1}$ respectively. We further require that
\begin{equation}\label{eqn:BaseConnectionCond}
    \fr_1 \sw_{1,1} + \fr_{2,1} \sw_{2,1} = 0,
\end{equation}
and that for $i = 2, \dots, r$,
\begin{equation}\label{eqn:FiberConnectionCond}
    \sw_{1,i} - \sw_{2,i} =  a_i \fr_{1,1} \sw_{1,1} = -a_i \fr_{2,1} \sw_{2,1}
\end{equation}
for some $a_i \in \bQ$.

Applying the dual basis construction as in Section \ref{sect:VertexData} and using the conditions \eqref{eqn:BaseConnectionCond}, \eqref{eqn:FiberConnectionCond}, we see that there exist lattice vectors $b_{1,1}', b_{2,1}', b_2', \dots, b_{r-1}' \in N'$ such that the sets
\begin{equation}\label{eqn:TwoVertexBases}
    \{b_{1,1} = b'_{1,1} + \sv, b_2 = b'_2 + \sv, \dots, b_{r-1} = b'_{r-1} + \sv, b_r = \sv\}, \quad \{b_{2,1} = b'_{2,1} + \spadesuit \sv, b_2, \dots, b_{r-1}, b_r\}
\end{equation}
are the $\bQ$-bases of $N_{\bQ}$ dual to the $\bQ$-bases
$$
    \{\sw_{1,1}, \dots, \sw_{1,r-1}, \sw_{1,r} + \su\}, \quad \{\sw_{2,1}, \dots, \sw_{2,r-1}, \sw_{2,r} + \su\}
$$
of $M_{\bQ}$ respectively. Let $\sigma_{v_1}$, $\sigma_{v_2}$ denote the simplicial $r$-cones spanned by the two sets in \eqref{eqn:TwoVertexBases} respectively. The two cones intersect at the facet $\tau_e$ spanned by $\{b_2, \dots, b_r\}$ and together fit into a simplicial fan in $N_{\bR}$. Here $\tau_e$ is identified with the facets $\tau_{1,1}$, $\tau_{2,1}$ constructed as in Section \ref{sect:VertexData}. Let
$$
    \cZ_e = \cZ_{v_1} \cup \cZ_{v_2}
$$
be the canonical $r$-dimensional toric orbifold of the fan, where the toric charts
$$
    \cZ_{v_1} \cong [\bC^r/G_{v_1}], \qquad \cZ_{v_2} \cong [\bC^r/G_{v_2}]
$$
are constructed as in Section \ref{sect:VertexData}. The orbit $\fl_e := \cV(\tau_e)$ in $\cZ_e$ corresponding to $\tau_e$, which is the union of $\fl_{1,1} \subset \cZ_{v_1}$ and $\fl_{2,1} \subset \cZ_{v_2}$, is a 1-dimensional complete toric Deligne-Mumford stack whose rigidification is isomorphic to the toric orbifold $\cC_{\fr_{1,1}, \fr_{2,1}}$ defined in \cite[Example 95]{Liu13} and whose coarse moduli space is isomorphic to $\bP^1$. The generic stabilizer group $G_e$ of $\fl_e$ is identified with the groups $G_{1,1}$, $G_{2,1}$ as constructed in Section \ref{sect:VertexData}.

In fact, $\cZ_e$ is the total space of the direct sum of $r-1$ orbifold line bundles on $\fl_e$ with degrees $a_2, \dots, a_r$ respectively. Equations \eqref{eqn:DoubleWRelation}, \eqref{eqn:BaseConnectionCond}, \eqref{eqn:FiberConnectionCond} together imply that
$$
    a_2 + \cdots + a_r = -\frac{1}{\fr_{1,1}} - \frac{\spadesuit}{\fr_{2,1}}.
$$
It follows that
$$
     K_{Z_e} \cong \cO_{Z_e}((1-\spadesuit) D_{2,1})
$$
where $Z_e$ denotes the coarse moduli space of $\cZ_e$. In other words, in the $r$-valent (or $\spadesuit=1$) case, $\cZ_e$ is Calabi-Yau, and in the univalent (or $\spadesuit = 0$) case, $(\cZ_e, \cD_{2,1})$ is log Calabi-Yau.

Finally, let
$$
    \hcZ_e
$$
denote the formal completion of $\cZ_e$ along the toric 1-skeleton $\cZ_e^1 = \fl_{\tau_e} \cup \bigcup_{i = 2}^r (\fl_{1,i} \cup \fl_{2,i})$. We have
$$
    \hcZ_e = \hcZ_{v_1} \cup \hcZ_{v_2}.
$$
In the univalent (or $\spadesuit = 0$) case, consider the divisor $\hcD_{v_2}$ of $\hcZ_{v_2}$ defined at the end of Section \ref{sect:VertexData} as a formal completion of $\cD_{2,1}$. Then $\hcD_{v_2}$ is a divisor in $\hcZ_e$ such that $(\hcZ_e, \hcD_{v_2})$ is log Calabi-Yau.

\subsection{Definition of FTCY graphs}\label{sect:FTCYGrahDef}
Now we are ready to define FTCY graphs. For the underlying graph, consider a graph $\hGa = (V(\hGa), E(\hGa))$ with finitely many vertices and edges, where $V(\hGa)$, $E(\hGa)$ denote the sets of vertices and edges respectively. Each edge in $E(\hGa)$ is either a compact edge containing two distinct vertices, or a ray containing only one vertex. Let $E(\hGa)_c \subset E(\hGa)$ denote the subset of compact edges. The set of flags in $\hGa$ is
$$
    F(\hGa) := \{(e,v) \in E(\hGa) \times V(\hGa) : v \in e\}.
$$
For each $v \in V(\hGa)$, let
$$
    E_v := \{e \in E(\hGa) : (e,v) \in F(\hGa)\}.
$$

\begin{definition}
A \emph{regular weakly $r$-valent} graph is a connected graph $\hGa = (V(\hGa), E(\hGa))$ as above such that $V(\hGa)$ admits a partition
$$
    V(\hGa) = V^r(\hGa) \sqcup V^1(\hGa)
$$
with $V^r(\hGa)$ non-empty, where
$$
    V^r(\hGa) := \{v \in V(\hGa) : |E_v| = r\}, \qquad V^1(\hGa) := \{v \in V(\hGa) : |E_v| = 1\}
$$
are the sets of $r-$valent and univalent vertices respectively.
\end{definition}

Note that every univalent vertex in a regular weakly $r$-valent graph is adjacent to a unique $r$-valent vertex.

\begin{definition}\label{def:FTCYGraph}
An \emph{FTCY graph} is a tuple $\sGa = (\hGa, \sfp, \sff_2, \dots, \sff_r)$ of a regular weakly $r$-valent graph $\hGa$ decorated with a \emph{position map}
$$
    \sfp: F(\hGa) \to M'_{\bQ} \setminus \{0\}
$$
and \emph{framing maps}
$$
    \sff_2, \dots, \sff_r: V^1(\hGa) \to M'_{\bQ} \setminus \{0\}
$$
such that the following conditions are satisfied:
\begin{itemize}

    \item At each $r$-valent vertex $v \in V^r(\hGa)$, the vectors $\{\sfp(e,v) : e \in E_v\}$ satisfy
    $$
        \sum_{e \in E_v} \sfp(e,v) = 0,
    $$
    and generate a rank-$(r-1)$ lattice in $M'_{\bQ}$ that contains $M'$ as a sublattice. 
    
    \item At each univalent vertex $v \in V^1(\hGa)$, the framing maps satisfy
    $$
        \sff_2(v) + \cdots + \sff_r(v) = 0.
    $$
    Moreover, let $E_v = \{e\}$. Then the vectors $\{\sfp(e,v), \sff_2(v), \dots, \sff_r(v)\}$ generate a rank-$(r-1)$ lattice in $M'_{\bQ}$ that contains $M'$ as a sublattice. 

\end{itemize}
From the discussion of Section \ref{sect:VertexData} and in particular \eqref{eqn:FTCYFlagSeq}, each flag $(e,v) \in F(\hGa)$ is assigned an integer $\fr{(e,v)} \in \bZ_{\ge 1}$.
\begin{itemize}
    \item For each compact edge $e \in E(\hGa)_c$ that contains two $r$-valent vertices $v', v'' \in V^r(\hGa)$, we have
    $$
        \fr{(e,v')} \sfp(e,v') + \fr{(e,v'')} \sfp(e,v'') = 0.
    $$
    Moreover, there is a bijection between $E_{v'} \setminus \{e\}$ and $E_{v''} \setminus \{e\}$ such that for any pair of corresponding edges $e' \in E_{v'} \setminus \{e\}$, $e'' \in E_{v''} \setminus \{e\}$, we have
    $$
        \sfp(e', v') - \sfp(e'', v'') = a \fr{(e,v')} \sfp(e,v') = -a \fr{(e,v'')} \sfp(e,v'')
    $$
    for some $a \in \bQ$.

    \item For each compact edge $e \in E(\hGa)_c$ that contains an $r$-valent vertex $v' \in V^r(\hGa)$ and a univalent vertex $v'' \in V^1(\hGa)$, we have
    $$
        \fr{(e,v')} \sfp(e,v') + \fr{(e,v'')} \sfp(e,v'') = 0.
    $$
    Moreover, there is a bijection between $E_{v'} \setminus \{e\}$ and $\{2, \dots, r\}$ such that for any edge $e' \in E_{v'} \setminus \{e\}$ corresponding to $i \in \{2, \dots, r\}$, we have
    $$
        \sfp(e', v') - \sff_i(v'') = a \fr{(e,v')} \sfp(e,v') = -a \fr{(e,v'')} \sfp(e,v'')
    $$
    for some $a \in \bQ$.
\end{itemize}
\end{definition}

Indeed, the conditions in Definition \ref{def:FTCYGraph} correspond to the conditions \eqref{eqn:WRelation}, \eqref{eqn:BaseConnectionCond}, \eqref{eqn:FiberConnectionCond} imposed on the local data in Sections \ref{sect:VertexData}, \ref{sect:EdgeData}.

\subsection{Relative FTCY orbifolds}
Given an FTCY graph $\sGa = (\hGa, \sfp, \sff_2, \dots, \sff_r)$, we now construct the associated FTCY orbifold $(\hcZ, \hcD)$ as follows. For each vertex $v \in V(\hGa)$, let
$$
    \hcZ_v
$$
be the formal toric orbifold constructed at the end of Section \ref{sect:VertexData}, and in the univalent case $v \in V^1(\hGa)$, let
$$
    \hcD_v
$$
be the divisor in $\hcZ_v$ constructed there. Moreover, for each compact edge $e \in E(\hGa)_c$ containing vertices $v_1, v_2 \in V(\hGa)$, let
$$
    \hcZ_e = \hcZ_{v_1} \cup \hcZ_{v_2}
$$
be the formal toric orbifold constructed at the end of Section \ref{sect:EdgeData}. We let
$$
    \hcZ := \bigcup_{v \in V(\hGa)} \hcZ_v
$$
where the charts $\{\hcZ_v\}_{v \in V(\hGa)}$ are glued together according to $\{\hcZ_e\}_{e \in E(\hGa)_c}$. Moreover, let
$$
    \hcD := \bigsqcup_{v \in V^1(\hGa)} \hcD_v.
$$
Then $\hcZ$ is an $r$-dimensional formal smooth Deligne-Mumford stack equipped with an action of the Calabi-Yau torus $T' \cong (\bC^*)^{r-1}$, and $\hcD$ is a divisor in $\hcZ$ with $|V^1(\hGa)|$ components such that it is invariant under the $T'$-action and that $(\hcZ, \hcD)$ is log Calabi-Yau.

\subsubsection{Curve classes}
Let $\hZ = \bigcup_{v \in V(\hGa)} \hZ_v$ be the coarse moduli space of $\hcZ$. We have
$$
	H_2(\hZ; \bZ) = \bigoplus_{e \in E(\hGa)_c} \bZ [l_e]
$$
where $l_e$ is the closed 1-dimensional torus orbit corresponding to the compact edge $e$. The subset of effective classes is
$$
	\Eff(\hZ) = \left\{\sum_{e \in E(\hGa)_c} d_e [l_e] : d_e  \ge 0 \text{ for all } e \right\}.
$$

\subsubsection{State spaces}
Moreover, we define
\begin{equation}\label{eqn:StateSpace}
    \cH(\hcZ) := \bigoplus_{v \in V(\hGa)} H_{\CR, T'}^*(\cZ_v; \bQ)    
\end{equation}
which may be viewed as the $T'$-equivariant Chen-Ruan cohomology of the FTCY orbifold $\hcZ$. We use the notation $\Box(v) := \Box(\cZ_v)$ for $v \in V(\hGa)$ and $\Box(\hcZ) = \bigcup_{v \in V(\hGa)} \Box(v)$. 

For each $v \in V^1(\vGa)$, since $\hcD_v$ is the formal completion of a toric Deligne-Mumford stack $\cD_v$, we put 
$$
    \cH(\hcD_v) := H_{\CR, T'}^*(\cD_v; \bQ).
$$
We use the notation $\Box(\cD_v)$.

\begin{example}\label{ex:GenuineFTCY} \rm{
Any toric Calabi-Yau $r$-orbifold $\cZ$ as in Section \ref{sect:CYCase} naturally gives rise to an FTCY graph with only $r$-valent vertices. The sets of vertices, edges, and flags in the graph correspond to $\Xi(r)$, $\Xi(r-1)$, and $F(\Xi)$ respectively. For a flag $(e,v)$ in the graph corresponding to $(\tau, \sigma) \in F(\Xi)$, the position vector $\sfp(e,v)$ is defined by the weight of the $T'$-action on $T_{\fp_\sigma}\fl_\tau$. The resulting FTCY orbifold $\hcZ$ is the formal completion of $\cZ$ along $\cZ^1$. There is a surjective map
$$
	\pi: H_2(\hZ; \bZ) \to H_2(Z; \bZ), \qquad [l_e] \mapsto [l_e]
$$
which preserves the sets of effective classes. Moreover, since $\cZ$ is equivariantly formal over $\bQ$, the natural injective map
$$
    H_{\CR, T'}^*(\cZ; \bQ) \to \cH(\hcZ)
$$
defined by restricting to affine charts is an isomorphism after base change to the field of fractions of $H_{T'}^*(\pt; \bQ)$ (cf. \cite{GKM98,LS20}).
}\end{example}



\section{Gromov-Witten theory of relative formal toric Calabi-Yau orbifolds}\label{sec:GW-FTCY}

In this section, we study the equivariant Gromov-Witten invariants of relative FTCY orbifolds. We refer to \cite{AGV02, AGV08, CR02} for foundations of orbifold Gromov-Witten theory and \cite{Liu13,LS20} for the cases of toric and GKM orbifolds. We refer to \cite{AF16} for relative orbifold Gromov-Witten theory.

\subsection{Relative stable maps to relative FTCY orbifolds}\label{sect:RelativeStableMaps}
Let $\sGa = (\hGa, \sfp, \sff_2, \dots, \sff_r)$ be an FTCY graph and $(\hcZ,\hcD)$ be its associated $r$-dimensional relative FTCY orbifold with the action of the Calabi-Yau torus $T' \cong (\bC^*)^{r-1}$, as in Section \ref{sect:FTCY}. Let $\hZ$ be the coarse moduli space of $\hcZ$. We write $V^1(\hGa) = \{v_1, \dots, v_s\}$ and $\hcD = \hcD_1 + \cdots + \hcD_s$. For each $i$, let $\fp_i$ denote the unique $T'$-fixed point of $\hcD_i$ and $G_i$ denote the group of generic stabilizers of $\fp_i$, which is in bijection with $\Box(\hcD_i)$.

\subsubsection{Target}
The target of a relative stable map can be an \emph{expanded pair}. To the univalent vertex $v_i$, recall from Section \ref{sect:VertexData} that we associated a pair $(\cZ_i, \cD_i)$ such that $\hcD_i$ is the formal completion of $\cD_i$ along its toric 1-skeleton. Let $\Delta_i$ be the total space of $\bP(\cO_{\cD_i} \oplus N_{\cD_i/\cZ_i})$. The action of $T'$ on $\cD_i$ extends to an action on $\Delta_i$ under which the fixed locus is the fiber over $\fp_i$. Now given $m_i \in \bZ_{\ge 0}$, let
$$
	\cZ_i(m_i) := \Delta_{i, 1} \cup \cdots \cup \Delta_{i, m_i}
$$
where each $\Delta_{i, j}$ is a copy of $\Delta_i$ with distinguished sections $\cD_{i, j-1} = \bP(\cO_{\cD_i} \oplus 0)$ and $\cD_{i, j} = \bP(0 \oplus N_{\cD_i/\cZ_i})$, and $\Delta_{i, j}, \Delta_{i, j+1}$ are glued along $\cD_{i, j}$. There is a natural projection
$$
	\pi_{i, m_i}: \cZ_i(m_i) \to \cD_{i, 0} = \cD_i.
$$
For each $j$, let $\fp_{i, j} := \pi_{i, m_i}^{-1}(\fp_i) \cap \cD_{i,j} \cong \fp_i$ and $\fl_{i,j} := \pi_{i, m_i}^{-1}(\fp_i) \cap \Delta_{i,j} \cong \fp_i \times \bP^1$. The torus $T'$ acts on $\cZ_i(m_i)$ with fixed locus
$$
	\fl_{i, 1} \cup \cdots \cup \fl_{i, m_i}.
$$
Moreover, there is an action of $(\bC^*)^{m_i}$ on $\cZ_i(m_i)$ which scales the fiber directions of the respective components $\Delta_{i,j}$. The projection $\pi_{i, m_i}$ is equivariant with respect to both torus actions. We let $\hcZ_i(m_i)$ (resp. $\hcD_{i,j}$) denote the formal completion of $\cZ_i(m_i)$ (resp. $\cD_{i, j}$) along the toric 1-skeleton. Then given $\bm = (m_1, \dots, m_s) \in \bZ_{\ge 0}^s$, we define the expanded pair $(\hcZ_{\bm},\hcD_{\bm})$ where
$$
	\hcZ_{\bm} := \hcZ \cup \hcZ_1(m_1) \cup \cdots \cup \hcZ_s(m_s), \qquad \hcD_{\bm} := \hcD_{1, m_1} \sqcup \cdots \sqcup \hcD_{s, m_s}
$$
and each $\hcZ_i(m_i)$ is glued to $\hcZ$ along $\hcD_{i, 0} = \hcD_i$. For $\bm = (0, \dots, 0)$, we have $(\hcZ_{\bm},\hcD_{\bm}) = (\hcZ,\hcD)$. We have a projection
$$
	\pi_{\bm}: (\hcZ_{\bm},\hcD_{\bm}) \to (\hcZ,\hcD).
$$

\subsubsection{Domain}
The domain of an orbifold relative stable map is a prestable \emph{twisted curve} with $n \in \bZ_{\ge 0}$ marked points, which is a 1-dimensional proper connected Deligne-Mumford stack $\cC$, with at most nodal singularities, together with disjoint closed substacks $\fx_1, \dots, \fx_n$ of the smooth locus $\cC^\sm$ of $\cC$, such that $\cC^\sm \setminus \bigcup_{i = 1}^n \fx_i$ is a scheme and that each node is a balanced node. The latter condition means that locally near a node, $\cC$ is isomorphic to
$$
[\Spec(\bC[x,y]/(xy))/\mu_r]
$$
for some $r \in \bZ_{\ge 1}$ where the action of $\mu_r$ is balanced in the sense that $\zeta \in \mu_r$ acts by $\zeta \cdot (x,y) = (\zeta x, \zeta^{-1}y)$. Similarly, locally near $\fx_i$, $\cC$ is isomorphic to
$$
	[\Spec(\bC[x])/\mu_{r_i}]
$$
for some $r_i \in \bZ_{>0}$ where $\zeta \in \mu_{r_i}$ acts by $\zeta \cdot x = \zeta x$. The integer $r_i$ is the orbifold index of $\fx_i$. Let $C$ denote the coarse moduli space of $\cC$ and $x_i$ denote the image of $\fx_i$. Then $(C, x_1, \dots, x_n)$ is an $n$-pointed prestable curve. We focus on the genus-zero case, where $C$ has arithmetic genus zero.

\subsubsection{Curve class and tangency profile}
For each $i$, let $e_i \in E(\hGa)_c$ be the unique compact edge incident to the univalent vertex $v_i$. Let 
$$
    \hbeta = \sum_{e \in E(\hGa)_c} d_e [l_e] \qquad \in \Eff(\hZ)
$$
be an effective class such that $d_i := d_{e_i} \ge 1$ for all $i = 1, \dots, s$. For degree-$\hbeta$ relative stable maps to $(\hcZ, \hcD)$, the tangency profile with the component $\hcD_i$ is in general specified by a partition of $d_i$ whose parts are weighted by elements in the stabilizer group $G_i$. We focus on the case of \emph{maximal tangency} where the partition has a single part $d_i$, weighted by $k_i \in G_i$. Let $\ord(k_i)$ denote the order of $k_i$ in $G_i$. Let $\age_{\hcD}(k_i)$ denote the age of $k_i$ viewed as an element of $\Box(\hcD_i)$, which is an integer. We write 
$$
    \bk = (k_1, \dots, k_s).
$$

\subsubsection{Map}
A genus-zero, $n$-pointed, degree-$(\hbeta, \bk)$ maximally tangent \emph{twisted relative stable map} to $(\hcZ, \hcD)$ is a representable morphism
$$
	u: (\cC, \underline{\fx}, \underline{\fy}) = (\cC, \fx_1, \dots, \fx_n, \fy_{n+1}, \dots, \fy_{n+s}) \to (\hcZ_{\bm}, \hcD_{\bm}).
$$
where:
\begin{itemize}
	\item $\bm \in \bZ_{\ge 0}^s$ and $(\hcZ_{\bm}, \hcD_{\bm})$ is the expanded pair.
	
	\item $(\cC, \underline{\fx}, \underline{\fy})$ is a genus-zero prestable twisted curve with $n+s$ marked points.
	
	\item On the level of coarse moduli, $(\pi_{\bm} \circ u)_*[C] = \hbeta$.
	
	\item For each $i = 1, \dots, s$, $u(\fy_{n+i}) \subset \hcD_{i, m_i}$ and $u^{-1}(\hcD_{i, m_i}) = d_i\ord(k_i) \fy_{n+i}$ as Cartier divisors. Here, $d_i$ is the intersection multiplicity of the map with $\hcD_{i, m_i}$ and $d_i\ord(k_i)$ is the contact order. The orbifold index of $\fy_{n+i}$ is $\ord(k_i)$ and $k_i$ specifies the monodromy around $\fy_{n+i}$.
	
	\item For each $i = 1, \dots, s$, $j = 1, \dots, m_i-1$, the preimage $u^{-1}(\hcD_{i,j})$ consists of nodes $\fq$ in $\cC$, and the on the two irreducible components of $\cC$ that intersect at $\fq$, $u$ has the same contact order with $\hcD_{i,j}$.
	
	\item $u$ has finite automorphism group. Here, an automorphism is a pair $(\alpha_1, \alpha_2)$ where $\alpha_1$ is an automorphism of $(\cC, \underline{\fx}, \underline{\fy})$ and $\alpha_2$ is an automorphism of $(\hcZ_{\bm}, \hcD_{\bm})$ invariant under $\pi_{\bm}$, such that $u \circ \alpha_1 = \alpha_2 \circ u$.
	
\end{itemize}

\subsection{Moduli space and torus-fixed locus}\label{sect:Moduli}
Let
$$
	\Mbar_{0,n}(\hcZ/\hcD, \hbeta, \bk)
$$
be the moduli space of genus-zero, $n$-pointed, degree-$(\hbeta, \bk)$ maximally tangent twisted relative stable maps to $(\hcZ, \hcD)$. It admits a perfect obstruction theory. There are evaluation maps
$$
	\ev_i: \Mbar_{0,n}(\hcZ/\hcD, \hbeta, \bk) \to \cI \hcZ, \qquad i = 1, \dots, n,
$$
$$
	\ev_{n+i}: \Mbar_{0,n}(\hcZ/\hcD, \hbeta, \bk) \to \cI \hcD_i, \qquad i = 1, \dots, s
$$
at the marked points $\fx_i, \fy_{n+i}$ respectively. The image of $\ev_{n+i}$ lies in the inertia component $(\hcD_i)_{k_i}$ of $\cI \hcD_i$ indexed by $k_i$. For $\bj = (j_1, \dots, j_n) \in \Box(\hcZ)^n$, let
$$
	\Mbar_{0,\bj}(\hcZ/\hcD, \hbeta, \bk) := \bigcap_{i = 1}^n \ev_i^{-1}(\hcZ_j),
$$
which has virtual dimension
$$
	r - 3 + n + s - \sum_{i = 1}^n \age(j_i) - \sum_{i = 1}^s \age_{\hcD}(k_i).
$$
We have
$$
	\Mbar_{0,n}(\hcZ/\hcD, \hbeta, \bk) = \bigsqcup_{\bj \in \Box(\hcZ)^n} \Mbar_{0,\bj}(\hcZ/\hcD, \hbeta, \bk).
$$

The action of $T'$ on $(\hcZ,\hcD)$ induces an action on $\Mbar_{0,n}(\hcZ/\hcD, \hbeta, \bk)$ under which the evaluation maps are $T'$-equivariant, the perfect obstruction theory is $T'$-equivariant, and the $T'$-fixed locus is proper. Components of the $T'$-fixed locus can be described by \emph{decorated graphs}.

\begin{definition}\rm{
A genus-zero, $(\bj, \bk)$-twisted, degree-$\hbeta$ \emph{decorated graph} for $(\hcZ,\hcD)$ is a tuple $\vGa = (\Gamma, \vf, \vd,\vs,\vk)$, where:
\begin{itemize}
	\item $\Gamma = (V(\Gamma), E(\Gamma))$ is a compact, connected, $1$-dimensional CW complex, where $V(\Gamma)$ is the vertex set and $E(\Gamma)$ is the edge set. Let $F(\Gamma)$ denote the set of flags, consisting of pairs $(e, v) \in E(\Gamma) \times V(\Gamma)$ with $v \in e$.
	
	\item $\vf: V(\Gamma) \sqcup E(\Gamma) \to V(\hGa) \sqcup E(\hGa)_c$ is the \emph{label map} that labels each $v \in V(\Gamma)$ by a vertex $\sigma_v \in V(\hGa)$ and each $e \in E(\Gamma)$ by a compact edge $\tau_e \in E(\hGa)_c$ such that for each flag $(e,v) \in F(\Gamma)$, $(\tau_e, \sigma_v)$ is a flag in $F(\hGa)$. And it induces the label map of flags $\vf: F(\Gamma)\to F(\vGa)$. We denote $G_v:= G_{\sigma_v}$ for $v \in V(\Gamma)$ and $G_e:= G_{\tau_e}$ for $e \in E(\Gamma)$. We partition $V(\Gamma)$ into two subsets
	$$
		V(\vGa)^{(0)} := \vf^{-1}(V^r(\hGa)), \qquad V(\vGa)^{(1)} := \vf^{-1}(V^1(\hGa)).
	$$
	
	\item $\vd$ is the \emph{degree map} that assigns to each edge $e \in E(\Gamma)$ an element $\gamma_e \in H_{\tau_e}$ such that $d_e:= \pi_{\tau_e}(\gamma_e)$ is a positive integer.
	
	\item $\vs: \{1, \dots, n, \dots, n+s\} \to V(\Gamma)$ is the \emph{marking map}.
	
	\item $\vk$ is the \emph{twisting map} that sends each flag $(e,v) \in F(\Gamma)$ to an element $k_{(e,v)} \in G_v$ and each $i = 1, \dots, n, \dots, n+s$ to an element $\vk(i) \in G_{\vs(i)}$.
\end{itemize}
The above data is required to satisfy the following conditions:
\begin{itemize}
		\item The underlying graph $\Gamma = (V(\Gamma), E(\Gamma))$ is a tree, i.e.
		$$
		|E(\Gamma)|-|V(\Gamma)|+1 = 0.
		$$
	
	\item We have $\displaystyle{   \sum_{e \in E(\Gamma)} d_e[l_{\tau_e}] = \hbeta. }$
	
	\item For any edge $e \in E(\Gamma)$ incident to vertices $v, v' \in V(\Gamma)$, we have
	$$
		\pi_{(\tau_e, \sigma_v)}(\gamma_e) = k_{(e,v)}, \quad \pi_{(\tau_e, \sigma_{v'})}(\gamma_e) = k_{(e,v')}.
	$$
	
	\item For any vertex $v \in V(\Gamma)$, we have
	$$
	\prod_{(e,v) \in F(\Gamma)} k_{(e,v)}^{-1} \prod_{i \in \vs^{-1}(v)}\vk(i) = 1
	$$
	in $G_v$.
	
	\item For $i = 1, \dots, n$, the pair $(\fp_{\sigma_{\vs(i)}}, \vk(i))$ represents a point in the inertia component $\hcZ_{j_i}$.
	
	\item For $i = 1, \dots, s$, we have $\vf \circ \vs(n+i) = v_i \in V^1(\hGa)$ and $\vk(n+i) = k_i$. 
\end{itemize}
}\end{definition}

Let $\Gamma_{0, \bj}(\hcZ/\hcD, \hbeta, \bk)$ be the set of all decorated graphs as above. We set up the following additional notation on a decorated graph $\vGa$:
\begin{itemize}
	\item For each $v \in V(\Gamma)$, let
	$$
	E_v := \{e \in E(\Gamma) : (e,v) \in F(\Gamma) \}, \quad \val(v):= |E_v|,  \quad S_v:= \vs^{-1}(v), \quad n_v := |S_v|,
	$$
	and $\vk_v:= (k_{(e,v)}^{-1}, \vk(i)) \in G_v^{E_v \cup S_v}$.
	
	\item Let
	$$
	V^S(\vGa):= \{v \in V(\Gamma): \val(v) + n_v -2 >0 \},
	$$
	which is the set of \emph{stable} vertices of $\Gamma$, and
	\begin{align*}
		V^1(\vGa) &:= \{v \in V(\Gamma): \val(v) = 1, n_v = 0 \},\\
		V^{1,1}(\vGa) &:= \{v \in V(\Gamma): \val(v) = n_v = 1 \},\\
		V^2(\vGa) &:= \{v \in V(\Gamma): \val(v) = 2, n_v = 0 \},
	\end{align*}
	which partition the set of \emph{unstable} vertices $V(\Gamma) \setminus V^S(\vGa)$.
	
	\item For $(e,v) \in F(\Gamma)$, let $r_{(e,v)}$ be the order of $k_{(e,v)}$ in $G_v$. For $v \in V^2(\vGa)$ with $E_v = \{e_1, e_2\}$, let $r_v := r_{(e_1,v)} = r_{(e_2, v)}$.
	
	\item For $i = 1, \dots, s$, let $v_i(\vGa)$ be the unique vertex in $V(\vGa)^{(1)}$ that is labeled by $v_i$. In particular, $V(\vGa)^{(1)} = \{v_1(\vGa), \dots, v_s(\vGa)\}$. We have $S_{v_i(\vGa)} = \{n+i\}$ and $n_{v_i(\vGa)} = 1$.

	\item For $i = 1, \dots, s$, let $\mu_i(\vGa)$ be the partition of $d_i$ determined by the degrees of edges in $E_{v_i(\vGa)} = \vf^{-1}(e_i)$. We have $\val(v_i(\vGa)) = l(\mu_i(\vGa))$, the length of $\mu_i(\vGa)$.
	
	\item Let $\Aut(\vGa)$ be the \emph{automorphism group} of $\vGa$, which consists of all automorphisms of $\Gamma$ under which $\vf, \vd, \vs, \vk$ are invariant. Let $A(\vGa)$ denote the subgroup of $\Aut(\vGa)$ of automorphisms which fixes any edge labeled by some $e_i$. Define the coefficient
	\begin{equation}\label{eqn:cGamma}
		c_{\vGa}:= \frac{1}{|A(\vGa)| \cdot \prod_{e \in E(\Gamma)}(d_e|G_e|)} \cdot \prod_{(e,v) \in F(\Gamma)} \frac{|G_v|}{r_{(e,v)}}.
	\end{equation}

\end{itemize}

Given a stable map $u: (\cC, \underline{\fx}, \underline{\fy}) \to (\hcZ_{\bm}, \hcD_{\bm})$ that represents a point in the $T'$-fixed locus $\Mbar_{0,\bj}(\hcZ/\hcD, \hbeta, \bk)^{T'}$, we can assign a decorated graph $\vGa \in \Gamma_{0, \bj}(\hcZ/\hcD, \hbeta, \bk)$ as follows. Let $\tu = \pi_{\bm} \circ u: \cC \to \hcZ$ denote the composition and $\bar{u}: C \to \hZ$ denote the induced map on the level of coarse moduli. 
\begin{itemize}
	\item The vertex set $V(\Gamma)$ is in bijection with the set of connected components of $\bar{u}^{-1}(\hZ^{T'})$. For $v \in V(\Gamma)$, let $C_v \subseteq C$ denote the component corresponding to $v$ and $\cC_v \subseteq \cC$ denote the preimage. Let $\vf(v) = \sigma_v \in V(\hGa)$ such that $\cC_v$ is mapped to $\fp_{\sigma_v}$ under $\tu$.
	
	\item The edge set $E(\Gamma)$ is in bijection with the set of irreducible components of $C$ that do not map constantly to $\hat{Z}$ under $\bar{u}$. For $e \in E(\Gamma)$, let $C_e \subseteq C$ denote the component corresponding to $e$ and $\cC_e \subseteq \cC$ denote the preimage. Let $\vf(e) = \tau_e \in E(\hGa)_c$ such that $\cC_e$ is mapped to $\fl_{\tau_e}$ under $\tu$.
	
	\item The flag set $F(\Gamma)$ consists of $(e,v)$ such that $\cC_e \cap \cC_v \neq \emptyset$. For $(e,v) \in F(\Gamma)$, let $\fn(e,v) := \cC_e \cap \cC_v$ which is a node or marked point of $\cC$. Let $k_{(e,v)} \in G_v$ be the image of the generator of the generic stabilizer group of $\fn(e,v)$ in $\cC_e$ under $\tu$.
	
	\item Let $e \in E(\Gamma)$ and $v, v' \in V(\Gamma)$ be incident vertices. We have $\cC_e \cong \cC_{r_{(e,v)}, r_{(e,v')}}$. Let $\gamma_e \in H_{\tau_e}$ be the element defined by $\tu|_{\cC_e}$ and let $d_e= \pi_{\tau_e}(\gamma_e)$. The compatibility conditions $\pi_{(\tau_e, \sigma_v)}(\gamma_e) = k_{(e,v)}$ and $\pi_{(\tau_e, \sigma_v')}(\gamma_e) = k_{(e,v')}$ are satisfied.
	
	\item For $i = 1, \dots, n$, let $\vs(i) = v \in V(\Gamma)$ such that $\fx_i \subseteq \cC_v$. Let $\vk(i) \in G_v$ such that $\fx_i$ is mapped to $(\fp_{\sigma_v}, \vk(i))$ under $\tu$. Then for each $v \in V^S(\vGa) \cap V(\vGa)^{(0)}$, $u \big|_{\cC_v}$ represents a point in $\Mbar_{0, \vk_v}(\cB G_v)$.
\end{itemize}
Under this assignment, we have that for $i = 1, \dots, s$, $l(\mu_i(\vGa)) = 1$ if and only if $m_i = 0$. If $m_i = 0$, we have $\cC_{v_i(\vGa)} = \fy_{n+i}$ and $v_i(\vGa) \in V^{1,1}(\vGa)$. If $m_i > 0$, the restriction of $u$ to the component $\cC_{v_i(\vGa)}$ represents a point in the moduli space
$$
	\Mbar_{0,0}(\bP^1 \times G_i, \mu_i(\vGa), (d_i)) \sslash \bC^*,
$$ 
of relative stable maps to the non-rigid $(\bP^1, 0, \infty) \times G_i$ with relative condition $\mu_i(\vGa)$ at $0$ and $(d_i)$ at $\infty$, where $(d_i)$ is weighted by $k_i \in G_i$ and $\mu_i(\vGa)$ is weighted by the elements $k_{(e, v_i(\vGa))} \in G_i$ as $e$ ranges through $E_{v_i(\vGa)}$. We refer to \cite{Zong12,Zong15} for a detailed description of this moduli space when $G_i$ is \emph{cyclic}, which is the case needed in the rest of this paper.

The above assignment gives a decomposition of the $T'$-fixed locus
$$
	\Mbar_{0,\bj}(\hcZ/\hcD, \hbeta, \bk)^{T'} = \bigsqcup_{\vGa \in \Gamma_{0,\bj}(\hcZ/\hcD, \hbeta, \bk)} \cF_{\vGa}
$$
into connected components $\cF_{\vGa}$, where up to a finite morphism, $\cF_{\vGa}$ is be identified with
$$
	\cM_{\vGa} := \prod_{v \in V^S(\vGa) \cap V(\vGa)^{(0)}} \Mbar_{0, \vk_v}(\cB G_v) \times \prod_{m_i > 0} \Mbar_{0,0}(\bP^1 \times G_i, \mu_i(\vGa), (d_i)) \sslash \bC^*,
$$
and in $A_*(\cM_{\vGa})$ we have
\begin{equation}\label{eqn:FGammaFiniteMorphism}
	[\cF_{\vGa}] = \frac{1}{|A(\vGa)| \cdot \prod_{e \in E(\Gamma)}(d_e|G_e|)} \cdot \prod_{v \in V^S(\vGa), e \in E_v} \frac{|G_v|}{r_{(e,v)}} \cdot \prod_{v \in V^2(\vGa)} \frac{|G_v|}{r_v} \cdot [\cM_{\vGa}].
\end{equation}

Finally, we set
$$
	\Gamma_{0,n}(\hcZ/\hcD, \hbeta, \bk) := \bigsqcup_{\bj \in \mathrm{Box}(\hcZ)^n} \Gamma_{0,\bj}(\hcZ/\hcD, \hbeta, \bk).
$$
We have
$$
	\Mbar_{0,n}(\hcZ/\hcD, \hbeta, \bk)^{T'} = \bigsqcup_{\vGa \in \Gamma_{0,n}(\hcZ/\hcD, \hbeta, \bk)} \cF_{\vGa}.
$$

\subsection{Formal relative Gromov-Witten invariants}\label{sect:RGWDef}
Gromov-Witten invariants of the relative FTCY orbifold $(\hcZ, \hcD)$ are virtual counts of the stable maps and may be defined by virtual localization \cite{GP99, GV05}. As insertions, let $\gamma_1, \dots, \gamma_n \in \cH(\hcZ)$, and for $i = 1, \dots, s$, let $\gamma_{n+i} \in H^*_{T'}((\cD_i)_{k_i}; \bQ) \subseteq \cH(\hcD_i)$ where $(\cD_i)_{k_i}$ denotes the inertia component of $\cD_i$ indexed by $k_i$. We define the formal relative Gromov-Witten invariant
\begin{equation}\label{eqn:RGWDef}
    \langle \gamma_1,\dots,\gamma_{n} \mid \gamma_{n+1},\dots,\gamma_{n+s} \rangle^{\hcZ/\hcD, T'}_{\hbeta,\bk} := \int_{[\Mbar_{0,n}(\hcZ/\hcD, \hbeta, \bk)^{T'}]^\vir} \frac{\iota^*\left(\prod_{i = 1}^n \ev_i^*(\gamma_i) \prod_{i = 1}^s \ev_{n+i}^*(\gamma_{n+i})\right)}{e_{T'}(N^\vir)} 
\end{equation}
where $\iota: \Mbar_{0,n}(\hcZ/\hcD, \hbeta, \bk)^{T'} \to \Mbar_{0,n}(\hcZ/\hcD, \hbeta, \bk)$ is the inclusion of the $T'$-fixed locus and $N^\vir$ is the virtual normal bundle. The invariant takes value in the field of fractions of $H^*_{T'}(\pt; \bQ)$.

\subsection{Specialization to the absolute case}\label{sect:Absolute}
Now we specialize to the case where the set $V^1(\hGa)$ of univalent vertices in the FTCY graph $\sGa$ is empty, in which case $\hcD = \emptyset$. The discussion of Sections \ref{sect:RelativeStableMaps}, \ref{sect:Moduli} concerns the moduli space $\Mbar_{0,n}(\hcZ, \hbeta)$ of usual genus-zero, $n$-pointed, degree-$\hbeta$ twisted stable maps
$$
	u: (\cC, \underline{\fx}) = (\cC, \fx_1, \dots, \fx_n) \to \hcZ
$$
to the FTCY orbifold $\hcZ$, which is a special case of formal GKM orbifolds treated in \cite{LS20}. 
Connected components of the $T'$-fixed locus of the moduli is indexed by the set
$$
	\Gamma_{0,n}(\hcZ, \hbeta) = \bigsqcup_{\bj \in \mathrm{Box}(\hcZ)^n} \Gamma_{0,\bj}(\hcZ, \hbeta)
$$
and we have a decomposition
$$
	\Mbar_{0,n}(\hcZ, \hbeta)^{T'} = \bigsqcup_{\vGa \in \Gamma_{0,n}(\hcZ, \hbeta)} \cF_{\vGa}.
$$
For each decorated graph $\vGa \in \Gamma_{0,n}(\hcZ, \hbeta)$, we have $V(\vGa)^{(1)} = \emptyset$ and $A(\vGa) = \Aut(\vGa)$. Equation \eqref{eqn:RGWDef} defines the usual formal Gromov-Witten invariant $\langle \gamma_1,\dots,\gamma_{n} \rangle^{\hcZ, T'}_{\hbeta}$.

Suppose furthermore that the FTCY graph $\sGa$ is induced by a toric Calabi-Yau $r$-orbifold $\cZ$, as in Example \ref{ex:GenuineFTCY}. For an effective curve class $\beta \in \Eff(Z)$, the $T'$-fixed locus of the moduli space $\Mbar_{0,n}(\cZ, \beta)$ of genus-zero, $n$-pointed, degree-$\beta$ twisted stable maps to $\cZ$ has a decomposition
$$
	\Mbar_{0,n}(\cZ, \beta)^{T'} = \bigsqcup_{\hbeta \in \Eff(\hZ), \pi(\hbeta) = \beta} \Mbar_{0,n}(\hcZ, \hbeta)^{T'}
$$
where $\pi: \Eff(\hZ) \to \Eff(Z)$ is the projection. Therefore, the connected components are indexed by the set
$$
	\Gamma_{0,n}(\cZ, \beta) :=  \bigsqcup_{\hbeta \in \Eff(\hZ), \pi(\hbeta) = \beta} \Gamma_{0,n}(\hcZ, \hbeta)
$$
of decorated graphs. The decomposition also respects the refinement by elements $\bj \in \Box(\cZ) = \Box(\hcZ)$ and we use analogous notation for $\cZ$. Moreover, for $\gamma_1, \dots, \gamma_n \in H^*_{\CR, T'}(\cZ; \bQ) \subseteq \cH(\hcZ)$, we have a decomposition 
$$
    \langle \gamma_1,\dots,\gamma_{n} \rangle^{\cZ, T'}_{\beta} = \sum_{\hbeta \in \Eff(\hZ), \pi(\hbeta) = \beta} \langle \gamma_1,\dots,\gamma_{n} \rangle^{\hcZ, T'}_{\hbeta}
$$
of the usual equivariant Gromov-Witten invariant of $\cZ$ as a sum of formal invariants. For later use, in this case, we introduce the notation
$$
	\bpsi_i := \epsilon^* \psi_i \in A^1(\Mbar_{0,n}(\cZ, \beta)), \quad i = 1, \dots, s
$$
for the descendant class where $\epsilon: \Mbar_{0,n}(\cZ, \beta) \to \Mbar_{0,n}(Z, \beta)$ is the natural map.


\section{Multi-component geometries}\label{sect:Geometry}
In this section, we describe the various geometries involved in the correspondences. The construction consists of the following three aspects:
\begin{itemize}
    \item Open geometry $(\cX, \cL, f)$ (Sections \ref{sect:TCY3}-\ref{sect:OpenGeometry}), where $\cX$ is a toric Calabi-Yau 3-orbifold and $\cL \subset \cX$ is an $s$-component Lagrangian suborbifold of Aganagic-Vafa type, equipped with the data of a framing $f$.

    \item Closed/local geometry $\tcX$ (Section \ref{sect:LocalGeometry}), which is a toric Calabi-Yau ($3+s$)-orbifold.
    
    \item Intermediate relative geometries $(\hcY^{(\ell)}, \hcD^{(\ell)})$ for $\ell = 0, \dots, s$ (Sections \ref{sect:RelativeGeometry}, \ref{sect:IntermediateGeometry}), which is a relative FTCY ($3+\ell$)-orbifold.

\end{itemize}
We describe the example $\cX = \bC^3$ explicitly in Section \ref{sect:example}. 


\subsection{Toric Calabi-Yau 3-orbifolds}\label{sect:TCY3}
Let $\cX$ be a toric Calabi-Yau 3-orbifold defined by an extended stacky fan $\bSi = (N, \Sigma, \alpha)$. We specialize the notation in Sections \ref{sect:ToricOrbifolds}, \ref{sect:CYCase} to this case as follows:
\begin{itemize}
    \item $X$ denotes the coarse moduli space of $\cX$.
    
    \item $N \cong \bZ^3$, $T:= N \otimes \bC^* \cong (\bC^*)^3$, $M := \Hom(N, \bZ)$.
    
    \item $\alpha: \bZ^R \to N$, $e_i \mapsto b_i \in N$, where $\{e_1, \dots, e_R\}$ is the standard basis of $\bZ^R$.

    \item $R' := |\Sigma(1)|$, $\Sigma(1) = \{\bR_{\ge 0}b_1, \dots, \bR_{\ge 0}b_{R'}\}$.

    \item Let $\su_3 \in M$ be the character determined by the Calabi-Yau condition, and $N' := \ker(\su_3) \cong \bZ^2$, $T' := N' \otimes \bC^* \cong (\bC^*)^2$, $M' := \Hom(N', \bZ)$.
    
\end{itemize}



In the 3-dimensional case, given a flag $(\tau, \sigma) \in F(\Sigma)$, $G_\tau$ is a cyclic subgroup of $G_\sigma$. We define
$$
    \fm(\tau, \sigma):= |G_\tau|.
$$
Moreover, we define
$$
    \sw(\tau, \sigma) := c_1^{T'}(T_{\fp_\sigma}\fl_\tau) \in H^2_{T'}(\pt;\bQ) \cong M'_{\bQ}.
$$
The Calabi-Yau condition implies that for any $\sigma \in \Sigma(3)$, we have
$$
    \sum_{(\tau, \sigma) \in F(\Sigma)} \sw(\tau, \sigma) = 0.
$$



Let $P$ be the cross section of the support $|\Sigma|$ of $\Sigma$ in the hyperplane $N'_{\bR} \times \{1\}$. Then $|\Sigma|$ is the cone over $P$, and $\Sigma$ induces a triangulation of $P$. We assume that all maximal cones in $\Sigma$ are 3-dimensional and that $P$ is a simple polygon. Moreover, we assume that the additional lattice points $b_{R'+1}, \dots, b_R$ are chosen in a way that $(b_1, \dots, b_R)$ is a listing of the points in $P \cap N$. In particular, the homomorphism $\alpha$ is surjective and fits into the following short exact sequence of lattices:
\begin{equation}\label{eqn:XFanSequence}
    \xymatrix{
        0 \ar[r] & \bL \ar[r]^{\psi} & \bZ^R \ar[r]^{\alpha} & N \ar[r] & 0
    }
\end{equation}
where $\bL := \ker(\alpha) \cong \bZ^{R-3}$.

Let $G:= \bL \otimes \bC^* \cong (\bC^*)^{R-3}$. Let $\{\epsilon_1, \dots, \epsilon_{R-3}\}$ be a basis for $\bL$, and for each $a = 1, \dots, R-3$, let
$$
    l^{(a)} = (l^{(a)}_1, \dots, l^{(a)}_R) := \psi(\epsilon_a) \in \bZ^R.
$$
The vectors $l^{(a)}$ are known as \emph{charge vectors}, which describe the linear action of $G$ on $\bC^R = \bZ^R \otimes \bC = \Spec(\bC[x_1, \dots, x_R])$ induced by the inclusion $\psi$, as follows:
\begin{equation}\label{eqn:GAction}
    (s_1, \dots, s_{R-3}) \cdot (x_1, \dots, x_R) = \left(\prod_{a=1}^{R-3}s_a^{l_1^{(a)}}x_1, \dots, \prod_{a=1}^{R-3}s_a^{l_R^{(a)}}x_R  \right),
\end{equation}
where $(s_1, \dots, s_{R-3})$ are coordinates on $G$ specified by the basis $\{\epsilon_1, \dots, \epsilon_{R-3}\}$. Under this action, $\cX$ can be presented as the quotient stack
\begin{equation}\label{eqn:XQuotient}
    \cX = \left[ ((\bC^{R'} \setminus Z(\Sigma)) \times (\bC^*)^{R-R'})/G \right],
\end{equation}
where $Z(\Sigma)$ is a closed subvariety of $\bC^{R'}$ defined by $\Sigma$.

In this paper, we do not assume that $\cX$ has semi-projective coarse moduli space $X$ or equivalently $P$ is convex. In general, $\cX$ can be embedded as an open substack of a toric Calabi-Yau 3-orbifold $\bar{\cX}$ with semi-projective coarse moduli space $\bar{X}$. The fan $\bar{\Sigma}$ of $\bar{X}$ contains $\Sigma$ as a subfan and satisfies that $\bar{\Sigma}(1) = \Sigma(1)$ and $|\bar{\Sigma}|$ is the cone over the convex hull of $P$.

\subsection{Aganagic-Vafa branes}\label{sect:AVBranes}\
In this section, we describe the symplectic structure on $\cX$ and define Aganagic-Vafa Lagrangian branes. We refer to \cite{FLT12} for additional details. 

\subsubsection{Semi-projective partial compactification}
For the purpose above, we fix an embedding of $\cX$ as an open substack of a toric Calabi-Yau 3-orbifold $\bar{\cX}$ with semi-projective coarse moduli space. We may apply the discussion in Section \ref{sect:TCY3} to describe the geometry of $\bar{\cX}$. The extended stacky fan of $\bar{\cX}$ has form $\bar{\bSi} = (N, \bar{\Sigma}, \bar{\alpha})$ where $\bar{\alpha}: \bZ^{\bar{R}} \to N$ is a surjective homomorphism. Here, $\bar{R}$ is the number of lattice points in $N$ that is contained in the convex hull of $P$, and may be greater than $R$ in general. Similar to \eqref{eqn:XFanSequence}, there is a short exact sequence
\begin{equation}\label{eqn:XspFanSequence}
    \xymatrix{
        0 \ar[r] & \bar{\bL} \ar[r]^{\bar{\psi}} & \bZ^{\bar{R}} \ar[r]^{\bar{\alpha}} & N \ar[r] & 0
    }
\end{equation}
where $\bar{\bL} := \ker(\bar{\alpha}) \cong \bZ^{\bar{R}-3}$. This induces a linear action of $\bar{G} := \bar{\bL} \otimes \bC^* \cong (\bC^*)^{\bar{R}-3}$ on $\bC^{\bar{R}}$ that can be described by charge vectors
$$
    \bar{l}^{(a)} = (\bar{l}^{(a)}_1, \dots, \bar{l}^{(a)}_{\bar{R}}) \in \bZ^{\bar{R}}, \qquad a = 1, \dots, \bar{R}-3.
$$
Similar to \eqref{eqn:XQuotient}, $\bar{\cX}$ can be presented as a quotient stack
$$    
    \bar{\cX} = \left[ ((\bC^{R'} \setminus Z(\bar{\Sigma})) \times (\bC^*)^{\bar{R}-R'})/\bar{G} \right],
$$
which is a GIT quotient.

\subsubsection{Symplectic structure}
Let $\bar{G}_{\bR} \cong U(1)^{\bar{R}-3}$ be the maximal compact subgroup of $G$, which carries a Hamiltonian action on $\bC^{\bar{R}}$. Let $\bar{\fg}_{\bR}^* \cong \bR^{\bar{R}-3}$ be the dual of the Lie algebra $\bar{\fg}_{\bR}$ of $\bar{G}_{\bR}$. The moment map $\bar{\mu}: \bC^{\bar{R}} \to \bar{\fg}_{\bR}^*$ of the $\bar{G}_{\bR}$-action can be described as
$$
    \bar{\mu}(x_1, \dots, x_{\bar{R}}) = \left(\sum_{i = 1}^{\bar{R}} \bar{l}_i^{(1)}|x_i|^2, \dots, \sum_{i = 1}^{\bar{R}} \bar{l}_i^{(\bar{R}-3)}|x_i|^2  \right).
$$

Applying $\Hom(-, \bZ)$ to \eqref{eqn:XspFanSequence}, we obtain a short exact sequence
$$
    \xymatrix{
        0 \ar[r] &  M \ar[r]^{\bar{\alpha}^\vee} & \bZ^{\bar{R}} \ar[r]^{\bar{\psi}^\vee} & \bar{\bL}^\vee \ar[r] & 0.
    }
$$
There is a canonical identification $\bar{\fg}_{\bR}^* \cong \bar{\bL}_{\bR}^\vee$. Let $\{e_1^\vee, \dots, e_{\bar{R}}^\vee\}$ be the basis dual to $\{e_1, \dots, e_{\bar{R}}\}$. For each $i = 1, \dots, \bar{R}$, define
$$
    \bar{D}_i := \bar{\psi}^\vee(e_i^\vee) \in \bar{\bL}^\vee.
$$
The \emph{extended nef cone} of $\bar{\cX}$ is defined to be
$$
    \tNef(\bar{\cX}):= \bigcap_{\sigma \in \bar{\Sigma}(3)} \sum_{i \in I_\sigma} \bR_{\ge 0} \bar{D}_i,
$$
which is an $(\bar{R}-3)$-dimensional simplicial cone in $\bar{\bL}_{\bR}^\vee$. Let $r = (r_1, \dots, r_{\bar{R}-3})$ be a point in the interior of $\tNef(\bar{\cX})$, which can be viewed as an extended K\"ahler class of $\bar{\cX}$. Then $\bar{\cX}$ is the symplectic quotient
$$
    [\bar{\mu}^{-1}(r)/\bar{G}_{\bR}]
$$
and the standard K\"ahler form
$$
    \frac{\sqrt{-1}}{2} \sum_{i=1}^{\bar{R}} dx_i \wedge d\bar{x}_i
$$
on $\bC^{\bar{R}}$ descends to a K\"ahler form on $\bar{\cX}$. It further restricts to a K\"ahler form on $\cX \subseteq \bar{\cX}$.

\subsubsection{Aganagic-Vafa branes}
An \emph{Aganagic-Vafa brane} \cite{AV00,FLT12} $\cL$ in $\bar{\cX}$ is a Lagrangian suborbifold of form
$$
    \cL:= \left[ \left. \left\{ (x_1, \dots, x_{\bar{R}}) \in \bar{\mu}^{-1}(r) : \sum_{i = 1}^{\bar{R}} \bar{l}_i' |x_i|^2 = c', \sum_{i=1}^{\bar{R}} \bar{l}_i''|x_i|^2 = c'', \arg(\prod_{i=1}^{\bar{R}} x_i) = c'''       \right\}  \middle/ \bar{G}_{\bR} \right. \right],
$$
where $c', c'', c''' \in \bR$ are constants and the vectors $\bar{l}' = (l_1', \dots, l_{\bar{R}}'), \bar{l}'' = (l_1'', \dots, l_{\bar{R}}'') \in \bZ^{\bar{R}}$ satisfy
$$
    \sum_{i=1}^{\bar{R}} \bar{l}_i' = \sum_{i=1}^{\bar{R}} \bar{l}_i'' = 0.
$$
The brane $\cL$ is preserved under the action of the maximal compact subgroup $T_{\bR}' \cong U(1)^2$ of the Calabi-Yau torus $T'$. Moreover, it intersects a unique $T$-invariant line $\fl_{\tau}$ in $\bar{\cX}$, where $\tau \in \bar{\Sigma}(2)$. Suppose $\tau \in \Sigma(2) \subseteq \bar{\Sigma}(2)$. Then $\cL$ is contained in $\cX$ as a Lagrangian suborbifold. Suppose furthermore that $\tau \in \Sigma(2) \setminus \Sigma(2)_c$, in which case $\cL$ is called an \emph{outer} brane in $\cX$. 

We have $\cL \cap \fl_{\tau} = \cL \cap \fo_{\tau} \cong S^1 \times \cB G_\tau$ and the inclusions $\cL \cap \fl_{\tau} \to \cL$, $\cL \cap \fl_\tau \to \fo_\tau$ are homotopy equivalences. We have
$$
    H_1(\cL;\bZ) \cong \pi_1(\cL) \cong \pi_1(\fo_\tau) = H_\tau \cong \bZ \times G_\tau.
$$
Let $L \subset X$ be the coarse moduli space of $\cL$. Then the intersection of $L$ and $l_\tau \cong \bC$ in $X$ bounds a disk $B$ in $l_\tau$. We orient $B$ by the holomorphic structure of $X$. We have $L \cap l_\tau = \partial B \cong S^1$, and
$$
    H_1(L; \bZ) \cong \bZ[\partial B], \qquad H_2(X, L; \bZ) \cong H_2(X;\bZ) \oplus \bZ[B].
$$
For the effective classes, we have
$$
    \Eff(X,L) \cong \Eff(X) \oplus \bZ_{\ge 0}[B].
$$

\subsubsection{Additional combinatorial data}\label{sect:BraneCombData}

Let $\sigma \in \Sigma(3)$ be the unique $3$-cone that contains $\tau$ as a facet. Take indices $i_1, i_2, i_3 \in \{1, \dots, R'\}$ such that
$$
    I'_{\tau} = \{i_2, i_3\}, \qquad I'_{\sigma} = \{i_1, i_2, i_3\}
$$
and that $b_{i_1}, b_{i_2}, b_{i_3}$ appears in $N' \times \{1\}$ in counterclockwise order. There is a $\bZ$-basis $\{\sv_1, \sv_2, \sv_3\}$ of $N$ under which we have
$$
    b_{i_1} = \fr(\tau, \sigma) \sv_1- \fs(\tau, \sigma) \sv_2 + \sv_3, \quad b_{i_2} = \fm(\tau, \sigma) \sv_2 + \sv_3, \quad b_{i_3} = \sv_3
$$
for a unique $\fs(\tau, \sigma) \in \{0, \dots, \fr(\tau, \sigma)-1\}$. The $\bZ$-basis of $M$ dual to $\{\sv_1, \sv_2, \sv_3\}$ takes form $\{\su_1, \su_2, \su_3\}$ for some $\su_1, \su_2 \in M$. Here, $\su_3 \in M$ is the element specified by the Calabi-Yau condition in Section \ref{sect:TCY3}. We still use $\su_1, \su_2$ to denote their images under the under the projection $M \to M'$ by an abuse of notation.  Let $\tau_2, \tau_3 \in \Sigma(2)$ denote the other two facets of $\sigma$ with
$$
    I'_{\tau_2} = \{i_1, i_3\}, \qquad I'_{\tau_3} = \{i_1, i_2\}.
$$
We have the following description of tangent $T'$-weights at $\fp_\sigma$:
$$
    \sw(\tau, \sigma) = \frac{1}{\fr} \su_1, \quad \sw(\tau_2,\sigma) = \frac{\fs}{\fr\fm}\su_1 + \frac{1}{\fm}\su_2 ,\quad \sw(\tau_3,\sigma) = -\frac{\fm + \fs}{\fr\fm}\su_1 - \frac{1}{\fm}\su_2.
$$
where $\fr = \fr(\tau, \sigma)$, $\fm = \fm(\tau, \sigma), \fs = \fs(\tau, \sigma)$.

\subsubsection{Framing}
A \emph{framing} on the brane $\cL$ is a choice of a rational number $f \in \bQ$. Writing $f = \frac{\fb}{\fa}$ for coprime integers $\fa \in \bZ_{\ge 1}, \fb \in \bZ$ and viewing $\su_1, \su_2 \in M$ above as characters of the Calabi-Yau torus $T'$, the framing determines a 1-dimensional subtorus $T_f := \ker(\fa \su_2 - \fb \su_1) \cong \bC^*$ of $T'$.

\subsection{Multi-brane open geometry}\label{sect:OpenGeometry}
Let
$$
    \cL = \cL_1 \sqcup \cdots \sqcup \cL_s
$$
be a disjoint union of $s$ Aganagic-Vafa branes $\cL_1, \dots, \cL_s$ in $\cX$. For $i =1, \dots, s$, let $\tau_i \in \Sigma(2)$ be such that $\fl_{\tau_i}$ is the unique $T$-invariant line in $\cX$ that $\cL_i$ intersects. We make the following assumption on the brane configuration (cf. \cite[Assumption 2.3]{LY21}).

\begin{assumption}\label{assump:Branes} \rm{
We assume that $\tau_1, \dots, \tau_s$ are distinct. Moreover, we assume that each $\cL_i$ is an outer brane in a (and hence any) toric Calabi-Yau semi-projective partial compactification $\bar{\cX}$ of $\cX$, as in Section \ref{sect:AVBranes}. Equivalently, if $\bar{\Sigma}$ is the fan of the coarse moduli $\bar{X}$ of $\bar{\cX}$, then we assume that $\tau_i \in \Sigma(2) \setminus \bar{\Sigma}(2)_c$ for $i = 1, \dots, s$. In particular, each $\cL_i$ is an outer brane in $\cX$.
}
\end{assumption}

For $i = 1, \dots, s$, let $L_i \subset X$ be the coarse moduli space of $\cL_i$. The intersection of $L_i$ and $l_{\tau_i} \cong \bC$ in $X$ bounds a disk $B_i$ in $l_\tau$. Writing $L = L_1 \sqcup \cdots \sqcup L_s$, we have
$$
    H_1(L; \bZ) \cong \bigoplus_{i = 1}^s \bZ[\partial B_i], \qquad H_2(X, L; \bZ) \cong H_2(X;\bZ) \oplus \bigoplus_{i = 1}^s \bZ[B_i].
$$
For the effective classes, we have
$$
    \Eff(X,L) \cong \Eff(X) \oplus \bigoplus_{i = 1}^s \bZ_{\ge 0}[B_i].
$$
Moreover, we introduce the following notation for the data associated to the brane $\cL_i$ as in Section \ref{sect:BraneCombData}:
\begin{itemize}
    \item $\sigma_i \in \Sigma(3)$ denotes the unique 3-cone that contains $\tau_i$;
    
    \item $i_1^{(\tau_i, \sigma_i)}, i_2^{(\tau_i, \sigma_i)}, i_3^{(\tau_i, \sigma_i)} \in \{1, \dots, R'\}$;
    
    \item $\tau_{2,i}, \tau_{3,i} \in \Sigma(2)$ denote the two other facets of $\sigma_i$ with $I'_{\tau_{2,i}} = \{i_1^{(\tau_i, \sigma_i)}, i_3^{(\tau_i, \sigma_i)}\}$, $I'_{\tau_{3,i}} = \{i_1^{(\tau_i, \sigma_i)}, i_2^{(\tau_i, \sigma_i)}\}$;

    \item $\fr_i = \fr(\tau_i, \sigma_i)$, $\fm_i = \fm(\tau_i, \sigma_i)$, $\fs_i = \fs(\tau_i, \sigma_i)$;

    \item $\bZ$-basis $\{\sv_{1,i}, \sv_{2,i}, \sv_{3,i}\}$ of $N$;

    \item dual $\bZ$-basis $\{\su_{1,i}, \su_{2,i}, \su_3\}$ of $M$ (note again the independence of $\su_3$ on $i$);
    
    
    \item images $\su_{1,i}, \su_{2,i} \in M'$ of $\su_{1,i}, \su_{2,i}$ respectively under the projection $M \to M'$ (by abusive notation).
\end{itemize}
For each $i$, we have
\begin{equation}\label{eqn:sigmaWts}
    \sw_{1,i} := \sw(\tau_i, \sigma_i) = \frac{1}{\fr_i} \su_{1,i}, \quad 
    \sw_{2,i} := \sw(\tau_{2,i},\sigma_i) = \frac{\fs_i}{\fr_i\fm_i}\su_{1,i} + \frac{1}{\fm_i}\su_{2,i} ,\quad 
    \sw_{3,i} := \sw(\tau_{3,i},\sigma_i) = -\frac{\fm_i + \fs_i}{\fr_i\fm_i}\su_{1,i} - \frac{1}{\fm_i}\su_{2,i}.
\end{equation}

Now we choose a collection of \emph{parallel framings} for the Aganagic-Vafa branes above, in the following sense. Take $\su_1, \su_2 \in M$ such that $\{\su_1, \su_2, \su_3\}$ is a $\bZ$-basis, which will be used as reference globally. By an abuse of notation, we still use $\su_1, \su_2$ to denote their images under the projection $M \to M'$, which then form a $\bZ$-basis of $M'$. Choose $f = \frac{\fb}{\fa} \in \bQ$ where $\fa \in \bZ_{\ge 1}, \fb \in \bZ$ are coprime. Consider the 1-dimensional subtorus $T_f := \ker(\fa \su_2 - \fb \su_1) \cong \bC^*$ of $T'$. Let $M_f := \Hom(T_f, \bC^*) \cong \bZ$ denote the character lattice of $T_f$. We use the following notation for the natural projection
$$
    \big|_{\su_2 - f\su_1 = 0} : M' \to M_f, \qquad \su_1 \mapsto \fa \su, \su_2 \mapsto \fb \su 
$$
where $\su$ is a generator of $M_f$.

\begin{assumption}\label{assump:Parallel} \rm{
We assume that $f \in \bQ$ is chosen such that
$$
    \su_{1,i} \big|_{\su_2 - f\su_1 = 0} \neq 0 \quad \in M_f \cong \bZ
$$
for all $i = 1, \dots, s$.
}
\end{assumption}

Note that this is achieved for a \emph{generic} choice of $f \in \bQ$, i.e. all but finitely many rational numbers. Then, we define the framing $f_i$ of the brane $\cL_i$ by
$$
    f_i := \frac{\su_{2,i}\big|_{\su_2 - f\su_1 = 0}}{\su_{1,i}\big|_{\su_2 - f\su_1 = 0}} \quad \in \bQ.
$$
In other words, the non-zero vectors
$$
    \su_2 - f \su_1, \su_{2,1} - f_1 \su_{1,i}, \dots, \su_{2,s} - f_s \su_{1,s} \quad \in M'_{\bQ}
$$
are all parallel, spanning a 1-dimensional linear subspace that does not contain any of $\su_{1,i}, \dots, \su_{1,s}$. The framing subtori of $T'$ defined by the framed branes $(\cL_i, f_i)$ are all identified with $T_f$. 
For $i = 1, \dots, s$, let $\fa_i \in \bZ_{\ge 1}$, $\fb_i \in \bZ$ be coprime integers such that
$$
    f_i = \frac{\fb_i}{\fa_i}.
$$

We adopt the following notation which view the characters of the tori as equivariant parameters:
$$
    H^*_T(\pt;\bZ) = \bZ[\su_1, \su_2, \su_3], \quad H^*_{T'}(\pt;\bZ) = \bZ[\su_1, \su_2], \quad H^*_{T_f}(\pt;\bZ) = \bZ[\su].
$$
Over $\bQ$, we introduce the following notation for the fields of fractions:
$$
    \cQ_T := \bQ(\su_1, \su_2, \su_3), \quad \cQ_{T'} := \bQ(\su_1, \su_2), \quad \cQ_{T_f} := \bQ(\su).
$$

From now on, we use the triple $(\cX, \cL, f)$ to denote the open geometry constructed here.

\subsection{Closed geometry}\label{sect:LocalGeometry}
In this section, we construct the closed geometry $\tcX$ corresponding to the open geometry $(\cX, \cL, f)$, which is a toric Calabi-Yau ($3+s$)-orbifold determined by an extended stacky fan $\tbSi = (\tN, \tSi, \talpha)$. Here,
$$
    \tN := N \oplus \bZ \sv_4 \oplus \cdots \oplus \bZ \sv_{3+s}
$$
is a lattice of rank $3+s$ and $\talpha: \bZ^{R + 2s} \to \tN$ is the homomorphism specified by the following images of the standard basis vectors $\tb_i := \talpha(e_i)$, $i = 1, \dots, R+2s$:
\begin{align*}
    &\tb_i = b_i && \text{for } i = 1, \dots, R;\\
    &\tb_{R+2i-1} = -\fa_i \sv_{1,i} - \fb_i \sv_{2,i} + \sv_{3,i} + \sv_{3+i}, \quad \tb_{R+2i} = \sv_{3,i} + \sv_{3+i} && \text{for } i = 1, \dots, s.
\end{align*}
The set of rays in the fan $\tSi$ is
$$
    \tSi(1) = \{\trho_1, \dots, \trho_{R'}, \trho_{R+1}, \dots, \trho_{R+2s}\}
$$
where $\trho_i := \bR_{\ge 0} \tb_i$. We fix a notation for the following subset of indices
$$
    I_s := \{R+2, R+4, \dots, R+2s\}.
$$
We describe the maximal cones in $\tSi$ below whose faces will then constitute the fan. 

\subsubsection{Maximal cones}
The set of maximal cones can be partitioned into two types as
$$
    \tSi(3+s) = \iota(\Sigma(3)) \sqcup \{\tsi_1, \dots, \tsi_s\}.
$$
First, for any $\sigma \in \Sigma(3)$, there is a maximal cone $\iota(\sigma) \in \tSi(3+s)$ defined by the index set
$$
    I'_{\iota(\sigma)} = I'_{\sigma} \sqcup I_s.
$$
Indeed, this definition is extended to all cones in $\Sigma$, yielding injective maps
$$
    \iota: \Sigma(d) \to \Sigma(d+s)
$$
for $d = 0, 1, 2, 3$. Second, for $i = 1, \dots, s$, there is an additional maximal cone $\tsi_i$ defined by the index set
$$
    I'_{\tsi_i} = \{i_2^{(\tau_i, \sigma_i)}, i_3^{(\tau_i, \sigma_i)}, R+2i-1\} \sqcup I_s.
$$
Note that $\tsi_i$ is the only maximal cone that contains the ray $\trho_{R+2i-1}$. Moreover, $\tsi_i$ shares the facet $\iota(\tau_i)$ with $\iota(\sigma_i)$, which makes $\iota(\tau_i)$ an element of $\tSi(2+s)_c$. Indeed, we have
$$
    \tSi(2+s)_c = \iota(\Sigma(2)_c \sqcup \{\tau_1, \dots, \tau_s\}).
$$

\subsubsection{Stabilizers}
Under the inclusion $N \subset \tN$, $\Sigma$ is contained in $\tSi$ as a subfan which corresponds to the inclusion
$$
    \iota: \cX \subset \cX \times \bC^s \subset \tcX.
$$
For any $\sigma \in \Sigma(3)$, $\Box(\iota(\sigma))$ is identified with $\Box(\sigma)$. The stabilizers of the extra maximal cones are described by the following result which is a straightforward generalization of \cite[Lemma 2.4]{LY22}.

\begin{lemma}\label{lem:Age}
For $i = 1, \dots, s$, $G_{\tsi_i}$ is a cyclic group of order $\fa_i \fm_i$. Elements of age at most 1 are exactly those contained in the subgroup $G_{\tau_i} \cong \mu_{\fm_i}$. All other elements have age $2$.
\end{lemma}

Therefore, there is an inclusion $\Box(\cX) \subseteq \Box(\tcX)$ where the additional elements in $\Box(\tcX)$ all have age $2$.

\subsubsection{Calabi-Yau subtorus}

Let $\{\sv_1, \sv_2, \sv_3\}$ be the $\bZ$-basis of $N$ dual to the reference basis $\{\su_1, \su_2, \su_3\}$ of $M$ chosen in Section \ref{sect:OpenGeometry}, which extends to the $\bZ$-basis $\{\sv_1, \sv_2, \sv_3, \sv_4, \dots, \sv_{3+s}\}$ of $\tN$. Let $\{\su_1, \su_2, \su_3, \su_4, \dots, \su_{3+s}\}$ be the dual $\bZ$-basis of $\tM:= \Hom(\tN, \bZ)$. Here, we again abuse notation since $\su_1, \su_2, \su_3 \in \tM$ are lifts of their counterparts in $M$ under the natural projection $\tM \to M$. Then $\inner{\su_3, \tb_i} = 1$ for any $i = 1, \dots, R+2s$, which implies that $\tcX$ is Calabi-Yau. Let $\tT := \tN \otimes \bC^* \cong (\bC^*)^{3+s}$ be the algebraic torus of $\tcX$. Let $\tN' := \ker(\su_3) \cong \bZ^{2+s}$, $\tT' := \tN' \otimes \bC^* \cong (\bC^*)^{2+s}$ be the Calabi-Yau subtorus, and $\tM' := \Hom(\tN', \bZ)$. There are commutative squares of inclusions and natural projections
$$
    \xymatrix{
        \tN & N \ar[l]\\
        \tN' \ar[u] & N', \ar[u] \ar[l]
    } \qquad \qquad 
    \xymatrix{
        \tM \ar[r] \ar[d] & M \ar[d]\\
        \tM' \ar[r] & M'
    }
$$
where the projection $\tM' \to M'$ is realized by the weight restriction $\big|_{\su_4 = \cdots = \su_{3+s} = 0}$. We have
$$
    H^*_{\tT}(\pt) = \bZ[\su_1, \su_2, \su_3, \su_4, \dots, \su_{3+s}], \quad 
    H^*_{\tT'}(\pt) = \bZ[\su_1, \su_2, \su_4, \dots, \su_{3+s}].
$$
Over $\bQ$, we let
$$
    \cQ_{\tT} := \bQ(\su_1, \su_2, \su_3, \su_4, \dots, \su_{3+s}), \quad \cQ_{\tT'} := \bQ(\su_1, \su_2, \su_4, \dots, \su_{3+s})
$$
denote the fields of fractions.




Given a flag $(\ttau, \tsi) \in F(\tSi)$, we define
$$
    \tsw(\ttau, \tsi) := c_1^{\tT'}(T_{\fp_{\tsi}}\fl_{\ttau}) \in H^2_{\tT'}(\pt;\bQ) \cong \tM'_{\bQ}.
$$
The Calabi-Yau condition implies that for all $\tsi \in \tSi(3+s)$, we have
\begin{equation}\label{eqn:CYtX}
    \sum_{(\ttau, \tsi) \in F(\tSi)} \tsw(\ttau, \tsi) = 0.
\end{equation}

\subsubsection{Tangent weights}
We now describe the tangent $\tT'$-weights in more detail. 
First, consider a maximal cone $\iota(\sigma) \in \tSi(3+s)$ for some $\sigma \in \Sigma(3)$. For each of the three facets $\tau \in \Sigma(2)$ of $\sigma$, we have
$$
    \text{$\fr(\iota(\tau), \iota(\sigma)) = \fr(\tau, \sigma)$}, \quad \tsw(\iota(\tau), \iota(\sigma)) \big|_{\su_4 = \cdots = \su_{3+s} = 0} = \sw(\tau, \sigma).
$$
Moreover, for $j = 1, \dots, s$, there is a facet $\delta_{3+j}(\sigma)$ of $\iota(\sigma)$ defined by the index set
$$
    I'_{\delta_{3+j}(\sigma)} = I'_{\sigma} \sqcup I_s \setminus \{R+2j\}.
$$
We have
$$
    \fr(\delta_{3+j}(\sigma), \iota(\sigma)) = 1, \quad \tsw(\delta_{3+j}(\sigma), \iota(\sigma)) = \su_{3+j}.
$$

Second, consider the maximal cone $\tsi_i \in \tSi(3+s)$ for $i = 1, \dots, s$. We define and study its facets
$$
    \{\iota(\tau_i), \delta_2(\tsi_i), \delta_3(\tsi_i), \delta_4(\tsi_i), \dots, \delta_{3+s}(\tsi_i)\}
$$
in steps below.
\begin{itemize}
    \item For the facet $\iota(\tau_i)$, we have
        $$
            \fr(\iota(\tau_i), \tsi_i) = \fa_i, \quad \tsw(\iota(\tau_i), \tsi_i) = -\frac{\fr_i}{\fa_i} \tsw(\iota(\tau_i), \iota(\sigma_i)) = -\frac{1}{\fa_i}\su_{1,i}.
        $$

    \item There are two other facets $\delta_2(\tsi_i), \delta_3(\tsi_i)$ defined by the index sets
        $$
            I'_{\delta_2(\tsi_i)} = \{i_3^{(\tau_i, \sigma_i)}, R+2i-1\} \sqcup I_s, \quad I'_{\delta_3(\tsi_i)} = \{i_2^{(\tau_i, \sigma_i)}, R+2i-1\} \sqcup I_s.
        $$
        We have $\fr(\delta_2(\tsi_i), \tsi_i) = \fr(\delta_3(\tsi_i), \tsi_i) = 1$, and
        $$
            \tsw(\delta_2(\tsi_i), \tsi_i) = \tsw(\iota(\tau_{2,i}), \iota(\sigma_i)) - \frac{\fr_i f_i -\fs_i}{\fr_i\fm_i} \su_{1,i}, \quad \tsw(\delta_3(\tsi_i), \tsi_i) = \tsw(\iota(\tau_{3,i}), \iota(\sigma_i)) - \frac{-\fr_i f_i - \fm_i + \fs_i}{\fr_i\fm_i}  \su_{1,i}.
        $$
        In particular, the vectors
        $$
            \tsw(\delta_2(\tsi_i), \tsi_i)  \big|_{\su_4 = \cdots = \su_{3+s} = 0} = -\frac{f_i}{\fm_i}\su_{1,i} + \frac{1}{\fm_i}\su_{2,i}, \quad \tsw(\delta_3(\tsi_i), \tsi_i)  \big|_{\su_4 = \cdots = \su_{3+s} = 0} = \frac{f_i}{\fm_i}\su_{1,i} - \frac{1}{\fm_i}\su_{2,i},
        $$
        add up to zero in $M'$ and are non-zero scalar multiples of $\su_2 - f\su_1$. 
        
    \item There is a fourth facet $\delta_{3+i}(\tsi_i)$ defined by the index set
        $$
            I'_{\delta_{3+i}(\tsi_i)} = \{i_2^{(\tau_i, \sigma_i)}, i_3^{(\tau_i, \sigma_i)}, R+2i-1\} \sqcup I_s \setminus \{R+2i\}.
        $$
        We have
        $$
            \fr(\delta_{3+i}(\tsi_i), \tsi_i) = 1, \quad \tsw(\delta_{3+i}(\tsi_i), \tsi_i) = \frac{1}{\fa_i}\su_{1,i} + \su_{3+i}.
        $$    
        We illustrate the four facets discussed so far in Figure \ref{fig:ExtraConeFacets23}.

        \begin{figure}[h]
            \begin{tikzpicture}
                \draw (-2.5,-1.5) -- (-1,0) -- (-1,1.5); 
                \draw (-3, 0.5) -- (-1,0) -- (3,0) -- (5,0.5);
                \draw (3, 1.5) -- (3, -1.5);
                \node at (-0.8,-0.3) {\small $\iota(\sigma_i)$};
                \node[below] at (1,0) {\small $\iota(\tau_i)$};
                \node at (2.8,-0.3) {\small $\tsi_i$};
                \node[right] at (3,-1.3) {\small $\delta_3(\tsi_i)$};
                \node[right] at (3,1.3) {\small $\delta_2(\tsi_i)$};
                \node[right] at (5,0.5) {\small $\delta_{3+i}(\tsi_i)$};
                \node[left] at (-1,1.3) {\small $\iota(\tau_{2,i})$};
                \node[left] at (-2.3, -1.3) {\small $\iota(\tau_{3,i})$};
                \node[left] at (-3, 0.5) {\small $\delta_{3+i}(\sigma_i)$};
                \node[above] at (-0.3, 0) {\small $\frac{1}{\fr_i}\su_{1,i}$};
                \node[above] at (2.2, 0) {\small $-\frac{1}{\fa_i}\su_{1,i}$};
                \node at (-1.8, 0) {\small $\su_{3+i}$};
                \node at (4.3, -0.1) {\small $\frac{1}{\fa_i}\su_{1,i} + \su_{3+i}$};
                \node at (1, 1.2) {\small $\cO\left(\frac{\fr_i f_i -\fs_i}{\fr_i\fm_i}\right)$};
                \node at (0.7, -1.2) {\small $\cO\left(-\frac{\fr_i f_i +\fm_i -\fs_i}{\fr_i\fm_i}\right)$};
            \end{tikzpicture}

            \caption{$\tT'$-invariant lines and tangent weights associated to facets $\iota(\tau_i)$, $\delta_2(\tsi_i)$, $\delta_3(\tsi_i)$, $\delta_{3+i}(\tsi_i)$ of $\tsi_i$ and their corresponding facets of $\iota(\sigma_i)$.}
            \label{fig:ExtraConeFacets23}
        \end{figure}
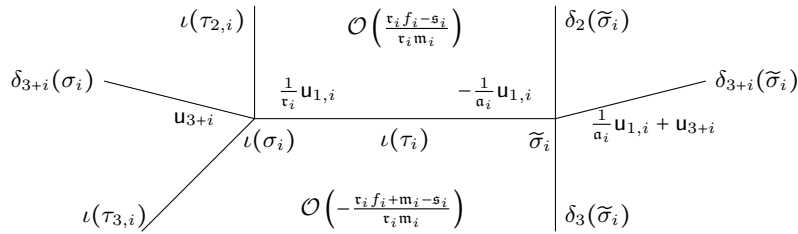

    \item Finally, for $j = 1, \dots, s$, $j \neq i$, there is a facet $\delta_{3+j}$ defined by the index set
    $$
        I'_{\delta_{3+j}(\tsi_i)} = \{i_2^{(\tau_i, \sigma_i)}, i_3^{(\tau_i, \sigma_i)}, R+2i-1\} \sqcup I_s \setminus \{R+2j\}.
    $$
    We have
    $$
        \fr(\delta_{3+j}(\tsi_i), \tsi_i) = 1, \quad \tsw(\delta_{3+j}(\tsi_i), \tsi_i) = \su_{3+j}.
    $$
\end{itemize}

\subsection{Formal relative geometry}\label{sect:RelativeGeometry}
We bridge the 3-dimensional open geometry $(\cX, \cL, f)$ and the ($3+s$)-dimensional closed geometry by a sequence of geometries in the intermediate dimensions. The starting point is a 3-dimensional formal relative geometry $(\hcY, \hcD)$. From the viewpoint of the open/relative Gromov-Witten correspondence, as in \cite{LLLZ09,FL13,FLT12,LY21}, the Lagrangian boundary conditions in $(\cX, \cL, f)$ can be interpreted as tangency conditions to the new (formal) toric divisor $\hcD$.

We construct $(\hcY, \hcD)$ as the relative FTCY 3-orbifold associated to an FTCY graph $\sGa = (\hGa, \sfp, \sff_2, \sff_3)$ defined as follows. For the underlying regular weakly trivalent graph
$$
    \hGa = \left( V(\hGa) = V^3(\hGa) \sqcup V^1(\hGa), E(\hGa) \right),
$$
we take
$$
    V^3(\hGa) := \Sigma(3), \quad E(\hGa) := \Sigma(2).
$$
Each flag $(\tau, \sigma) \in F(\Sigma)$ determines a flag $(e,v)$ in $F(\hGa)$. Moreover, we take $V^1(\hGa)$ to be a set with $s$ vertices, denoted by
$$
    V^1(\hGa) := \{\tv_1, \dots, \tv_s\}.
$$
For $i = 1, \dots, s$, let $v_i \in V^3(\hGa)$ (resp. $e_i \in E(\hGa)$) denote the trivalent vertex (resp. edge) corresponding to the cone $\sigma_i \in \Sigma(3)$ (resp. $\tau_i \in \Sigma(2)$) determined by $\cL_i$. We let $e_i$ be the only edge incident to the univalent vertex $\tv_i$. In other words, we have
$$
    F(\hGa) = F(\Sigma) \sqcup \{(e_1, \tv_1), \dots, (e_s, \tv_s)\}.
$$
The set of compact edges in $\hGa$ is then
$$
    E(\hGa)_c = \Sigma(2)_c \sqcup \{e_1, \dots, e_s\}.
$$

The position map $\sfp: F(\hGa) \to M'_{\bQ} \setminus \{0\}$ is defined as follows. For $(e,v) \in F(\hGa)$ corresponding to some $(\tau, \sigma) \in F(\Sigma)$, we set
$$
    \sfp(e,v) := \sw(\tau, \sigma).
$$
For the flag $(e_i, \tv_i)$ at the univalent vertex $\tv_i$, we set
$$
    \sfp(e_i, \tv_i) := -\frac{1}{\fa_i}\su_{1,i}.
$$
The framing maps $\sff_2, \sff_3: V^1(\hGa) \to M'_{\bQ} \setminus \{0\}$ are defined by
$$
    \sff_2(\tv_i) := -\frac{f_i}{\fm_i}\su_{1,i} + \frac{1}{\fm_i}\su_{2,i}, \quad \sff_3(\tv_i) := \frac{f_i}{\fm_i}\su_{1,i} - \frac{1}{\fm_i}\su_{2,i}.
$$
The compatibility conditions at the vertices and compact edges in the FTCY graph are directly verified. In particular, along the compact edge $e_i$, the bijection between $E_{v_i} \setminus \{e_i\}$ and $\{2,3\}$ identifies the edge corresponding to $\tau_{2,i}$ (resp. $\tau_{3,i}$) with 2 (resp. 3).

Then, we define $(\hcY, \hcD)$ to be the relative FTCY 3-orbifold associated to $\sGa$, where $\hcD$ has $s$ connected components
$$
    \hcD = \hcD_1 \sqcup \cdots \sqcup \hcD_s.
$$

\subsection{Formal intermediate geometries}\label{sect:IntermediateGeometry}
Based on the relative FTCY 3-orbifold $(\hcY, \hcD)$, we now construct a sequence of intermediate geometries
$$
    (\hcY^{(0)}, \hcD^{(0)}) = (\hcY, \hcD), (\hcY^{(1)}, \hcD^{(1)}), \dots, (\hcY^{(s)}, \hcD^{(s)})
$$
where for $\ell = 0, \dots, s$, $(\hcY^{(\ell)}, \hcD^{(\ell)})$ is a relative FTCY ($3+\ell$)-orbifold and the divisor $\hcD^{(\ell)}$ has $s-\ell$ connected components. Non-equivariantly, for $\ell = 1, \dots, s$, $\hcY^{(\ell)}$ can be viewed as the geometric formal vector bundle associated to the rank $\ell$ locally free sheaf
$$
    \bigoplus_{i = 1}^\ell \cO_{\hcY}(-\hcD_i)
$$
on the formal scheme $\hcY$, and 
$$
    \hcD^{(\ell)} = \hcD^{(\ell)}_{\ell+1} \sqcup \cdots \hcD^{(\ell)}_{s}
$$
where $\hcD^{(\ell)}_i$ is the preimage of the divisor $\hcD_i$ under the projection $\hcY^{(\ell)} \to \hcY$. Indeed, this is motivated by the construction of intermediate geometries in the log/local correspondence as described in \cite[Section 1.4.1]{BBvG20} for log pairs whose boundary divisor has multiple components.\footnote{We note that our indexing convention is slightly different from that in \cite{BBvG20}.} At each intermediate step, $\hcY^{(\ell)}$ can be viewed as the geometric formal line bundle associated to the invertible sheaf
$$
    \cO_{\hcY^{(\ell-1)}}(-\hcD^{(\ell-1)}_\ell)
$$
on $\hcY^{(\ell-1)}$, where the relative condition $\hcD^{(\ell-1)}_\ell$ is replaced by a local condition. At the end, $\hcD^{(s)}$ is empty and $\hcY^{(s)}$ is an FTCY ($3+s$)-orbifold which is indeed identified with the formal completion of the toric Calabi-Yau ($3+s$)-orbifold $\tcX$ along its toric 1-skeleton $\tcX^1$.


Now we give the details of the construction via FTCY graphs. In the process, we clarify the tangent weights of the Calabi-Yau tori actions on the intermediate geometries. First, for $\ell = 1, \dots, s$, let
$$
    N_{(\ell)} := N \oplus \bZ \sv_4 \oplus \cdots \oplus \bZ \sv_{3+\ell}, \qquad T_{(\ell)} := N_{(\ell)} \otimes \bC^*, \qquad M_{(\ell)} := \Hom(N_{(\ell)}, \bZ).
$$
In particular, $N_{(s)} = \tN$, $T_{(s)} = \tT$, $M_{(s)} = \tM$. There is a sequence of projections
$$
    \xymatrix{
        \tM = M_{(s)} \ar[r] & M_{(s-1)} \ar[r] & \cdots \ar[r] & M_{(1)} \ar[r] & M.
    }
$$
For $N_{(\ell)}$, we use the $\bZ$-basis $\{\sv_1, \dots, \sv_{3+\ell}\}$, whose dual basis of $M_{(\ell)}$ is $\{\su_1, \dots, \su_{3+\ell}\}$ which are the images of the corresponding elements of $\tM$ under the projection above. We define 
$$
    N_{(\ell)}' := \ker(\su_3) \subset N_{(\ell)}, \qquad T_{(\ell)}' := N_{(\ell)}' \otimes \bC^*, \qquad M_{(\ell)}' := \Hom(N_{(\ell)}', \bZ).
$$
There is a sequence of projections
$$
    \xymatrix{
        \tM' = M_{(s)}' \ar[r]^{\big|_{\su_{3+s} = 0}} & M_{(s-1)}' \ar[r]^{\big|_{\su_{2+s} = 0}} & \cdots \ar[r] & M_{(1)}' \ar[r]^{\big|_{\su_4 = 0}} & M'.
    }
$$

For $\ell = 1, \dots, s$, we construct $(\hcY^{(\ell)}, \hcD^{(\ell)})$ as the relative FTCY ($3+\ell$)-orbifold associated to an FTCY graph $\sGa^{(\ell)} = (\hGa^{(\ell)}, \sfp^{(\ell)}, \sff_2^{(\ell)}, \dots, \sff_{3+\ell}^{(\ell)})$ defined below. Extending the definitions to $\ell = 0$ will recover the graph $\sGa^{(0)} = \sGa$ and $(\hcY^{(0)}, \hcD^{(0)}) = (\hcY, \hcD)$. First, we define the underlying regular weakly ($3+\ell$)-valent graph
$$
    \hGa^{(\ell)} = \left( V(\hGa^{(\ell)}) = V^{3+\ell}(\hGa^{(\ell)}) \sqcup V^1(\hGa^{(\ell)}), E(\hGa^{(\ell)}) \right).
$$
For vertices, we take
$$
    V(\hGa^{(\ell)}) = V(\hGa) = \Sigma(3) \sqcup \{\tv_1, \dots, \tv_s\}.
$$
It is convenient to identify this set with the set of maximal cones
$$
    \tSi(3+s) = \iota(\Sigma(3)) \sqcup \{\tsi_1, \dots, \tsi_s\}
$$
in the fan $\tSi$ of $\tcX$. We take the partition
$$
    V^{3+\ell}(\hGa^{(\ell)}) := \iota(\Sigma(3)) \sqcup \{\tsi_1, \dots, \tsi_\ell\}, \quad V^1(\hGa^{(\ell)}) := \{\tsi_{\ell+1}, \dots, \tsi_s\}.
$$
For edges, we take the following subset of facets
$$
    E(\hGa^{(\ell)}) := \iota(\Sigma(2)) \sqcup \{\delta_{3+j}(\sigma) : \sigma \in \Sigma(3), 1 \le j \le \ell\} \sqcup \{\delta_{3+j}(\tsi_i) : 1 \le i,j \le \ell\}.
$$
The set of flags $F(\hGa^{(\ell)})$ is determined by the subset of $F(\tSi)$ involving the cones above. Note that the choice of $E(\hGa^{(\ell)})$ satisfies the valency requirements of all the vertices. The set of compact edges in $\hGa^{(\ell)}$ is
$$
    E(\hGa^{(\ell)})_c = \tSi(2+s)_c = \iota(\Sigma(2)_c \sqcup \{\tau_1, \dots, \tau_s\}),
$$
which is identified with $E(\hGa)_c$.

The position map $\sfp^{(\ell)}: F(\hGa^{(\ell)}) \to M'_{(\ell), \bQ} \setminus \{0\}$ is defined as follows. For $(e,v) \in F(\hGa)$ corresponding to some $(\ttau, \tsi) \in F(\tSi)$, we set
$$
    \sfp^{(\ell)}(e,v) := \tsw(\ttau, \tsi) \big|_{\su_{4+\ell} = \cdots = \su_{3+s} = 0}.
$$
For $j = -1, 0, 1, \dots, \ell$, the framing map $\sff_{3+j}^{(\ell)}: V^1(\hGa^{(\ell)}) \to M'_{(\ell),\bQ} \setminus \{0\}$ is defined on $\tsi_i$, $i = \ell+1, \dots, s$, by
$$
    \sff_{3+j}^{(\ell)}(\tsi_i) := \tsw(\delta_{3+j}(\tsi_i), \tsi_i) \big|_{\su_{4+\ell} = \cdots = \su_{3+s} = 0}.
$$

The compatibility conditions at the vertices follow from the Calabi-Yau condition on $\tcX$, specifically \eqref{eqn:CYtX}. At each $\sigma \in \Sigma(3)$, we use $\tsw(\delta_{3+j}(\sigma), \iota(\sigma)) = \su_{3+j}$ for $j = \ell+1, \dots, s$ to verify that
$$
    \sum_{(\ttau, \iota(\sigma)) \in F(\hGa^{(\ell)})} \tsw(\ttau, \iota(\sigma)) \big|_{\su_{4+\ell} = \cdots = \su_{3+s} = 0} = \sum_{(\ttau, \iota(\sigma)) \in F(\tSi)}  \tsw(\ttau, \iota(\sigma)) \big|_{\su_{4+\ell} = \cdots = \su_{3+s} = 0} = 0.
$$
At the univalent vertex $\tsi_i$, $i = \ell+1, \dots, s$, we use
$$
    \tsw(\delta_{3+j}(\tsi_i), \tsi_i) = \begin{cases}
            -\tsw(\iota(\tau_i), \tsi_i) + \su_{3+j} & \text{if } j = i\\
            \su_{3+j} & \text{if } j \in \{\ell+1, \dots, s\} \setminus \{i\}
        \end{cases}
$$
to verify that
$$
    \sum_{j = -1}^\ell \sff_{3+j}^{(\ell)}(\tsi_i) = \sum_{j = -1}^\ell \tsw(\delta_{3+j}(\tsi_i), \tsi_i) \big|_{\su_{4+\ell} = \cdots = \su_{3+s} = 0} = \sum_{(\ttau, \tsi_i) \in F(\tSi)}  \tsw(\ttau, \tsi_i) \big|_{\su_{4+\ell} = \cdots = \su_{3+s} = 0} = 0.
$$

The compatibility conditions along compact edges are also directly verified. In particular, at $\iota(\tau_i)$, $i = 1, \dots, \ell$, the bijection between $E_{\iota(\sigma_i)} \setminus \{\iota(\tau_i)\}$ and $E_{\tsi_i} \setminus \{\iota(\tau_i)\}$ identifies $\iota(\tau_{2,i})$ with $\delta_2(\tsi_i)$, $\iota(\tau_{3,i})$ with $\delta_3(\tsi_i)$, and $\delta_{3+j}(\sigma_i)$ with $\delta_{3+j}(\tsi_i)$ for $j = 1, \dots, \ell$. At $\iota(\tau_i)$, $i = \ell+1, \dots, s$, the bijection between $E_{\iota(\sigma_i)} \setminus \{\iota(\tau_i)\}$ and $\{2,3, \dots, 3+\ell\}$ identifies $\iota(\tau_{2,i})$ with $2$, $\iota(\tau_{3,i})$ with $3$, and $\delta_{3+j}(\sigma_i)$ with $3+j$ for $j = 1, \dots, \ell$. 

Then, we define $(\hcY^{(\ell)}, \hcD^{(\ell)})$ to be the relative FTCY ($3+\ell$)-orbifold associated to $\sGa^{(\ell)}$, where $\hcD^{(\ell)}$ has $s-\ell$ connected components
$$
    \hcD^{(\ell)} = \hcD_{\ell+1}^{(\ell)} \sqcup \cdots \sqcup \hcD_s^{(\ell)}.
$$
In particular, for $\ell = s$, it follows from the definitions that $\hcD^{(s)}$ is empty and $\hcY^{(s)}$ is the formal completion of $\tcX$ along $\tcX^1$ (cf. Example \ref{ex:GenuineFTCY}).


\section{Open/relative/local Gromov-Witten invariants}\label{sec:GW}
In this section, we define the different types of Gromov-Witten invariants involved in the correspondences and compute them by virtual localization \cite{GP99, GV05, Liu13}.

\subsection{Open Gromov-Witten invariants}\label{sect:OpenGW}
We start with the open Gromov-Witten invariants of $(\cX, \cL, f)$ which are virtual counts of twisted open stable maps. We follow the definition of \cite{FLT12,LY22}, although as in the case of relative invariants, we will only consider stable maps from genus-zero domains and with maximal winding at each of the $s$ components of $\cL$.


Let $n\in\bZ_{\geq 0}$, $\beta \in \Eff(X)$, and $(d_i, \lambda_i) \in \bZ_{\ge 1} \times G_{\tau_i} \subset H_1(\cL_i;\bZ)$ for $i = 1, \dots, s$. We write
$$
	\bd := ((d_1, \lambda_1), \dots, (d_s, \lambda_s)).
$$
Moreover, we define the effective class
$$
	\beta' := \beta + \sum_{i = 1}^s d_i[B_i] \qquad \in \Eff(X,L).
$$
Let
$$
	\Mbar_{(0,s),n}(\cX,\cL, \beta',\bd)
$$
be the moduli space of maps
$$u:((\cC,\underline{\fx}),\partial\cC)\to(\cX,\cL)$$
where:
\begin{itemize}
    \item $(\cC,\underline{\fx}) = (\cC,\fx_1,\dots,\fx_n)$ is a prestable bordered orbifold Riemann surface of topological type $(0,s)$ with $n$ interior marked points $\fx_1,\dots,\fx_n$, where stacky points are only allowed at the nodes (which are all balanced) and the interior marked points. Let $C$ denote the coarse moduli space of $\cC$ and $x_i$ denote the image of $\fx_i$. 
	Then $(C,x_1,\dots,x_n)$ is a prestable bordered Riemann surface of topological type $(0,s)$ with $n$ interior marked points.
    \item Let $\nu:\hat{\cC}\to\cC$ be the normalization map, where $\hat{\cC}$ is a possibly disconnected bordered orbifold Riemann surface with no nodes. Then $\nu\circ u:\hat{\cC}\to\cX$ is holomorphic.
    \item The automorphism group of $u$ is finite.
    \item Let $\bar{u}:(C,\partial C)\to(X,L)$ be the induced map on the level of coarse moduli. Then $\bar{u}_*[C]=\beta'$. 
	
	\item There is a labeling of the $s$ connected components of $\partial \cC$ as $(\partial \cC)_1, \dots, (\partial \cC)_s$ such that $u((\partial\cC)_i) \subset \cL_i$ and $u_*[(\partial\cC)_i]=(d_i,\lambda_i)$.
\end{itemize}
There are evaluation maps
$$
	\ev_i: \Mbar_{(0,s),n}(\cX,\cL, \beta',\bd) \to \cI \cX, \qquad i = 1, \dots, n
$$
at the marked points $\fx_i$. The action of $\bT'_{\bR}$ on $(\cX, \cL)$ induces an action on $\Mbar_{(0,s),n}(\cX,\cL, \beta',\bd)$ under which the evaluation maps are $\bT'_{\bR}$-equivariant and the $\bT'_{\bR}$-fixed locus is proper.

\begin{definition}\label{def:OpenGW}
Let $\gamma_1, \dots, \gamma_n \in H^2_{\CR,T'_{\bR}}(\cX;\bQ)=H^2_{\CR, T'}(\cX; \bQ)$. We define the open Gromov-Witten invariant
$$
	\langle \gamma_1,\dots,\gamma_n  \rangle^{\cX, \cL, T_f}_{\beta',\bd} := \int_{[\Mbar_{(0,s),n}(\cX,\cL, \beta',\bd)^{T'_{\bR}}]^\vir} \frac{\iota^*\left(\prod_{i = 1}^n \ev_i^*(\gamma_i)\right)}{e_{T'_{\bR}}(N^\vir)} \bigg|_{\su_2 - f\su_1 = 0} \qquad \in \bQ
$$
where $\iota: \Mbar_{(0,s),n}(\cX,\cL, \beta',\bd)^{T'_{\bR}} \to \Mbar_{(0,s),n}(\cX,\cL, \beta',\bd)$ is the inclusion of the $T'_{\bR}$-fixed locus and $N^\vir$ is the virtual normal bundle. 
\end{definition}

Note that the integral above is an element of $\cQ_{T'}$ that is homogenous of degree $0$, and the weight restriction $\big|_{\su_2 - f\su_1 = 0}$ is an element of $\cQ_{T_f}$ that is homogenous of degree $0$. Therefore, the invariant takes value in $\bQ$.

In \cite{FLT12}, the open Gromov-Witten invariant was computed via localization and related to a descendant Gromov-Witten invariant of $\cX$. We now summarize the computation. For $i = 1, \dots, s$, based on \eqref{eqn:sigmaWts}, set
$$
	w_{1,i} := \frac{1}{\fr_i},\quad 
	w_{2,i}:=\frac{\fs_i+\fr_i f_i}{\fr_i\fm_i},\quad 
	w_{3,i}:=-\frac{\fm_i+\fs_i+\fr_i f_i}{\fr_i\fm_i}.
$$
Moreover, let
$$
	h(d_i, \lambda_i) := \pi_{(\tau_i, \sigma_i)}(d_i, \lambda_i) \in G_{\sigma_i}
$$
(see \eqref{eqn:FlagFundGroup}) and 
\begin{equation}\label{eqn:AdditionalInsertions}
	\gamma_{n+i} := 
	\iota_{\sigma_i, *}(\one_{h(d_i, \lambda_i)^{-1}}) \quad \in H^*_{\CR, T'}(\cX; \bQ)
\end{equation}
denote the $T'$-equivariant Poincar\'e dual of $(\fp_{\sigma_i}, h(d_i, \lambda_i)^{-1})$. Let $\epsilon_{2,i}, \epsilon_{3,i} \in \bQ \cap [0,1)$ such that $h(d_i, \lambda_i)$ acts on $T_{\fp_{\sigma_i}}\fl_{\tau_{2,i}}, T_{\fp_{\sigma_i}}\fl_{\tau_{3,i}}$ by multiplication by $e^{2\pi\sqrt{-1}\epsilon_{2,i}}, e^{2\pi\sqrt{-1}\epsilon_{3,i}}$ respectively. We define the disk factor (following the sign convention of \cite[Section 3.4.2]{LY22})
\begin{equation}\label{equ:diskinvariant}
	D_{d_i,\lambda_i} = (-1)^{\floor{d_i w_{3,i}-\epsilon_{3,i}} + \ceil{\frac{d_i}{\fa_i}}} \left(\frac{\su_{1,i}}{d_i}\right)^{\age(h(d_i,\lambda_i))-1}\cdot{\frac{1}{d_i\fm_i\lfloor d_iw_{1,i}\rfloor!}}\cdot \prod_{a=1}^{\lfloor d_i w_{1,i}\rfloor+\age(h(d_i,\lambda_i))-1} \left(\frac{d_i\sw_{2,i}}{\su_{1,i}}+a-\epsilon_{2,i} \right).
\end{equation}
Then we have
\begin{align*}
	\langle \gamma_1,\dots,\gamma_n  \rangle^{\cX, \cL, T_f}_{\beta', \bd} & = \prod_{i = 1}^s \text{$\fr_i \fm_i$} D_{d_i,\lambda_i} \left\langle \gamma_1,\dots,\gamma_n , \frac{\gamma_{n+1}}{\frac{\su_{1,1}}{d_1} - \bpsi_{n+1}}, \dots, \frac{\gamma_{n+s}}{\frac{\su_{1,s}}{d_s} - \bpsi_{n+s}}  \right\rangle^{\cX, T'}_{\beta} \bigg|_{\su_2 - f\su_1 = 0} \\
	& = \prod_{i = 1}^s \text{$\fr_i \fm_i$} D_{d_i,\lambda_i} \int_{[\Mbar_{0,n+s}(\cX, \beta)^{T'}]^\vir} \frac{\iota^*\left(\prod_{i = 1}^n \ev_i^*(\gamma_i) \prod_{i = 1}^s \ev_{n+i}^*(\gamma_{n+i}) \right)}{e_{T'}(N^\vir) \prod_{i = 1}^s \left(\frac{\su_{1,i}}{d_i} - \bpsi_{n+i} \right)} \bigg|_{\su_2 - f\su_1 = 0}.
\end{align*}
The localization computation is given as follows, using the description of the $T'$-fixed locus of the moduli space $\Mbar_{0,n+s}(\cX, \beta)$ in terms of decorated graphs given in Section \ref{sect:Moduli}.


\begin{theorem}[\cite{FLT12}]
	We have
    \begin{equation}\label{equ:OpenInvariant}
    \begin{aligned}
      \langle \gamma_1,\dots,\gamma_n  \rangle^{\cX, \cL, T_f}_{\beta', \bd}
      = & \prod_{i=1}^{s}\fr_i\fm_i D_{d_i,\lambda_i} \sum_{\vGa\in \Gamma_{0,n+s}(\cX,\beta)} c_{\vGa} \prod_{e\in E(\Gamma)} \bh(e) \prod_{(e,v)\in F(\Gamma)} \bh(e,v) \\
    &\prod_{v\in V(\Gamma)}\int_{\Mbar_{0,\vk_v}(\cB G_v)}\frac{\bh(v) \prod_{i \in S_v}\iota^*_{\sigma_v} (\gamma_i) }{\prod_{n+i \in S_v} (\frac{\su_{1,i}}{d_i} - \bpsi_{n+i}) \prod_{e\in E_v}\left(\bw_{(e,v)}-\frac{\bpsi_{(e,v)}}{r_{(e,v)}}\right)} \bigg|_{\su_2-f\su_1=0}.
    \end{aligned}
    \end{equation}
\end{theorem}

Here, to define the quantities involved in the theorem, for each $\vGa\in\Gamma_{0,n+s}(\cX,\beta)$, we pick a stable map $u:(\cC,\fx_1,\dots,\fx_{n+s})\to\cX$ whose associated decorated graph is $\vGa$ (see Section \ref{sect:Moduli}). The definitions do not depend on the choice of $u$.
\begin{itemize}
    \item For each $e\in E(\Gamma)$, define $$\bh(e)=\frac{e_{T'}(H^1(\cC_e,(u|_{\cC_e})^*T\cX)^m)}{e_{T'}(H^0(\cC_e,(u|_{\cC_e})^*T\cX)^m)}$$
    where the superscript ``$m$'' represents the \emph{moving part}: Any complex representation $V$ of a torus $(\mathbb{C}^*)^r$ can be decomposed as $V=V^f\oplus V^m$ where the fixed part $V^f$ is the direct sum of all trivial 1-dimensional subrepresentations and the moving part $V^m$ is the direct sum of all nontrivial ones.

    \item For each $(e,v)\in F(\Gamma)$, define
	$$
		\bh(e,v):=e_{T'}((T_{\fp_{\sigma_v}}\cX)^{k_{(e,v)}})=\prod_{\substack{(\tau,\sigma_v)\in F(\Sigma), k_{(e,v)}\in G_{\tau}}}\sw(\tau,\sigma_v)
	$$
    where $(T_{\fp_{\sigma_v}}\cX)^{k_{(e,v)}}$ is the maximal subspace of $T_{\fp_{\sigma_v}}\cX$ that is invariant under the action of $k_{(e,v)}$.

	\item For each stable vertex $v \in V^S(\vGa)$, define
	$$
		\bh(v):= \frac{e_{T'}(H^1(\cC_v,(u|_{\cC_v})^*T\cX)^m)}{e_{T'}(H^0(\cC_v,(u|_{\cC_v})^*T\cX)^m)},
	$$
	and for each unstable vertex $v \not \in V^S(\hGa)$, define
	$$
		\bh(v):= \bh(e,v)^{-1} 
	$$
	for $e \in E_v$. 

    \item For each $v\in V(\Gamma)$, the marked points and corresponding descendant classes of $\Mbar_{0,\vk_v}(\cB G_v)$ are indexed by $E_v\cup S_v$. We take the following integration convention for unstable vertices
    $$
		\int_{\Mbar_{0,(1)}(\cB G)} \frac{1}{\bw_1 - \bpsi_1} = \frac{\bw_1}{|G|}, \qquad 	\int_{\Mbar_{0,(k, k^{-1})}(\cB G)} \frac{1}{\bw_1 - \bpsi_1} = \frac{1}{|G|},
	$$
	$$
		\int_{\Mbar_{0,(k, k^{-1})}(\cB G)} \frac{1}{(\bw_1 - \bpsi_1)(\bw_2 - \bpsi_2)} = \frac{1}{|G|(\bw_1 + \bw_2)}.
	$$
    
    \item For each $(e,v)\in F(\Gamma)$, define
    $$\mathbf{w}_{(
    e,v)}:=e_{T'}(T_{\fn(e,v)}\cC_e)=\frac{\fr(\tau_e,\sigma_v)\sw(\tau_e,\sigma_v)}{r_{(e,v)}d_e},$$
    where the $T'$-action on $T_{\fn(e,v)}\cC_e$ is induced from that on $T_{\fp_{\sigma_v}}\fl_{\tau_e}.$
\end{itemize}

\subsection{Corresponding data}\label{sect:CorrData}
Under the setup of open invariants in Section \ref{sect:OpenGW}, we define the corresponding data on the relative and local sides of the correspondences. 

\subsubsection{Corresponding curve classes}
Let $\tX$ be the coarse moduli space of $\tcX$. Recall that
$$
    \tSi(2+s)_c = \iota(\Sigma(2)_c \sqcup \{\tau_1, \dots, \tau_s\}).
$$
The inclusion $\iota: \cX \to \tcX$ induces an isomorphism
$$
    H_2(X,L; \bZ) \to H_2(\tX; \bZ), \qquad [B_i] \mapsto [l_{\tau_i}]
$$
under which $\Eff(X, L)$ is identified with $\Eff(\tX)$.

Let $\beta, d_1, \dots, d_s$, and $\beta'$ be as in Section \ref{sect:OpenGW}. Then we define
$$
	\tbeta \in \Eff(\tX)
$$
to be the class corresponding to $\beta' \in \Eff(X, L)$.

For the formal intermediate geometries, recall that the sets of compact edges
$$
	E(\hGa^{(0)})_c = \cdots = E(\hGa^{(s)})_c
$$
in the FTCY graphs $\sGa^{(0)}, \dots, \sGa^{(s)}$ are all identified with $\tSi(2+s)_c$. For $\ell = 0, \dots, s$, let $\hY^{(\ell)}$ denote the coarse moduli space of $\hcY^{(\ell)}$. In particular, $\hY = \hY^{(0)}$ denotes the coarse moduli space of $\hcY = \hcY^{(0)}$. There is then an identification of effective classes
$$
	\Eff(\hY^{(0)}) = \cdots = \Eff(\hY^{(s)})
$$
and a commutative diagram
$$
	\xymatrix{
		\Eff(\hY^{(0)}) \ar@{=}[r] \ar[d]^\pi & \Eff(\hY^{(s)}) \ar[d]^\pi \\
		\Eff(X, L) \ar[r]^\sim & \Eff(\tX).
	}
$$
By an abuse of notation, we will use $\hbeta$ to denote a class in $\Eff(\hY^{(0)}) = \cdots = \Eff(\hY^{(s)})$ that maps to the corresponding classes $\beta', \tbeta$ under the projections $\pi$ in the above.

\subsubsection{Corresponding insertions}\label{sect:Insertions}
Let $\ell = 0, \dots, s$. Recall that the set of vertices in the FTCY graph $\sGa^{(\ell)}$ is
$$
	V(\hGa^{(\ell)}) = \Sigma(3) \sqcup \{\tv_1, \dots, \tv_s\}.
$$

Let $\gamma_1, \dots, \gamma_n \in H^2_{\CR, T'}(\cX; \bQ)$. For $i = 1, \dots, n$, we take a lift $\gamma_i^{(\ell)} \in \cH(\hcY^{(\ell)})$ of $\gamma_i$ in the following way. Let $\iota_v^*(-)$ denotes taking the component in \eqref{eqn:StateSpace} indexed by $v \in V(\hGa^{(\ell)})$. Given $\gamma \in H^2_{\CR, T'}(\cX; \bQ)$, if $\gamma \in H^2_{T'}(\cX; \bQ)$, we take the lift $\gamma^{(\ell)} \in \cH(\hcY^{(\ell)})$ determined by that 
$$
	\iota_{\sigma}^*(\gamma^{(\ell)}) \big|_{\su_4 = \cdots \su_{3+\ell} = 0} = \iota_\sigma^*(\gamma)
$$
for all $\sigma \in \Sigma(3)$ and $\iota_{\tv_j}^*(\gamma^{(\ell)}) = 0$ for $j = 1, \dots, s$. On the other hand, if $\gamma = \one_j$ for some twisted sector $j \in \Box(\cX)$, we take the lift $\gamma^{(\ell)}$ determined by that
$$
	\iota_v^*(\gamma^{(\ell)}) = \begin{cases}
		\one_j & \text{if } j \in \Box(v),\\
		0 & \text{if } j \not \in \Box(v).
	\end{cases}
$$


Let $(d_1, \lambda_1), \dots, (d_s, \lambda_s)$ be as in Section \ref{sect:OpenGW}. For $i = 1, \dots, s$, we define the element
$$
	\tk_i := \pi_{(\iota(\tau_i), \tsi_i)}(d_i, \lambda_i) \in G_{\tsi_i}
$$
with respect to the closed geometry $\tcX$. Its $\age$ in $\tcX$ is described by Lemma \ref{lem:Age}.

Now for $(\hcY^{(\ell)}, \hcD^{(\ell)})$, if $i \le \ell$, the stabilizer group of the chart $\hcY^{(\ell)}_{\tv_i}$ is identified with $G_{\tsi_i}$ in a way that preserves $\age$. We define the insertion
$$
	\gamma_{n+i}^{(\ell)} \in \cH(\hcY^{(\ell)})
$$
such that
$$
	\iota_{\tv_i}^*(\gamma_{n+i}^{(\ell)})
	= \begin{cases}
		\sfp^{(\ell)}(\iota(\tau_i), \tsi_i) \sff_2^{(\ell)}(\tsi_i) & \text{if } \tk_i = 1,\\
		\sfp^{(\ell)}(\iota(\tau_i), \tsi_i) \one_{\tk_i^{-1}} & \text{if } \tk_i \in G_{\tau_i} \setminus \{1\},\\
		\one_{\tk_i^{-1}} & \text{if } \tk_i \not \in G_{\tau_i},
	\end{cases}
$$
and $\iota_{v}^*(\gamma_{n+i}^{(\ell)}) = 0$ for any $v \neq \tv_i$. This is a homogenous element of degree $4$.

On the other hand, if $i \ge \ell+1$, the stabilizer group of the divisor component $\hcD_i^{(\ell)}$ is also identified with $G_{\tsi_i}$ except that for $\tk \in G_{\tsi_i}$, we have $\age_{\hcD}(\tk) \in \{0, 1\}$ and it is equal to 0 if and only if $\tk = 1$. We define the insertion
$$
	\gamma_{n+i}^{(\ell)} := \begin{cases}
		\sff_2^{(\ell)}(\tsi_i) & \text{if } \tk_i = 1,\\
		\one_{\tk_i^{-1}} & \text{if } \tk_i \neq 1
	\end{cases} \quad \in \cH(\hcD_i^{(\ell)}),
$$
which is a homogenous element of degree $2$. Let
$$
	\bk^{(\ell)} := (\tk_{\ell+1}, \dots, \tk_s).
$$
For $\ell = 0$ we will write $\bk := \bk^{(0)}$. Note also that $\bk^{(s)}$ is empty.

When $\ell = s$, where $\hcY^{(s)}$ is the formal completion of $\tcX$ along $\tcX^1$, the elements $\gamma_1^{(s)}, \dots, \gamma_{n+s}^{(s)}$ defined above lie in the image of the injective map
$$
    H_{\CR, \tT'}^*(\tcX; \bQ) \to \cH(\hcY^{(s)}).
$$
We let
$$
	\tgamma_1, \dots, \tgamma_{n+s} \in H^*_{\CR, \tT'}(\tcX; \bQ)
$$
be the corresponding elements. The elements $\tgamma_{n+1}, \dots, \tgamma_{n+s}$ lie in $H^4_{\CR, \tT'}(\tcX; \bZ)$ and can be described alternatively as follows. For $i = 1, \dots, R+2s$, let $[\cV(\trho_i)] \in H^2_{\tT'}(\tcX; \bZ)$ denote the $\tT'$-equivariant Poincar\'e dual of the divisor $\cV(\trho_i)$. Then, for $i = 1, \dots, s$, we have
$$
	\tgamma_{n+i} = \begin{cases}
		[\cV(\trho_{R+2i-1})] [\cV(\trho_{i_2^{(\tau_i,\sigma_i)}})] & \text{if } \tk_i = 1,\\
		[\cV(\trho_{R+2i-1})] \one_{\tk_i^{-1}} & \text{if } \tk_i \in G_{\tau_i} \setminus \{1\},\\
		\one_{\tk_i^{-1}} & \text{if } \tk_i \not \in G_{\tau_i}.
	\end{cases}
$$
It admits a natural non-equivariant limit in $H^4_{\CR}(\tcX; \bZ)$ defined by the same formula, with the notation $[\cV(\trho_i)]$ used to represent the non-equivariant divisor classes.

\subsection{Formal relative Gromov-Witten invariants}\label{sec:relativeGW}
Let $\ell = 0, \dots, s$. Now we consider the relevant Gromov-Witten invariant of the relative FTCY orbifold $(\hcY^{(\ell)},\hcD^{(\ell)})$.

\begin{definition}\label{def:RelGW}
Let $\bd, \beta'$ be as in Section \ref{sect:OpenGW} and $\gamma_1, \dots, \gamma_n \in H^2_{\CR, T'}(\cX; \bQ)$. Let $\hbeta \in \Eff(\hY^{(\ell)})$ such that $\pi(\hbeta) = \beta'$, and $\bk^{(\ell)}$, $\gamma_1^{(\ell)}, \dots, \gamma_{n+s}^{(\ell)}$ be defined as in Section \ref{sect:Insertions}. We define the formal relative Gromov-Witten invariant (cf. \eqref{eqn:RGWDef})
\begin{align*}
	& \langle \gamma_1^{(\ell)},\dots,\gamma_{n+\ell}^{(\ell)} \mid \gamma_{n+\ell+1}^{(\ell)},\dots,\gamma_{n+s}^{(\ell)} \rangle^{\hcY^{(\ell)}/\hcD^{(\ell)}, T_f}_{\hbeta,\bk^{(\ell)}} \\
	& \qquad := \langle \gamma_1^{(\ell)},\dots,\gamma_{n+\ell}^{(\ell)} \mid \gamma_{n+\ell+1}^{(\ell)},\dots,\gamma_{n+s}^{(\ell)} \rangle^{\hcY^{(\ell)}/\hcD^{(\ell)}, T'_{(\ell)}}_{\hbeta,\bk^{(\ell)}}  \bigg|_{\su_2 - f\su_1 = \su_4 = \cdots = \su_{3+\ell} = 0} \quad \in \bQ.
\end{align*}
\end{definition}
	
We now give the localization computation of the formal relative invariant, following computations in \cite{FLT12,Zong12,Zong15}. We use the description of the $T'_{(\ell)}$-fixed locus of the moduli space $\Mbar_{0,n+\ell}(\hcY^{(\ell)}/\hcD^{(\ell)}, \hbeta, \bk^{(\ell)})$ in terms of decorated graphs given in Section \ref{sect:Moduli}.
	

\begin{proposition}\label{prop:RelativeLocalization}
We have	
\begin{equation}\label{eqn:RelativeLocalization}
\begin{aligned}
	& \langle \gamma_1^{(\ell)},\dots,\gamma_{n+\ell}^{(\ell)} \mid \gamma_{n+\ell+1}^{(\ell)},\dots,\gamma_{n+s}^{(\ell)} \rangle^{\hcY^{(\ell)}/\hcD^{(\ell)}, T_f}_{\hbeta,\bk^{(\ell)}} \\
	& = \sum_{\vGa \in \Gamma_{0,n+\ell}(\hcY^{(\ell)}/\hcD^{(\ell)},\hbeta,\bk^{(\ell)})}  c_{\vGa} \prod_{e\in E(\Gamma)}\bh^{(\ell)}(e) \prod_{i=\ell+1}^s A_v^i \cdot \iota_{\tsi_i}^*(\gamma_{n+i}^{(\ell)})\\
	& \prod_{v\in V(\vGa)^{(0)}} \prod_{(e,v)\in F(\Gamma)}   \bh^{(\ell)}(e,v)\int_{\Mbar_{0,\vk_v}(\cB G_v)}\frac{\bh^{(\ell)}(v) \prod_{i\in S_v}\iota^*_{\sigma_v}(\gamma_i^{(\ell)}) }{\prod_{e\in E_v}(\bw^{(\ell)}_{(e,v)}-\frac{\overline{\psi}_{(e,v)}}{r_{(e,v)}})} \bigg|_{\su_2 - f\su_1 = \su_4 = \cdots = \su_{3+\ell} = 0} .
\end{aligned}
\end{equation}
\end{proposition}

Here, to define the quantities involved in the proposition, for each $\vGa \in \Gamma_{0,n+\ell}(\hcY^{(\ell)}/\hcD^{(\ell)},\hbeta,\bk^{(\ell)})$, we pick a stable map $u:(\cC,\fx_1,\dots,\fx_{n+\ell}, \fy_{n+\ell+1}, \dots, \fy_s) \to (\hcY^{(\ell)}_{\bm}, \hcD^{(\ell)}_{\bm})$ whose associated decorated graph is $\vGa$ (see Section \ref{sect:Moduli}). The definitions do not depend on the choice of $u$.
\begin{itemize}
	\item For each $e\in E(\hGa)$, define 
	$$
		\bh^{(\ell)}(e)=\frac{e_{T'_{(\ell)}}(H^1(\cC_e,(u|_{\cC_e})^*\Omega_{\hcY^{(\ell)}}(\log(\hcD^{(\ell)}))^\vee)^m)}{e_{T'_{(\ell)}}(H^0(\cC_e,(u|_{\cC_e})^*\Omega_{\hcY^{(\ell)}}(\log(\hcD^{(\ell)}))^\vee)^m)}.
	$$

	\item For each $(e,v)\in F(\Gamma)$, $v\in V(\Gamma)^{(0)}$, define
    \begin{align*}
        \bh^{(\ell)}(e,v):=\prod_{\substack{(\tau,\sigma_v)\in F(\hGa), k_{(e,v)}\in G_{\tau}}}\sfp^{(\ell)}(\tau,\sigma_v),
    \end{align*}
		
	

	\item For each $v \in V(\vGa)^{(0)} \cap V^S(\hGa)$, define
	$$
	\bh^{(\ell)}(v):= \frac{e_{T'_{(\ell)}}(H^1(\cC_v,(u|_{\cC_v})^*\Omega_{\hcY^{(\ell)}}(\log(\hcD^{(\ell)}))^\vee)^m)}{e_{T'_{(\ell)}}(H^0(\cC_v,(u|_{\cC_v})^*\Omega_{\hcY^{(\ell)}}(\log(\hcD^{(\ell)}))^\vee)^m)},
	$$
	and for each $v \in V(\vGa)^{(0)} \setminus V^S(\hGa)$, define
	$$
		\bh^{(\ell)}(v):= \bh^{(\ell)}(e,v)^{-1} 
	$$
	for $e \in E_v$. 

	\item For each $(e,v)\in F(\Gamma)$, define
	$$
		\bw^{(\ell)}_{(e,v)} := \frac{\fr(\tau_e,\sigma_v)\sfp^{(\ell)}(\tau_e,\sigma_v)}{r_{(e,v)}d_e}.
	$$
	
	
   \item For $i = \ell+1, \dots, s$, define $A^i_v = 1$ if $l(\mu_i(\vGa))=1$, and 
   \begin{align}\label{eqn: Avi}
       	A^i_v = \prod_{j=1}^{l(\mu_i(\vGa))}(\mu_i(\vGa))_j\int_{{\Mbar}^{(i)}_{\vGa}}\frac{1}{-\frac{\su_{1,i}}{\fa_i}-\psi_i}\left(\prod_{\substack{-1\le j \le \ell\\k_{(e,\tv_i)}\in G_{\delta_{3+j}(\tsi_i
        )}}}\sff^{(\ell)}_{j+3}(\tsi_i)\right)
        \frac{\prod_{-1\le j\le \ell }\Lambda^{\vee}_{\chi(\delta_{3+j}(\tsi_i),\tsi_i)}(\sff^{(\ell)}_{j+3}(\tsi_i))}{\prod_{-1\le j\le\ell, \hat{G}_i\subset G_{\delta_{3+j}(\tsi_i)}} \sff^{(\ell)}_{j+3}(\tsi_i)}
   \end{align}
	if $l(\mu_i(\vGa))>1$, where 
\begin{itemize}[label=$\circ$]
    \item $\psi_i$ is the target psi class on $\Mbar_{\vGa}^{(i)}:=\Mbar_{0,0}(\bP^1 \times G_i, \mu_i(\vGa), (d_i)) \sslash \bC^*$.
    

    \item $\hat{G}_i$ denote the subgroup of $G_{\tv_i}$ generated by the monodromies of the $G_{\tv_i}$-cover $\widetilde{\cC} \to \cC$, pulled back from $\pt \to \cB G_{\tv_i}$ along a morphism $\cC \to \fp_{\tsi_i}$ that represents a point in $\Mbar_{0, \vk_{\tv_i}}(\cB G_{\tv_i})$. Note that in this case, $\tv_i$ is a stable vertex.
    
\end{itemize}
	 
\end{itemize}



\subsection{Closed Gromov-Witten invariants}\label{sec: LocalGW}
Finally, we consider the relevant Gromov-Witten invariant of the closed geometry $\tcX$.

\begin{definition}\label{def:ClosedGW}
Let $\bd, \beta'$ be as in Section \ref{sect:OpenGW} and $\gamma_1, \dots, \gamma_n \in H^2_{\CR, T'}(\cX; \bQ)$. Let $\tbeta \in \Eff(\tX)$ and $\tgamma_1, \dots, \tgamma_{n+s}$ be defined as in Section \ref{sect:CorrData}. We define the Gromov-Witten invariant
$$
	\langle \tgamma_1, \dots, \tgamma_{n+s} \rangle^{\tcX, T_f}_{\tbeta}
	:= \langle \tgamma_1, \dots, \tgamma_{n+s} \rangle^{\tcX, \tT'}_{\tbeta}  \bigg|_{\su_2 - f\su_1 = \su_4 = \cdots = \su_{3+s} = 0} \qquad \in \bQ.
$$
\end{definition}

As in Section \ref{sect:Absolute}, the closed invariant is related to the formal invariants of $(\hcY^{(s)},\hcD^{(s)} = \emptyset)$ by
$$
    \langle \tgamma_1, \dots, \tgamma_{n+s} \rangle^{\tcX, T_f}_{\tbeta} = \sum_{\hbeta \in \Eff(\hY^{(s)}), \pi(\hbeta) = \tbeta} \langle \gamma_1^{(s)},\dots, \gamma_{n+s}^{(s)} \rangle^{\hcY^{(s)}, T_f}_{\hbeta}.
$$
Then, Proposition \ref{prop:RelativeLocalization} applied with $\ell = s$ provides a localization computation of the closed invariant.

\section{The correspondences}\label{sect:Correspondence}
In this section, we prove the correspondences among the open, relative, and local Gromov-Witten invariants. Let $\bd, \beta'$ be as in Section \ref{sect:OpenGW} and $\gamma_1, \dots, \gamma_n \in H^2_{\CR, T'}(\cX; \bQ)$.

\begin{theorem}\label{thm:OpenClosedStatement}
	We have the open/relative/local correspondence
	$$
		\langle \gamma_1,\dots,\gamma_n  \rangle^{\cX, \cL, T_f}_{\beta',\bd} 
		= (-1)^{\sum_{i = 1}^s (\ceil{\frac{d_i}{\fa_i}} -1)} \sum_{\hbeta \in \Eff(\hY), \pi(\hbeta) = \beta'}\langle \gamma_1^{(0)},\dots,\gamma_n^{(0)} \mid \gamma_{n+1}^{(0)},\dots,\gamma_{n+s}^{(0)} \rangle^{\hcY/\hcD, T_f}_{\hbeta,\bk} 
		= \langle \tgamma_1, \dots, \tgamma_{n+s} \rangle^{\tcX, T_f}_{\tbeta}
	$$
	where the data $\bk$, $\gamma_1^{(0)}, \dots, \gamma_{n+s}^{(0)}$, $\tbeta$, $\tgamma_1, \dots, \tgamma_{n+s}$ are defined in Section \ref{sect:CorrData}.
\end{theorem}

Theorem \ref{thm:OpenClosedStatement} is based on a correspondence among the formal relative Gromov-Witten invariants of the intermediate geometries, stated inductively as follows.


\begin{theorem}\label{thm:RelativeStatement}

	Let $\ell = 0, \dots, s-1$, and $\hbeta \in \Eff(\hY^{(\ell)}) = \Eff(\hY^{(\ell+1)})$ such that $\pi(\hbeta) = \beta'$. We have the correspondence
	$$
		\langle \gamma_1^{(\ell)},\dots,\gamma_{n+\ell}^{(\ell)} \mid \gamma_{n+\ell+1}^{(\ell)},\dots,\gamma_{n+s}^{(\ell)} \rangle^{\hcY^{(\ell)}/\hcD^{(\ell)}, T_f}_{\hbeta,\bk^{(\ell)}} 
		= (-1)^{\ceil{\frac{d_{\ell+1}}{\fa_{\ell+1}}} -1} \langle \gamma_1^{(\ell+1)},\dots,\gamma_{n+\ell+1}^{(\ell+1)} \mid \gamma_{n+\ell+2}^{(\ell+1)},\dots,\gamma_{n+s}^{(\ell+1)} \rangle^{\hcY^{(\ell+1)}/\hcD^{(\ell+1)}, T_f}_{\hbeta,\bk^{(\ell+1)}}
	$$
	where the data $\bk^{(\ell)}$, $\bk^{(\ell+1)}$, $\gamma_1^{(\ell)}, \dots, \gamma_{n+s}^{(\ell)}$, $\gamma_1^{(\ell+1)}, \dots, \gamma_{n+s}^{(\ell+1)}$ are defined in Section \ref{sect:CorrData}.
\end{theorem}

\subsection{Contributions of external edges}
We first give a more detailed analysis of the formal relative Gromov-Witten invariant
$$
	\langle \gamma_1^{(\ell)},\dots,\gamma_{n+\ell}^{(\ell)} \mid \gamma_{n+\ell+1}^{(\ell)},\dots,\gamma_{n+s}^{(\ell)} \rangle^{\hcY^{(\ell)}/\hcD^{(\ell)}, T_f}_{\hbeta,\bk^{(\ell)}}
$$
of $(\hcY^{(\ell)},\hcD^{(\ell)})$, for $\ell = 0, \dots, s$, as computed by localization in Proposition \ref{prop:RelativeLocalization}. In this subsection, we explicitly compute the contributions of the edges whose labels are among $\iota(\tau_1), \dots,\iota(\tau_s)$.

\begin{lemma}\label{lem:WtExternalEdge}
In \eqref{eqn:RelativeLocalization}, let $e\in E(\Gamma)$ be an edge with label $\iota(\tau_i)$ for some $1 \le i \le s$ and degree $(d, \lambda) \in \bZ_{\ge 1} \times G_{\tau_i}$. Let $h(d, \lambda) \in G_{\sigma_i}$ and $\epsilon_2, \epsilon_3 \in \bQ \cap [0,1)$ be defined as in Section \ref{sect:OpenGW} and $\tk := \pi_{(\iota(\tau_i), \tsi_i)}(d, \lambda) \in G_{\tsi_i}$.
\begin{itemize}[leftmargin=*]

	\item Suppose $1\le i \le \ell$. Then
	$$
		\bh^{(\ell)}(e) = \bh' \prod_{1 \le j \le \ell, j \neq i} \su_{3+j}^{-1}
	$$
	where none of $\su_4, \dots, \su_{3+\ell}$ is a factor of $\bh'$. Moreover, if $\tk=1$, we have
	\begin{align*}
		&\bh' \big|_{\su_4 = \dots =\su_{3+\ell}= 0} \\
		&	= \frac{(-1)^{\floor{dw_{3,i}-\epsilon_{3,i}}+ \frac{d}{\fa_i} + 1}}{\floor{dw_{1,i}}!} \left(\frac{\su_{1,i}}{d} \right)^{\age(h(d,\lambda)) - 1} 
		\left(\frac{\su_{1,i}}{\fa_i}\right)^{-1}
		\left( \frac{\su_{2,i}-f_i\su_{1,i}}{\fm_i} \right)^{-1} \cdot \prod_{a = 1}^{\floor{dw_{1,i}}+\age(h(d,\lambda))-1} \left(\frac{d\sw_{2,i}}{\su_{1,i}} + a - \epsilon_{2,i} \right).
	\end{align*}
	If $\age(\tk)=1$, we have
	$$
		\bh' \big|_{\su_4 = \dots =\su_{3+\ell}= 0} 
		= \frac{(-1)^{\floor{dw_{3,i}-\epsilon_{3,i}}+ \frac{d}{\fa_i} + 1}}{\floor{dw_{1,i}}!} \left(\frac{\su_{1,i}}{d} \right)^{\age(h(d,\lambda)) - 1} \left(\frac{\su_{1,i}}{\fa_i}\right)^{-1} \cdot \prod_{a = 1}^{\floor{dw_{1,i}}+\age(h(d,\lambda))-1} \left(\frac{d\sw_{2,i}}{\su_{1,i}} + a - \epsilon_{2,i} \right).
	$$
	If $\age(\tk)=2$, we have
	$$
		\bh' \big|_{\su_4 = \dots =\su_{3+\ell}= 0} 
		= \frac{(-1)^{\floor{dw_{3,i}-\epsilon_{3,i}}+ \floor{\frac{d}{\fa_i}} + 1}}{\floor{dw_{1,i}}!} \left(\frac{\su_{1,i}}{d} \right)^{\age(h(d,\lambda)) - 1}  \cdot \prod_{a = 1}^{\floor{dw_{1,i}}+\age(h(d,\lambda))-1} \left(\frac{d\sw_{2,i}}{\su_{1,i}} + a - \epsilon_{2,i} \right).
	$$

	\item  Suppose otherwise $\ell+1 \le i \le s$. Then
	$$
		\bh^{(\ell)}(e) = \bh' \prod_{1 \le j \le \ell} \su_{3+j}^{-1}
	$$
	where none of $\su_4, \dots, \su_{3+\ell}$ is a factor of $\bh'$. Moreover, if $\tk=1$, we have
	\begin{align*}
		&\bh' \big|_{\su_4 = \dots =\su_{3+\ell}= 0} \\
		&= \frac{(-1)^{\floor{dw_{3,i}-\epsilon_{3,i}}+ 
				1}}{\floor{dw_{1,i}}!} \left(\frac{\su_{1,i}}{d} \right)^{\age(h(d,\lambda)) - 1} 
		\left(\frac{\su_{2,i}-f_i\su_{1,i}}{\fm_i}\right)^{-1}
		\cdot\prod_{a = 1}^{\floor{dw_{1,i}}+\age(h(d,\lambda))-1} \left(\frac{d\sw_{2,i}}{\su_{1,i}} + a - \epsilon_{2,i} \right).
	\end{align*} 
	If $\tk\neq 1$, we have
	$$
		\bh' \big|_{\su_4 = \dots =\su_{3+\ell}= 0} 
		= \frac{(-1)^{\floor{dw_{3,i}-\epsilon_{3,i}}  + 1}}{\floor{dw_{1,i}}!} \left(\frac{\su_{1,i}}{d} \right)^{\age(h(d,\lambda)) - 1} 
	\cdot\prod_{a = 1}^{\floor{dw_{1,i}}+\age(h(d,\lambda))-1} \left(\frac{d\sw_{2,i}}{\su_{1,i}} + a - \epsilon_{2,i} \right).
	$$
	
\end{itemize}
\end{lemma}

\begin{proof}
We refer to the explicit computation of $\bh^{(\ell)}(e)$ in \cite[Lemma 130]{Liu13} and a similar computation in \cite[Lemma B.10]{LY22}.  
We have a splitting
$$
	\Omega_{\hcY^{(\ell)}}(\log(\hcD^{(\ell)}))^\vee \big|_{\fl_{\iota(\tau_i)}} = 
	\Omega_{\fl_{\iota(\tau_i)}}(\log(\fl_{\iota(\tau_i)} \cap \hcD^{(\ell)}))^\vee \oplus N_{\fl_{\iota(\tau_i)}/\hcY^{(\ell)}}
$$
and the normal bundle $N_{\fl_{\iota(\tau_i)}/\hcY^{(\ell)}}$ further splits as a direct sum of $T'_{(\ell)}$-equivariant line bundles $L_2, L_3, \dots, L_{3+\ell}$ corresponding to the 2-dimensional torus closed substacks of $\hcY^{(\ell)}$ that contain $\fl_{\iota(\tau_i)}$.
We set
$$
	\bb_1 := \frac{1}{e_{T'_{(\ell)}}(H^0(\cC_e,(u|_{\cC_e})^*\Omega_{\fl_{\iota(\tau_i)}}(\log(\fl_{\iota(\tau_i)} \cap \hcD^{(\ell)}))^\vee)^m)},
$$
$$
	\bb_j := \frac{e_{T'_{(\ell)}}(H^1(\cC_e,(u|_{\cC_e})^*L_j)^m)}{e_{T'_{(\ell)}}(H^0(\cC_e,(u|_{\cC_e})^*L_j)^m)}, \quad j = 2, \dots, 3+\ell.
$$
Then, we have
$$
	\bh^{(\ell)}(e) =  \bb_1 \bb_2 \bb_3 \cdots \bb_{3+\ell}.
$$

First suppose $1\le i \le \ell$. In this case, $\Omega_{\fl_{\iota(\tau_i)}}(\log(\fl_{\iota(\tau_i)} \cap \hcD^{(\ell)}))^\vee$ is the tangent bundle $T_{{\iota(\tau_i)}}$. Moreover, for $j = 1, \dots, \ell$, $j \neq i$, $L_{3+j}$ is a trivial line bundle with weight $\su_{3+j}$ and $\bb_{3+j} = \su_{3+j}^{-1}$. Then
$$
	\bh' = \bb_1\bb_2\bb_3\bb_{3+i}
$$
and the lemma statement follows from the computation
$$
	\bb_1 =\frac{(-1)^{\floor{\frac{d}{\fa_i}}}}{\floor{d w_{1,i}}!\floor{\frac{d}{\fa_i}}!}
	\left(\frac{\text{$\su_{1,i}$}}{d}\right)^{-\floor{d w_{1,i}}-\floor{\frac{d}{\fa_i}}}, \quad 
	\bb_{3+i}\big{|}_{u_4=\cdots=u_{3+\ell}=0} = \left(\frac{\su_{1,i}}{d}\right)^{\ceil{\frac{d}{\fa_i}-1}} \ceil{\frac{d}{\fa_i}-1}!,
$$
$$
	\bb_2 = \begin{cases}
		\prod_{a=0}^{\floor{dw_{2,i}-\epsilon_{2,i}}} \left(\sw_{2,i}- \frac{a+\epsilon_{2,i}}{d}\su_{1,i} \right)^{-1} & \text{if } w_{2,i} \ge 0,\\
		\prod_{a=1}^{\floor{\epsilon_{2,i}-dw_{2,i}-1}} \left(\sw_{2,i}+ \frac{a-\epsilon_{2,i}}{d}\su_{1,i} \right) & \text{if } w_{2,i} < 0, 
	\end{cases}
	\quad 
	\bb_3 = \begin{cases}
		\prod_{a=0}^{\floor{dw_{3,i}-\epsilon_{3,i}}} \left(\sw_{3,i}- \frac{a+\epsilon_{3,i}}{d}\su_{1,i} \right)^{-1} & \text{if } w_{3,i} \ge 0,\\
		\prod_{a=1}^{\floor{\epsilon_3-dw_{3,i}-1}} \left(\sw_{3,i}+ \frac{a-\epsilon_{3,i}}{d}\su_{1,i} \right) & \text{if } w_{3,i} < 0.
	\end{cases}
$$


Suppose otherwise $\ell+1 \le i \le s$. In this case, for $j = 1, \dots, \ell$, $L_{3+j}$ is a trivial line bundle with weight $\su_{3+j}$ and $\bb_{3+j} = \su_{3+j}^{-1}$. Then
$$
	\bh' = \bb_1\bb_2\bb_3
$$
where $\bb_2, \bb_3$ are the same as in the previous case and 
$$
	\bb_1 = \frac{1}{\floor{d w_{1,i}}!} \left(\frac{\su_{1,i}}{d}\right)^{-\floor{d w_{1,i}}}.
$$
\end{proof}

\subsection{Vanishing arguments}\label{sect:Vanishing}
The computation of $\langle \gamma_1^{(\ell)},\dots,\gamma_{n+\ell}^{(\ell)} \mid \gamma_{n+\ell+1}^{(\ell)},\dots,\gamma_{n+s}^{(\ell)} \rangle^{\hcY^{(\ell)}/\hcD^{(\ell)}, T_f}_{\hbeta,\bk^{(\ell)}}$ in Proposition \ref{prop:RelativeLocalization} expresses the invariant as a sum of contributions from decorated graphs in $\Gamma_{0,n+\ell}(\hcY^{(\ell)}/\hcD^{(\ell)},\hbeta,\bk^{(\ell)})$. In this subsection, we show that the contribution of a decorated graph vanishes unless it belongs to the following subset.


\begin{definition}\label{def:ContributingGraphs}
\rm{
Let $\ell = 0, \dots, s$. We define $\Gamma_{0,n+\ell}(\hcY^{(\ell)}/\hcD^{(\ell)},\hbeta,\bk)^{(0)}$ to be the subset of graphs $\vGa$ in $\Gamma_{0,n+\ell}(\hcY^{(\ell)}/\hcD^{(\ell)},\hbeta,\bk^{(\ell)})$ that satisfies the following conditions:
\begin{enumerate}[label=(\alph*)]
	\item \label{item:RelEdge} For $i = \ell+1, \dots, s$, $l(\mu_i(\vGa))=1$.
		
	\item \label{item:LocalEdge} For $i = 1, \dots, \ell$, there exists a unique edge in $E(\Gamma)$ labeled by $\iota(\tau_i)$.
		
	\item \label{item:Insertion} 
	For $i=1,\dots,\ell$, $\vf\circ\vs(n+i)=\tsi_{i}$ and $\vk(n+i)=h(d_{i},\lambda_{i})^{-1}$.

	\item \label{item:LocalVertex}
	For $i=1,\dots,\ell$, there is a unique vertex $\tv_i$ with $\vf(\tv_i)=\tsi_{i}$. In addition, $S_{\tv_i} = \{n+i\}$.
\end{enumerate}
}
\end{definition}

For $\vGa \in \Gamma_{0,n+\ell}(\hcY^{(\ell)}/\hcD^{(\ell)},\hbeta,\bk^{(\ell)})$, let
$
	C_{\vGa}^{(\ell)}
$
denote the contribution of $\vGa$ to $\langle \gamma_1^{(\ell)},\dots,\gamma_{n+\ell}^{(\ell)} \mid \gamma_{n+\ell+1}^{(\ell)},\dots,\gamma_{n+s}^{(\ell)} \rangle^{\hcY^{(\ell)}/\hcD^{(\ell)}, T'_{(\ell)}}_{\hbeta,\bk^{(\ell)}}$ in \eqref{eqn:RelativeLocalization} and before the weight restriction.

\begin{proposition}\label{prop:Vanishing}
Let $\ell = 0, \dots, s$. For $\vGa \notin \Gamma_{0,n+\ell}(\hcY^{(\ell)}/\hcD^{(\ell)},\hbeta,\bk)^{(0)}$, we have 
$$
	C^{(\ell)}_{\vGa} \big|_{\su_2 - f\su_1 = \su_4 = \cdots = \su_{3+\ell} = 0} = 0.
$$
In particular, we have
$$
	\langle \gamma_1^{(\ell)},\dots,\gamma_{n+\ell}^{(\ell)} \mid \gamma_{n+\ell+1}^{(\ell)},\dots,\gamma_{n+s}^{(\ell)} \rangle^{\hcY^{(\ell)}/\hcD^{(\ell)}, T_f}_{\hbeta,\bk^{(\ell)}} = \sum_{\vGa \in \Gamma_{0,n+\ell}(\hcY^{(\ell)}/\hcD^{(\ell)},\hbeta,\bk)^{(0)}} C^{(\ell)}_{\vGa} \big|_{\su_2 - f\su_1 = \su_4 = \cdots = \su_{3+\ell} = 0}.
$$
\end{proposition}

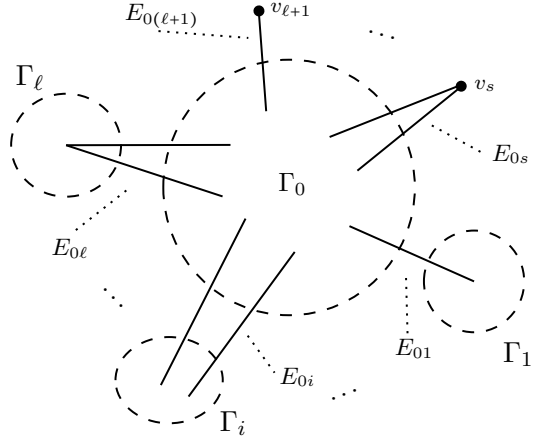
\begin{figure}[htbp]
	\centering

\tikzset{every picture/.style={line width=0.75pt}} 

\tikzset{every picture/.style={line width=0.75pt}} 

\tikzset{every picture/.style={line width=0.75pt}} 

\begin{tikzpicture}[x=0.75pt,y=0.75pt,yscale=-1,xscale=1]
	
	\draw [draw opacity=1 ]   (216.7,131.79) ;
	\draw [draw opacity=1 ] [dash pattern={on 0.84pt off 2.51pt}]  (255.38,171.35) -- (256.18,201.28) ;
	\draw  [draw opacity=1 ][dash pattern={on 4.5pt off 4.5pt}] (264.4,173.07) .. controls (264.4,158.9) and (275.72,147.41) .. (289.68,147.41) .. controls (303.64,147.41) and (314.96,158.9) .. (314.96,173.07) .. controls (314.96,187.25) and (303.64,198.74) .. (289.68,198.74) .. controls (275.72,198.74) and (264.4,187.25) .. (264.4,173.07) -- cycle ;
	\draw  [draw opacity=1 ][dash pattern={on 4.5pt off 4.5pt}] (132.95,125.66) .. controls (132.95,90.04) and (161.35,61.18) .. (196.38,61.18) .. controls (231.41,61.18) and (259.81,90.04) .. (259.81,125.66) .. controls (259.81,161.27) and (231.41,190.14) .. (196.38,190.14) .. controls (161.35,190.14) and (132.95,161.27) .. (132.95,125.66) -- cycle ;
	\draw [draw opacity=1 ]   (226.84,145.18) -- (289.68,173.07) ;
	\draw  [draw opacity=1 ][dash pattern={on 4.5pt off 4.5pt}] (108.84,222.92) .. controls (108.84,210.9) and (120.95,201.16) .. (135.89,201.16) .. controls (150.83,201.16) and (162.94,210.9) .. (162.94,222.92) .. controls (162.94,234.94) and (150.83,244.68) .. (135.89,244.68) .. controls (120.95,244.68) and (108.84,234.94) .. (108.84,222.92) -- cycle ;
	\draw  [draw opacity=1 ][dash pattern={on 4.5pt off 4.5pt}] (56.21,104.33) .. controls (56.21,89.98) and (68.62,78.35) .. (83.93,78.35) .. controls (99.23,78.35) and (111.64,89.98) .. (111.64,104.33) .. controls (111.64,118.67) and (99.23,130.3) .. (83.93,130.3) .. controls (68.62,130.3) and (56.21,118.67) .. (56.21,104.33) -- cycle ;
	\draw [draw opacity=1 ]   (145.7,231.19) -- (199.23,158.2) ;
	\draw [draw opacity=1 ]   (131.71,225.4) -- (174.43,141.38) ;
	\draw [draw opacity=1 ]   (83.93,104.33) -- (166.84,104.18) ;
	\draw [draw opacity=1 ]   (83.93,104.33) -- (162.84,130.18) ;
	\draw [draw opacity=1 ] [dash pattern={on 0.84pt off 2.51pt}]  (88.84,147.18) -- (116.84,125.18) ;
	\draw [draw opacity=1 ] [dash pattern={on 0.84pt off 2.51pt}]  (175.5,200.45) -- (189.08,219.16) ;
	\draw  [draw opacity=1 ][fill=black,fill opacity=1 ] (180.59,38.73) .. controls (179.39,38.4) and (178.68,37.15) .. (179.01,35.94) .. controls (179.34,34.74) and (180.59,34.03) .. (181.8,34.36) .. controls (183,34.69) and (183.71,35.94) .. (183.38,37.15) .. controls (183.05,38.35) and (181.8,39.06) .. (180.59,38.73) -- cycle ;
	\draw  [draw opacity=1 ][fill=black,fill opacity=1 ] (280.94,74.76) .. controls (280.82,73.52) and (281.73,72.41) .. (282.98,72.29) .. controls (284.22,72.17) and (285.33,73.08) .. (285.45,74.33) .. controls (285.57,75.57) and (284.66,76.68) .. (283.41,76.8) .. controls (282.17,76.92) and (281.06,76.01) .. (280.94,74.76) -- cycle ;
	\draw [draw opacity=1 ][fill opacity=1 ]   (181,35.36) -- (184.84,87.18) ;
	\draw [draw opacity=1 ]   (280.94,74.76) -- (216.84,100.18) ;
	\draw [draw opacity=1 ]   (281.94,75.76) -- (230.84,117.18) ;
	\draw [draw opacity=1 ] [dash pattern={on 0.84pt off 2.51pt}]  (130.51,46.73) -- (175.84,50.18) ;
	\draw [draw opacity=1 ] [dash pattern={on 0.84pt off 2.51pt}]  (268.38,94.35) -- (293.84,108.18) ;
	
	\draw (189.48,116.92) node [anchor=north west][inner sep=0.75pt]    {$\Gamma _{0}$};
	\draw (248.18,200.68) node [anchor=north west][inner sep=0.75pt]  [font=\small,opacity=1 ]  {$E_{01}$};
	\draw (190,215.05) node [anchor=north west][inner sep=0.75pt]  [font=\small ,opacity=1 ]  {$E_{0i}$};
	\draw (75.08,150.28) node [anchor=north west][inner sep=0.75pt]  [font=\small,opacity=1 ]  {$E_{0\ell}$};
	\draw (102.88,173.43) node [anchor=north west][inner sep=0.75pt]  [opacity=1 ,rotate=-44.81]  {$\dots$};
	\draw (214.99,231.05) node [anchor=north west][inner sep=0.75pt]  [opacity=1 ,rotate=-339.35]  {$\dots$};
	\draw (159.94,237.15) node [anchor=north west][inner sep=0.75pt]  [font=\large,opacity=1 ]  {$\Gamma_{i}$};
	\draw (302.89,202.06) node [anchor=north west][inner sep=0.75pt]  [font=\large,opacity=1 ]  {$\Gamma_{1}$};
	\draw (57,62.4) node [anchor=north west][inner sep=0.75pt]  [font=\large,opacity=1 ]  {$\Gamma_{\ell }$};
	\draw (235.06,44.08) node [anchor=north west][inner sep=0.75pt]  [opacity=1 ,rotate=-15.63]  {$\dots$};
	\draw (185.59,31.13) node [anchor=north west][inner sep=0.75pt]  [font=\footnotesize,opacity=1 ]  {$v_{\ell +1}$};
	\draw (288.45,69.73) node [anchor=north west][inner sep=0.75pt]  [font=\footnotesize,opacity=1 ]  {$v_{s}$};
	\draw (112.08,31.28) node [anchor=north west][inner sep=0.75pt]  [font=\small,opacity=1 ]  {$E_{0(\ell +1)}$};
	\draw (297.1,102.5) node [anchor=north west][inner sep=0.75pt]  [font=\small,opacity=1 ]  {$E_{0s}$};

\end{tikzpicture}

	\caption{Partitions of $V(\Gamma)$ and $E(\Gamma)$.}
	\label{fig:SubgraphsOuter}
\end{figure}

We prove the Proposition \ref{prop:Vanishing} in the rest of this subsection. We first note that by our choice of insertions, any decorated graph that does not satisfy condition \ref{item:Insertion} of Definition \ref{def:ContributingGraphs} will have zero contribution. Now let $\vGa \in \Gamma_{0,n+\ell}(\hcY^{(\ell)}/\hcD^{(\ell)},\hbeta,\bk^{(\ell)})$ such that it satisfies condition \ref{item:Insertion}. We introduce some notation. We partition the vertex set $V(\Gamma)$ as
$$
	V(\Gamma) = V_0 \sqcup V_1 \sqcup \cdots \sqcup V_\ell \sqcup \{\tv_{\ell+1}, \dots, \tv_{s} \}
$$
where:
\begin{itemize}
	\item $V_0 := \{v \in V(\Gamma) : \vf(v) \in \iota(\Sigma(3))\}$. 
	
	\item For $i = 1, \dots, \ell$, $V_i := \{ v\in V(\Gamma): \vf(v) =  \tsi_i\}$.
	
	\item For $i = \ell+1, \dots, s$, $\tv_i$ is the unique vertex with label $\vf(\tv_i) = \tsi_i$. (In Section \ref{sect:Moduli}, it is denoted by $v_{i-\ell}(\vGa)$.)
\end{itemize}
Note that condition \ref{item:Insertion} of Definition \ref{def:ContributingGraphs} implies that $V_i$ is non-empty for $i = 1, \dots, \ell$. We partition the edge set $E(\Gamma)$ as
$$
	E(\Gamma) = E_0 \sqcup E_{01} \sqcup \cdots \sqcup E_{0\ell} \sqcup E_{0(\ell+1)} \sqcup \cdots \sqcup E_{0s}
$$
where:
\begin{itemize}
	\item $E_0$ is the subset of edges that connects two vertices in $V_0$. 

	\item For $i = 1, \dots, \ell$, $E_{0i}$ is the subset of edges between $V_0$ and $V_i$.
	
	\item For $i = \ell+1, \dots, s$, $E_{0i}$ is the subset of edges between $V_0$ and $\tv_i$.
\end{itemize}
Note, in particular, that there are no edges between any two vertices in $V(\Gamma) \setminus V_0$. Let $\Gamma_0$ denote the subgraph $(V_0, E_0)$. See Figure \ref{fig:SubgraphsOuter} for an illustration of these partitions.
We further define the following subsets of flags:
\begin{itemize}
	\item $F_0 := \{(e,v) \in F(\Gamma): v \in V_0\}$.
	
	\item For $i = 1, \dots, \ell$, $F_i := \{(e,v) \in F(\Gamma): v \in V_i\}$.
	
	\item For $i = \ell+1, \dots, s$, $F_i := \{(e,\tv_i) \in F(\Gamma)\}$.

\end{itemize}

\begin{figure}[htbp]
	\centering

	\tikzset{every picture/.style={line width=0.75pt}} 

	\tikzset{every picture/.style={line width=0.75pt}} 

	\tikzset{every picture/.style={line width=0.75pt}} 
	
	\begin{tikzpicture}[x=0.75pt,y=0.75pt,yscale=-1,xscale=1]
		
		\draw    (140,154) ;
		\draw    (143.78,105.56) -- (283.16,137.39) ;
		\draw [shift={(283.16,137.39)}, rotate = 12.86] [color={rgb, 255:red, 0; green, 0; blue, 0 }  ][fill={rgb, 255:red, 0; green, 0; blue, 0 }  ][line width=0.75]      (0, 0) circle [x radius= 2.01, y radius= 2.01]   ;
		\draw  [dash pattern={on 0.84pt off 2.51pt}]  (226.78,158.56) -- (227.78,190.56) ;
		\draw  [dash pattern={on 4.5pt off 4.5pt}] (244.94,137.39) .. controls (244.94,116.29) and (262.05,99.18) .. (283.16,99.18) .. controls (304.26,99.18) and (321.37,116.29) .. (321.37,137.39) .. controls (321.37,158.5) and (304.26,175.61) .. (283.16,175.61) .. controls (262.05,175.61) and (244.94,158.5) .. (244.94,137.39) -- cycle ;
		\draw  [dash pattern={on 4.5pt off 4.5pt}] (44,145.28) .. controls (44,101.25) and (79.69,65.56) .. (123.72,65.56) .. controls (167.75,65.56) and (203.44,101.25) .. (203.44,145.28) .. controls (203.44,189.31) and (167.75,225) .. (123.72,225) .. controls (79.69,225) and (44,189.31) .. (44,145.28) -- cycle ;
		\draw   (116.68,92.95) .. controls (121.06,83.55) and (136.74,81.57) .. (151.71,88.54) .. controls (166.67,95.5) and (175.25,108.77) .. (170.88,118.18) .. controls (166.5,127.58) and (150.82,129.56) .. (135.85,122.59) .. controls (120.89,115.62) and (112.3,102.35) .. (116.68,92.95) -- cycle ;
		\draw   (95.89,177.56) .. controls (95.89,168.02) and (109.27,160.28) .. (125.78,160.28) .. controls (142.29,160.28) and (155.67,168.02) .. (155.67,177.56) .. controls (155.67,187.11) and (142.29,194.84) .. (125.78,194.84) .. controls (109.27,194.84) and (95.89,187.11) .. (95.89,177.56) -- cycle ;
		\draw    (125.78,177.56) -- (283.16,137.39) ;
		\draw  [dash pattern={on 0.84pt off 2.51pt}]  (132.78,131.56) -- (128.78,155.56) ;
		
		\draw (107,238.4) node [anchor=north west][inner sep=0.75pt]    {$\Gamma_{i}^c$};
		\draw (216.42,196.45) node [anchor=north west][inner sep=0.75pt]    {$E_{0i}$};
		\draw (273,188.4) node [anchor=north west][inner sep=0.75pt]    {$\Gamma_{i}$};
		\draw (277,143.4) node [anchor=north west][inner sep=0.75pt]    {$\tv_{i}$};

	\end{tikzpicture}

	\caption{Subgraphs $\Gamma_i$ and $\Gamma_i^c$. In the situation of Remark \ref{rem:U3+iVanish}, the unique vertex $\tv_i$ in $\Gamma_i$ is connected to each component of $\Gamma^c_i$ by an edge in $E_{0i}$.}
	\label{fig:SpecialCase}
\end{figure}
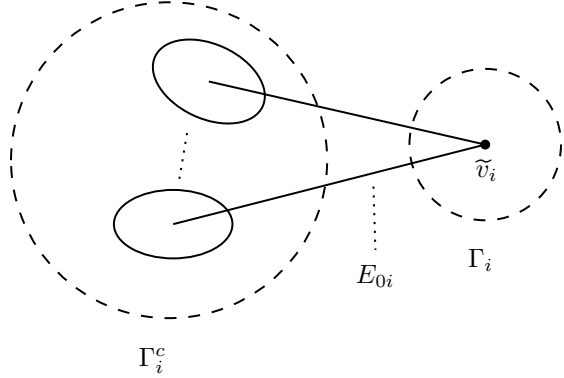

Moreover, for $i = 1, \dots, \ell$, let $V_i^c := V(\Gamma) \setminus V_i$ and $E_i^c := E(\Gamma) \setminus E_{01}$. Let $\Gamma_i^c$ denote the subgraph $(V_i^c, E_i^c)$. See Figure \ref{fig:SpecialCase}. We have
\begin{equation}\label{eqn:FiComplement}
	|F(\Gamma) \setminus F_i| = 2|E_i^c| + |E_{0i}|.
\end{equation} 
Let
$$
	c_i = |V_i^c| - |E_i^c|
$$
denote the number of connected components of $\Gamma_i^c$, which is upper-bounded by $|E_{0i}|$ since $V_i$ is non-empty. We prove the following preparatory lemma.


\begin{lemma}\label{lem:U3+iVanish}
For $\vGa \in \Gamma_{0,n+\ell}(\hcY^{(\ell)}/\hcD^{(\ell)},\hbeta,\bk^{(\ell)})$ satisfying condition \ref{item:Insertion} of Definition \ref{def:ContributingGraphs}, unless $c_i = |E_{0i}|$ for all $i = 1, \dots, \ell$, we have
$$
	C^{(\ell)}_{\vGa} \big|_{\su_4 = \dots = \su_{3+\ell} = 0} = 0.
$$
\end{lemma}

\begin{proof}
We consider the power of $\su_{i+3}$ in $C^{(\ell)}_{\vGa}$ for $i = 1, \dots, \ell$. By a computation analogous to \cite[Appendix B]{LY22} (cf. Lemma \ref{lem:WtExternalEdge}), each edge in $E(\Gamma_i^c)$ or vertex in $V(\Gamma_i^c)$ contributes a power of $-1$, each flag not in $F_i$ contributes a power of $1$, and there are no other contributions. Then, by \eqref{eqn:FiComplement}, the total power of $\su_{i+3}$ is 
$$
	-|E(\Gamma_i^c)|-|V(\Gamma_i^c)|+|F(\Gamma_i^c)| = |E_{0i}|-c_i.
$$
Taking all $i = 1, \dots, \ell$ into account, we obtain the lemma statement. 
\end{proof}

\begin{remark}\label{rem:U3+iVanish}
The condition $c_i = |E_{0i}|$, $i= 1,\dots,\ell$, of Lemma \ref{lem:U3+iVanish} means that $|V_i|=1$ and the unique vertex in $V_i$ is connected to each component of $\Gamma_i^c$ by an edge in $E_{0i}$. We denote the unique vertex in $V_i$ by $\tv_i$ as in condition \ref{item:LocalVertex} of Definition \ref{def:ContributingGraphs}. See Figure \ref{fig:SpecialCase}.
\end{remark}

We now prove the main statement on vanishing.

\begin{proof}[Proof of Proposition \ref{prop:Vanishing}]
	

Let $\vGa \in \Gamma_{0,n+\ell}(\hcY^{(\ell)}/\hcD^{(\ell)},\hbeta,\bk^{(\ell)})$ such that it satisfies condition \ref{item:Insertion} of Definition \ref{def:ContributingGraphs}. In view of Lemma \ref{lem:U3+iVanish}, we also assume that $c_i = |E_{0i}|$ for all $i = 1, \dots, \ell$. We show that unless $\vGa$ satisfies all of conditions \ref{item:RelEdge}, \ref{item:LocalEdge}, and \ref{item:LocalVertex}, we have $C^{(\ell)}_{\vGa} \big|_{\su_2 - f\su_1 = \su_4 = \cdots = \su_{3+\ell} = 0} = 0$. 
We consider the power of $\su_2-f\su_1$ in $C^{(\ell)}_{\vGa} \big|_{\su_4 = \cdots = \su_{3+\ell} = 0}$. Since $f \in \bQ$ is generic (Assumption \ref{assump:Parallel}), the powers may only come from the vertices $\tv_1, \dots, \tv_s$ and their incident edges and associated flags, as well as the insertions.




We first analyze the powers of $\su_2-f\su_1$ associated to $\tv_i$ for each $i=1,\dots,\ell$, again using results analogous to \cite[Appendix B]{LY22}. We show in several cases below that if either condition \ref{item:LocalEdge} or \ref{item:LocalVertex} is not satisfied by $i$, there is a strictly positive contribution of powers. In this situation, either $|E_{0i}| >1$, or $|E_{0i}|=1$ and $S_{\tv_i}$ contains a marking other than $n+i$. In either case, $\tv_i$ is a stable vertex. Moreover, by our choice of lifts $\gamma_1^{(\ell)},\dots,\gamma_{n}^{(\ell)}$ in Section \ref{sect:Insertions}, we must have $k_j \neq 1$ for all $j \in S_{\tv_i} \setminus \{n+i\}$ in order for $C^{(\ell+1)}_{\vGa}$ to be non-zero.
\begin{itemize}
	
	\item \textit{Case I: $\tk_i = 1$ and $\vk_{\tv_i} = (1)^{E_{0i} \cup S_{\tv_i}}$.} The latter condition means that $k_{(e, \tv_{i})} = 1$ for all $e \in E_{0i}$ and $S_{\tv_{i}} = \{n+i\}$. By definition, the insertion $\iota_{\tsi_{i}}^*(\gamma_{n+i}^{(\ell)})$ contributes a power of $1$. Moreover, $\tv_{i}$ contributes a power of $-2$, each edge $e \in E_{0i}$ contributes a power of $-1$, and the corresponding flag $(e, \tv_{i})$ contributes a power of $2$. Thus, the power of $\su_2 - f\su_1$ is
	$$
		|E_{0i}|-1 \ge 1.
	$$

	\item  \textit{Case II: $\tk_i = 1$ and $\vk_{\tv_i} \neq (1)^{E_{0i} \cup S_{\tv_i}}$.} Again $\iota_{\tsi_{i}}^*(\gamma_{n+i}^{(\ell)})$ contributes a power of $1$. Moreover, $\tv_{i}$ has no contribution, any edge $e \in E_{0i}$ with $k_{(e, \tv_{i})} = 1$ contributes a power of $-1$, and the corresponding flag $(e, \tv_{i})$ contributes a power of $2$. Thus, the power of $\su_2 - f\su_1$ is
	$$
		1 + |\{e \in E_{0i} : k_{(e, \tv_{i})} = 1\}| \ge 1.
	$$
	
	\item \textit{Case III: $\tk_{i} \neq 1$ and there exists $e \in E_{0i}$ with $k_{(e, \tv_{i})} = 1$.} In this case, there is no contribution from $\iota_{\tsi_{i}}^*(\gamma_{n+i}^{(\ell)})$. Moreover, $k_j \neq 1$ for all $j \in S_{\tv_{i}}$. Similar to above, there is no contribution from $\tv_{i}$, and the total contribution from edges in $E_{0i}$ and the associated flags is
	$$
		|\{e \in E_{0i} : k_{(e, \tv_{i})} = 1\}| \ge 1.
	$$

	\item \textit{Case IV: $\tk_{i} \neq 1$ and $k_{(e, \tv_{i})} \neq 1$ for all $e \in E_{0i}$.} This case is analogous to Case III except that $\tv_{i}$ always contributes a positive power of $\su_2 - f\su_1$ (cf. \cite[Lemma B.8]{LY22}).
	
\end{itemize}
Summarizing the cases, we see that in order for $\tv_i$ and its associated edges, flags, and insertions to contribute a zero power of $\su_2-f\su_1$, both conditions \ref{item:LocalEdge} and \ref{item:LocalVertex} must be satisfied by $i$.



The computation of the powers of $\su_2-f\su_1$ associated to $\tv_i$ for each $i=\ell+1,\dots,s$ follows from a parallel case analysis. We see that if condition \ref{item:RelEdge} is not satisfied by $i$, i.e. $l(\mu_i(\vGa))>1$, there is a strictly positive contribution of powers. In this situation, the nontrivial term $A^i_v$ in \eqref{eqn:RelativeLocalization} contributes powers. More specifically, each flag $(e, \tv_i)$ with $k_{(e,\tv_i)}=1$ contributes a power of $2$. If $\vk_{\tv_i} = (1)^{E_{0i} \cup \{n+i\}}$, 
the vertex term contributes a power of $-2$. On the other hand, if none of the components of $\vk_{\tv_i}$ is $1$, the vertex term contributes a positive power.
\end{proof}

\subsection{Open/relative correspondence}
We now prove the first equality in Theorem \ref{thm:OpenClosedStatement} which is a correspondence between open invariants of $(\cX, \cL, f)$ and relative invariants of $(\hcY,\hcD)$. The correspondence essentially follows from the computation in \cite{FLT12} (when each brane $\cL_i$ has integral framing $f_i \in \bZ$). We provide a proof for completeness.


\begin{proposition}\label{prop:OpenRelativeStatement}
Under the setup of Theorem \ref{thm:OpenClosedStatement}, we have the open/relative correspondence
$$
	\langle \gamma_1,\dots,\gamma_n  \rangle^{\cX, \cL, T_f}_{\beta',\bd} 
	= (-1)^{\sum_{i = 1}^s (\ceil{\frac{d_i}{\fa_i}} -1)} \sum_{\hbeta \in \Eff(\hY), \pi(\hbeta) = \beta'}\langle \gamma_1^{(0)},\dots,\gamma_n^{(0)} \mid \gamma_{n+1}^{(0)},\dots,\gamma_{n+s}^{(0)} \rangle^{\hcY/\hcD, T_f}_{\hbeta,\bk}.
$$
\end{proposition}

\begin{proof}
By Proposition \ref{prop:Vanishing}, for any $\hbeta \in \Eff(\hY)$ with $\pi(\hbeta) = \beta'$, only the decorated graphs in the set
$
	\Gamma_{0,n}(\hcY/\hcD,\hbeta,\bk)^{(0)}
$
contribute to the relative invariant. We set up a one-to-one correspondence
$$
	\epsilon: \Gamma_{0,n+s}(\cX,\beta) \rightarrow \bigsqcup_{\hbeta \in \Eff(\hY), \pi(\hbeta) = \beta'} \Gamma_{0,n}(\hcY/\hcD,\hbeta,\bk)^{(0)}
$$
as follows: 
Given $\vGa\in\Gamma_{0,n+s}(\cX,\beta)$, to get the corresponding graph $\epsilon(\vGa)$, 
for $i = 1, \dots, s$, we replace the ($n+i$)-th marked point in $\Gamma_{0,n+s}(\cX,\beta)$ with marking $v_{n+i}=\vs(n+i)$ by a new vertex $\tv_i$ and a new edge $e_i$ in $\epsilon(\vGa)$, with labels $\vf(\tv_i) = \tsi_i$, $\vf(e_{i})=\iota(\tau_{i})$, degree $\vd(e_{i})=(d_{i},\lambda_{i})$, and twisting $\vk(n+i)=\tk_i\in G_{\tsi_i}.$ The twistings of new flags are $k_{(e_i,v_{n+1})}=h(d_i,\lambda_i)$, $k_{(e_i,\tv_i)}=\tk_i$. See Figure \ref{fig:open-relative graph}. The labeling and twisting data of the original vertices and edges in $\vGa$ are kept, since the FTCY graph of $\cX$ is a subgraph of the FTCY graph $\sGa = \sGa^{(0)}$ of $(\hcY, \hcD)$. The degree map $\vd$ of $\vGa$ and $d_1, \dots, d_s$ determine a class $\hbeta \in \Eff(\hY)$ with $\pi(\hbeta) = \beta'$, and we have $\epsilon(\vGa) \in \Gamma_{0,n}(\hcY/\hcD,\hbeta,\bk)^{(0)}$.

\begin{figure}[htbp]
\tikzset{every picture/.style={line width=0.75pt}} 
\begin{tikzpicture}[x=0.75pt,y=0.75pt,yscale=-1,xscale=1]
	
	\draw    (17.46,103.96) -- (54.39,103.96) ;
	\draw [shift={(54.39,103.96)}, rotate = 0] [color={rgb, 255:red, 0; green, 0; blue, 0 }  ][fill={rgb, 255:red, 0; green, 0; blue, 0 }  ][line width=0.75]      (0, 0) circle [x radius= 3.35, y radius= 3.35]   ;
	\draw    (54.39,103.96) -- (86.7,149.76) ;
	\draw [shift={(86.7,149.76)}, rotate = 54.8] [color={rgb, 255:red, 0; green, 0; blue, 0 }  ][fill={rgb, 255:red, 0; green, 0; blue, 0 }  ][line width=0.75]      (0, 0) circle [x radius= 3.35, y radius= 3.35]   ;
	\draw    (86.7,149.76) -- (50.93,185.6) ;
	\draw [shift={(50.93,185.6)}, rotate = 134.94] [color={rgb, 255:red, 0; green, 0; blue, 0 }  ][fill={rgb, 255:red, 0; green, 0; blue, 0 }  ][line width=0.75]      (0, 0) circle [x radius= 3.35, y radius= 3.35]   ;
	\draw    (18.62,185.6) -- (50.93,185.6) ;
	\draw    (86.7,149.76) -- (103.38,149.6) ;
	
	\draw   (207.75,142.38) -- (247.44,142.38) -- (247.44,135.01) -- (270.58,147.75) -- (247.44,160.5) -- (247.44,153.13) -- (207.75,153.13) -- cycle ;
	
	\draw    (304.13,106.03) -- (341.06,106.03) ;
	\draw [shift={(341.06,106.03)}, rotate = 0] [color={rgb, 255:red, 0; green, 0; blue, 0 }  ][fill={rgb, 255:red, 0; green, 0; blue, 0 }  ][line width=0.75]      (0, 0) circle [x radius= 3.35, y radius= 3.35]   ;
	\draw    (341.06,106.03) -- (373.37,151.83) ;
	\draw [shift={(373.37,151.83)}, rotate = 54.8] [color={rgb, 255:red, 0; green, 0; blue, 0 }  ][fill={rgb, 255:red, 0; green, 0; blue, 0 }  ][line width=0.75]      (0, 0) circle [x radius= 3.35, y radius= 3.35]   ;
	\draw    (373.37,151.83) -- (337.6,187.68) ;
	\draw [shift={(337.6,187.68)}, rotate = 134.94] [color={rgb, 255:red, 0; green, 0; blue, 0 }  ][fill={rgb, 255:red, 0; green, 0; blue, 0 }  ][line width=0.75]      (0, 0) circle [x radius= 3.35, y radius= 3.35]   ;
	\draw    (305.29,187.68) -- (337.6,187.68) ;
	\draw    (419.79,151.67) -- (436.48,151.51) ;
	\draw    (373.37,151.83) -- (419.79,151.67) ;
	\draw [shift={(419.79,151.67)}, rotate = 359.8] [color={rgb, 255:red, 0; green, 0; blue, 0 }  ][fill={rgb, 255:red, 0; green, 0; blue, 0 }  ][line width=0.75]      (0, 0) circle [x radius= 3.35, y radius= 3.35]   ;
	
	\draw (34.45,157.18) node [anchor=north west][inner sep=0.75pt]   [align=left] {};
	\draw (84.9,164.37) node [anchor=north west][inner sep=0.75pt]    {$v_{n+i}$};
	\draw (110.22,139.75) node [anchor=north west][inner sep=0.75pt]    {$n+i$};
	\draw (321.12,159.26) node [anchor=north west][inner sep=0.75pt]   [align=left] {};
	\draw (365.79,164.15) node [anchor=north west][inner sep=0.75pt]    {$v_{n+i}$};
	\draw (443.29,143.07) node [anchor=north west][inner sep=0.75pt]    {$n+i$};
	\draw (413.24,158.79) node [anchor=north west][inner sep=0.75pt]    {$\tv_{i}$};
	\draw (391.02,130.65) node [anchor=north west][inner sep=0.75pt]    {$e_{i}$};

\end{tikzpicture}

\caption{New vertex and edge near marking $n+i$ in the construction of the decorated graph $\epsilon(\vGa)$ from $\vGa$.}
\label{fig:open-relative graph}
\end{figure}
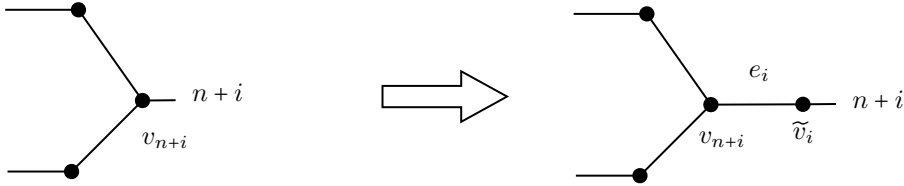

Now let $\vGa \in \Gamma_{0,n+s}(\cX,\beta)$ with underlying graph  $\Gamma$. Let $C_{\vGa}$ denote the contribution of $\vGa$ in the localization computation \eqref{equ:OpenInvariant} of the open invariant and before the weight restriction. Recall that $C^{(0)}_{\epsilon(\vGa)}$ denotes the contribution of $\epsilon(\vGa)$ in \eqref{eqn:RelativeLocalization} (for the case $\ell = 0$). We show that
$$
	C^{(0)}_{\epsilon(\vGa)} = (-1)^{\sum_{i = 1}^s (\ceil{\frac{d_i}{\fa_i}} -1)} C_{\vGa}
$$
by a comparison of the terms as follows:
\begin{itemize}
\item(Coefficients) By definition \eqref{eqn:cGamma}, we have
	$$
		c_{\epsilon(\vGa)} = c_{\vGa} \prod_{i=1}^s \frac{\frac{G_{v_{n+i}}}{r_{e_i,v_{n+i}}}\cdot\frac{G_{\tv_i}}{r_{(e_i,\tv_i)}}}{d_{e_i}|G_{e_i}|} = c_{\vGa} \prod_{i=1}^s  \frac{\fa_i \fr_{i}\fm_i}{d_i r_{(e_i,v_{n+i})}r_{(e_i,\tv_i)}}.
	$$
    \item (Edges) For each $e \in E(\Gamma)$, we have 
	$$\bh^{(0)}(e)=\bh(e).$$
	In addition, for $i=1, \dots, s$, we have by Lemma \ref{lem:WtExternalEdge} that
    $$
		\bh^{(0)}(e_i) 
		= \begin{cases}
			\displaystyle{     
			(-1)^{\frac{d_i}{\fa_i}-1} d_i\fm_i  \left( \frac{\su_{2,i}-f_i\su_{1,i}}{\fm_i} \right)^{-1} D_{d_i,1} }& \text{if } \tk_i = 1,\\
			\displaystyle{     (-1)^{\ceil{\frac{d_i}{\fa_i}}-1}d_i\fm_i  D_{d_{i},\lambda_{i}}  }& \text{if } \tk_i\neq 1.	
		\end{cases}
    $$

    \item(Flags) For each $(e,v)\in F(\Gamma)$, we have $$\bh^{(0)}(e,v)=\bh(e,v).$$
    In addition, let $i = 1, \dots, s$.
	For the flag $(e_{i},v_{n+i})$, we have
	$$
		\bh^{(0)}(e_{i},v_{n+i}) = \iota^*_{\sigma_{i}}(\gamma_{n+i}).
	$$
	For the flag $(e_{i},\tv_i)$, $\bh^{(0)}(e_{i},\tv_{i})$ can be accounted for by the term $\bh^{(0)}(\tv_{i})$ computed below.\\
    
	\item(Insertions) For $i=1,\dots,n$, the contribution of $\gamma_i^{(0)}$ is identified with that of $\gamma_i$. In addition, for $i = 1, \dots, s$, the contribution of $\gamma_{n+i}$ is accounted for above, and 
    $$
    	\iota_{\tv_i}^*(\gamma^{(0)}_{n+i}) = \begin{cases}
			\displaystyle{
			\sff_2^{(0)}(\tsi_i)=\frac{\su_{2,i}-f_i\su_{1,i}}{\fm_i} }&\text{ if } \tk_i=1 ,\\
			\displaystyle{
			\one_{\tk_i^{-1}} }&\text{ if } \tk_i\neq 1.
	    \end{cases}
    $$

	\item(Integral over vertices)  For any $v \in V(\Gamma) \setminus \{v_{n+1},\dots,v_{n+s}\}$, we have
	$$
		\int_{\Mbar_{0, \vk_v}(\cB G_v)} \frac{\bh^{(0)}(v)}{\prod_{e \in E_v} \left(\bw^{(0)}_{(e,v)} - \frac{\bar{\psi}_{(e,v)}}{r_{(e,v)}}\right)}  
		= \int_{\Mbar_{0, \vk_v}(\cB G_v)} \frac{\bh(v)}{\prod_{e \in E_v} \left(\bw_{(e,v)} - \frac{\bar{\psi}_{(e,v)}}{r_{(e,v)}}\right)}.
	$$
	Let $i = 1, \dots, s$.
	For the vertex $v_{n+i}$, note that the $(n+i)$-th marked point in $\vGa$ is replaced by the flag $(e_{i}, v_{n+i})$ in $\epsilon(\vGa)$, and the twisting $k_{n+i}$ in $\vGa$ is identified with $k_{(e_i, v_{n+i})}^{-1} = h(d_{i}, \lambda_{i})^{-1}$. Therefore, the integral
	$$
		\int_{\Mbar_{0, \vk_{v_{n+i}}}(\cB G_{v_{n+i}})} \frac{\bh^{(0)}(v_{n+i})}{\prod_{e \in E_{v_{n+i}}} \left(\bw^{(0)}_{(e,v_{n+i})}- \frac{\bar{\psi}_{(e,v_{n+i})}}{r_{(e,v_{n+i})}}\right)} 
	$$
	in $C^{(0)}_{\epsilon(\vGa)}$ is identified with
	$$
	r_{(e_{i}, v_{n+i})} \int_{\Mbar_{0, \vk_{v_{n+i}}}(\cB G_{v_{n+i}})} \frac{\bh(v_{n+i})}{\left(\frac{\su_{1,i}}{d_i} - \bar{\psi}_{n+i}\right) \cdot \prod_{e \in E_{v_{n+i}}} \left(\sw_{(e,v_{n+i})} - \frac{\bar{\psi}_{(e,v_{n+i})}}{r_{(e,v_{n+i})}}\right)}
	$$
	where the integral is the one in $C_{\vGa}$.
	Finally, for the univalent vertex $\tv_{i}$, we have by definition that $\bh^{(0)}(\tv_{i})  = \left( \bh^{(0)}(e_{i}, \tv_{i}) \right)^{-1}$, and
	\begin{align*}
		\int_{\Mbar_{0, \vk_{\tv_{i}}}(\cB G_{\tv_{i}})} \frac{\bh^{(0)}(\tv_{i})}{\bw^{(0)}_{(e,\tv_i)} - \frac{\bar{\psi}_{(e_{i},\tv_{i})}}{r_{(e_{i},\tv_{i})}}} 
		= \frac{r_{(e_{i}, \tv_{i})}}{|G_{\tv_{i}}|} \bh^{(0)}(\tv_{i})  = \frac{r_{(e_{i}, \tv_{i})}}{\fa_i\fm_{i}}  \left( \bh^{(0)}(e_{i}, \tv_{i})\right)^{-1}.
	\end{align*}
    
\end{itemize}
Summarizing the comparisons of the pieces above, we have
\begin{align*}
    C^{(0)}_{\epsilon(\vGa)}=& C_{\vGa}  \prod_{i=1}^s \left(
		\frac{\fa_i}{d_i r_{(e_i,v_{n+i})}r_{(e_i,\tv_i)}}  (-1)^{\ceil{\frac{d_i}{\fa_i}}-1}d_i\fm_i r_{(e_i,v_{n+i})} \frac{r_{(e_i,\tv_i)}}{\fa_i\fm_i}
    \right)
	\prod_{\tk_i=1} \left(\left( \frac{\su_{2,i}-f_i\su_{1,i}}{\fm_i}\right)^{-1} \left( \frac{\su_{2,i}-f_i\su_{1,i}}{\fm_i}\right)\right)\\
    =&(-1)^{\sum_{i=1}^s (\ceil{\frac{d_i}{\fa_i}}-1)}C_{\vGa}.
\end{align*}
\end{proof}

\subsection{Correspondence of intermediate relative invariants}
In this subsection, we prove Theorem \ref{thm:RelativeStatement}. Theorem \ref{thm:OpenClosedStatement} will then directly follow. 

\begin{proof}[Proof of Theorem \ref{thm:RelativeStatement}]
Let $\ell = 0, \dots, s-1$. By Proposition \ref{prop:Vanishing}, only the decorated graphs in the sets
$$
	\Gamma_{0,n+\ell}(\hcY^{(\ell)}/\hcD^{(\ell)},\hbeta,\bk)^{(0)}, \quad
	\Gamma_{0,n+\ell+1}(\hcY^{(\ell+1)}/\hcD^{(\ell+1)},\hbeta,\bk)^{(0)}
$$
contribute to the relative invariants. There is a natural one-to-one correspondence between the two sets
$$
	\epsilon:\Gamma_{0,n+\ell}(\hcY^{(\ell)}/\hcD^{(\ell)},\hbeta,\bk)^{(0)} \rightarrow \Gamma_{0,n+\ell+1}(\hcY^{(\ell+1)}/\hcD^{(\ell+1)},\hbeta,\bk)^{(0)}
$$
since the FTCY graphs $\sGa^{(\ell)}$ and $\sGa^{(\ell+1)}$ have the same underlying graph $\hGa$ and identified twisting data. We note specifically that for the marking $n+\ell+1$, the vertex $\vf \circ \vs(n+\ell+1) = \tsi_{\ell+1}$ changes from a univalent vertex to a $(4+\ell)$-valent vertex. In addition, the twisting $\vk(n+\ell+1) = \tk_{\ell+1}$ indexes an inertia component of the divisor $\hcD_{\ell+1}^{(\ell)}$ before but an inertia component of $\hcY^{(\ell+1)}$ after.



Now let $\vGa\in\Gamma_{0,n+\ell}(\hcY^{(\ell)}/\hcD^{(\ell)},\hbeta,\bk)^{(0)}$, and let $\Gamma$ be the common underlying graph of $\vGa$ and $\epsilon(\vGa)$. Recall from Lemma \ref{lem:U3+iVanish} that the total power of $\su_{4+\ell}$ in $C^{(\ell+1)}_{\epsilon(\vGa)}$ is $0$. We show that
$$
	C^{(\ell+1)}_{\epsilon(\vGa)} \big|_{\su_{4+\ell} = 0}
	= (-1)^{\ceil{\frac{d_{\ell+1}}{\fa_{\ell+1}}}-1} C^{(\ell)}_{\vGa}
$$
by a comparison of the terms as follows:
	\begin{itemize}
		\item(Coefficient) By definition \eqref{eqn:cGamma}, we have $c_{\epsilon(\vGa)}=c_{\vGa}$.
		
		\item(Edges) For the unique edge $e_{\ell+1}$ labeled by $\iota(\tau_{\ell+1})$, we have by Lemma \ref{lem:WtExternalEdge} that		
		$$
			\bh^{(\ell+1)}(e_{\ell+1}) \big|_{\text{$\su_{4+\ell}$} = 0} = \begin{cases}
				\displaystyle{     (-1)^{\frac{d_{\ell+1}}{\fa_{\ell+1}}}\left(\frac{\su_{1,\ell+1}}{\fa_{\ell+1}} \right)^{-1} \cdot\bh^{(\ell)}(e_{\ell+1}) }& \text{if } \tk_{\ell+1} = 1 \text{ or } \age(\tk_{\ell+1})=1,\\
				
				\displaystyle{(-1)^{\floor{\frac{d_{\ell+1}}{\fa_{\ell+1}}}}\cdot\bh^{(\ell)}(e_{\ell+1})}& \text{if } \age(\tk_{\ell+1})=2.
			\end{cases} 
		$$
		For any other edge $e\in E(\Gamma)$, we have
		$$
			\left(\su_{4+\ell} \bh^{(\ell+1)}(e) \right) \big|_{\su_{4+\ell} = 0} = \bh^{(\ell)}(e).
		$$

	\item(Flags) For each $(e,v)\in F(\Gamma)$ , we have $$\frac{\bh^{(\ell+1)}(e,v)}{\su_{4+\ell}}\bigg|_{\su_{4+\ell}=0}=\bh^{(\ell)}(e,v).$$
	
	\item(Insertions) 
	By the definition of insertions in Section \ref{sect:Insertions}, we have
	$$
	\iota^*_{\tv_{\ell+1}}(\gamma_{n+\ell+1}^{(\ell+1)}) = \begin{cases}
		\displaystyle{-\frac{\su_{1,i}}{\fa}}\iota^*_{\tv_{\ell+1}}(\gamma_{n+\ell+1}^{(\ell)}) & \text{if } \tk_{\ell+1}=1 \text{ or } \age(\tk_{\ell+1})=1,\\
		\iota^*_{\tv_{\ell+1}}(\gamma_{n+\ell+1}^{(\ell)}) &\text{if } \age(\tk_{\ell+1})=2.
	\end{cases}
	$$
	The contributions of the other insertions are identified.

	\item(Integral over vertices)  For any $v \in V(\vGa)^{(0)}$, we have
	$$
		\int_{\Mbar_{0, \vk_v}(\cB G_v)} \frac{\su_{4+\ell}\bh^{(\ell+1)}(v)}{\prod_{e \in E_v} \left(\bw^{(\ell+1)}_{(e,v)} - \frac{\bar{\psi}_{(e,v)}}{r_{(e,v)}}\right)} \bigg|_{\su_{4+\ell} = 0} 
		= \int_{\Mbar_{0, \vk_v}(\cB G_v)} \frac{\bh^{(\ell)}(v)}{\prod_{e \in E_v} \left(\bw^{(\ell)}_{(e,v)} - \frac{\bar{\psi}_{(e,v)}}{r_{(e,v)}}\right)}.
	$$

\end{itemize}
Summarizing the comparisons of the pieces above, we have for $\tk_{\ell+1}=1$ or $\age(\tk_{\ell+1})=1$ that
$$
	C^{(\ell+1)}_{\epsilon(\vGa)} \big|_{\su_{4+\ell} = 0} = C^{(\ell)}_{\vGa} \cdot \ (-1)^{\frac{d_{\ell+1}}{\fa_{\ell+1}}} \left(\frac{\su_{1,\ell+1}}{\fa} \right)^{-1} \cdot \left( -\frac{\su_{1,\ell+1}}{\fa} \right)  
	=(-1)^{\frac{d_{\ell+1}}{\fa_{\ell+1}}-1}C^{(\ell)}_{\vGa},
$$
and for $\age(\tk_{\ell+1})=2$ that
$$
	C^{(\ell+1)}_{\epsilon(\vGa)} \big|_{\su_{4+\ell} = 0} =C^{(\ell)}_{\vGa} \cdot (-1)^{\floor{\frac{d_{\ell+1}}{\fa_{\ell+1}}}}=C^{(\ell)}_{\vGa}\cdot (-1)^{\ceil{\frac{d_{\ell+1}}{\fa_{\ell+1}}}-1}.
$$
\end{proof}


Finally, Theorem \ref{thm:OpenClosedStatement} follows from Proposition \ref{prop:OpenRelativeStatement} and an inductive application of Theorem \ref{thm:RelativeStatement}.


\section{Application to BPS integrality}\label{sect:BPS}
In this section, in the case where $\cX$ is smooth and each brane $\cL_i$ has integral framing $f_i \in \bZ$, we apply the correspondence of open and closed Gromov-Witten invariants (Theorem \ref{thm:OpenClosedStatement}) to a correspondence of BPS invariants. The integrality of the open BPS invariants then implies the integrality of the closed BPS invariants. This section is a generalization of \cite[Section 4]{Yu24} which concerns the single-brane ($s=1$) case.

\subsection{Open and closed BPS invariants}\label{sect:OpenBPS}
We take the following specifications:
\begin{itemize}
    \item $\cX = X$ is a smooth toric Calabi-Yau 3-fold.

    \item $f \in \bQ$ is chosen such that the induced framing $f_i$ of each brane $\cL_i$ is an integer.
\end{itemize}
Then, each $\cL_i = L_i$ is smooth and open stable maps to $(X, L)$ has trivial monodromy along the boundary components. Moreover, the corresponding closed geometry $\tcX = \tX$ is a smooth toric Calabi-Yau ($3+s$)-fold.

As in Section \ref{sect:OpenGW}, let $\beta \in \Eff(X)$ and $d_1, \dots, d_s \in \bZ_{\ge 1}$. In this case, we have
$$
    \bd = ((d_1, 1), \dots, (d_s, 1)).
$$
Let
$$
	\beta' = \beta + \sum_{i = 1}^s d_i[B_i] \qquad \in \Eff(X,L).
$$
Consider the genus-zero, degree-$(\beta', \bd)$ open Gromov-Witten invariant
$$
    N^{X, L, f}_{\beta', d_1, \dots, d_s} := \langle \phantom{1}  \rangle^{X, L, T_f}_{\beta',\bd}
$$
without any interior insertions (i.e. $n = 0$), as defined in Definition \ref{def:OpenGW}. 

Open BPS invariants may be defined from the Gromov-Witten invariants from a resummation formula of Labastida, Mari\~no, Ooguri, and Vafa (LMOV) \cite{OV00,LM00,LMV00,MV02}. In the genus-zero case with maximal winding profile at each brane $\cL_i$, the LMOV formula is
\begin{equation}\label{eqn:LMOVOpen}
    N^{X, L, f}_{\beta', \bd} = (-1)^s \sum_{k \mid \beta, \bd} k^{s-3} n^{X, L, f}_{0,\frac{\beta'}{k}, \frac{\bd}{k}}.
\end{equation}
Here, we say that $k \in \bZ_{\ge 1}$ divides $\beta$ if $\frac{\beta}{k} \in H_2(X; \bZ)$, and $k$ divides $\bd$ if $k \mid d_i$ for each $i$, in which case we write
$$
    \frac{\bd}{k} = \left( \left(\frac{d_1}{k}, 1 \right), \dots, \left(\frac{d_s}{k}, 1 \right) \right).
$$
The coefficient
$$
    n^{X, L, f}_{\beta', \bd}
$$
is uniquely determined from the above and is referred to as the genus-zero, degree-$(\beta', \bd)$ \emph{open BPS invariant}
.

On the other hand, for the smooth toric Calabi-Yau ($3+s$)-fold $\tX$, let $\tbeta \in \Eff(\tX)$ and $\tgamma_1, \dots, \tgamma_s$ be defined as in Section \ref{sect:CorrData}. In this case, for $i = 1, \dots, s$,
$$
    \tgamma_i = [V(\trho_{R+2i-1})] [V(\trho_{i_2^{(\tau_i,\sigma_i)}})]
$$
can be viewed as a non-equivariant class in $H^4(\tX; \bZ)$. Consider the genus-zero, degree-$\tbeta$ Gromov-Witten invariant
$$
    N^{\tX}_{\tbeta}(\tgamma_1, \dots, \tgamma_s) := \langle \tgamma_1, \dots, \tgamma_s \rangle^{\tX, T_f}_{\tbeta},
$$
as defined in Definition \ref{def:ClosedGW}.

For genus-zero, $s$-pointed invariants of Calabi-Yau manifolds of dimension at least 4, closed BPS invariants may be defined from the Gromov-Witten invariants from a resummation formula of Klemm-Pandharipande \cite{KP08}, which generalizes the Aspinwall-Morrison multiple covering formula \cite{AM93}:
\begin{equation}\label{eqn:KPClosed}
    N^{\tX}_{\tbeta}(\tgamma_1, \dots, \tgamma_s) = \sum_{k \mid \tbeta} k^{s-3} n^{\tX}_{\frac{\tbeta}{k}}(\tgamma_1, \dots, \tgamma_s).
\end{equation}
The coefficient
$$
    n^{\tX}_{\tbeta}(\tgamma_1, \dots, \tgamma_s)
$$
is uniquely determined from the above and is referred to as the genus-zero, degree-$\tbeta$ \emph{closed BPS invariant}. We note that $k \mid \tbeta$ if and only if $k \mid \beta, \bd$.

\subsection{BPS correspondence and integrality}
Comparing the resummation formulas \eqref{eqn:LMOVOpen}, \eqref{eqn:KPClosed}, we may use the open/closed correspondence for Gromov-Witten invariants (Theorem \ref{thm:OpenClosedStatement}) to deduce the following correspondence for BPS invariants.

\begin{theorem}\label{thm:BPSCorrespondence}
We have
$$
    n^{X, L, f}_{\beta', \bd} = (-1)^s n^{\tX}_{\tbeta}(\tgamma_1, \dots, \tgamma_s).
$$
\end{theorem}

The integrality of the open BPS invariant
$$
    n^{X, L, f}_{\beta', \bd} \in \bZ
$$
is shown in \cite[Theorem 1.1]{Yu24}. Therefore, Theorem \ref{thm:BPSCorrespondence} implies the integrality of closed BPS invariants of $\tX$ with insertions $\tgamma_1, \dots, \tgamma_s$, which verifies the conjecture of Klemm-Pandharipande \cite[Conjecture 0]{KP08} in this case.

\begin{corollary}\label{cor:ClosedBPS}
We have
$$
    n^{\tX}_{\tbeta}(\tgamma_1, \dots, \tgamma_s) \in \bZ.
$$
\end{corollary}


\section{Example: $\bC^3$}\label{sect:example}

In this section, we demonstrate our constructions and results in the example $\cX = \bC^3$ relative to two Aganagic-Vafa branes with parallel framings. Moreover, we indicate a non-formal construction of the intermediate relative geometries.

\subsection{3-dimensional open geometry}
In $N_{\bR} \cong \bR^3$, the fan $\Sigma$ of $\cX = \bC^3$ consists of a single 3-cone $\sigma_{\{123\}}$ generated by the vectors
$$
	b_1 = (1,0,1), \quad b_2 = (0,1,1), \quad b_3 = (0,0,1).
$$
Here and in this rest of this section, we use the shorthand notation $\sigma_I$ (resp. $\tau_I$) for a index set $I$ to denote the cone such that $I'_{\sigma_I}$ (resp. $I'_{\tau_I}$) is equal to $I$. Let $s = 2$ and $\cL_1, \cL_2$ be Aganagic-Vafa branes with
$$
	\tau_1 = \tau_{\{23\}}, \qquad \tau_2 = \tau_{\{13\}}.
$$
We adopt the coordinate system under which
$$
	\su_1 = \sw(\tau_{\{23\}}, \sigma_{\{123\}}) = \su_{1,1}, \quad
	\su_2 = \sw(\tau_{\{13\}}, \sigma_{\{123\}}) = \su_{1,2} = \su_{2,1}, \quad
	-\su_1 - \su_2 = \sw(\tau_{\{12\}}, \sigma_{\{123\}}) = \su_{2,2}.
$$
Moreover, we choose the parallel framings
$$
	f = f_1 = 1, \qquad f_2 = -2.
$$
Note that we may expand the example by adding in a third brane corresponding to the 2-cone $\tau_{\{12\}}$, which would have the fractional framing $-\frac{1}{2}$ under the choice above.

\subsection{3-dimensional relative geometries}
\begin{figure}[h]
	\begin{tikzpicture}
		\draw (-2.3,-1.3) -- (-1,0);
		\draw (-1,0) -- (-1,2);
		\draw (-1,0) -- (1,0);
		\draw[dashed] (0.2,0.8) -- (2,-1);
		\draw[dashed] (-2,3) -- (-0.2,1.2);
		\node at (-1,0) {$\bullet$};
		\node at (1,0) {$\bullet$};
		\node at (-1,2) {$\bullet$};
		\node at (-0.8,0.2) {\small $v_0$};
		\node at (1.1,0.25) {\small $\tv_1$};
		\node at (-0.7,2.1) {\small $\tv_2$};
		\node[below] at (-0.5, 0) {\small $\su_1$};
		\node[below] at (0.5, 0) {\small $-\su_1$};
		\node[left] at (-1,1.5) {\small $-\su_2$};
		\node[left] at (-1,0.5) {\small $\su_2$};
		\node[left] at (-2.3,-1.3) {\small $-\su_1-\su_2$};
		\node[right] at (0.2,0.8) {\small $-\su_1+\su_2$};
		\node[right] at (2,-1) {\small $\su_1-\su_2$};
		\node[right] at (-0.2,1.3) {\small $\su_1-\su_2$};
		\node[left] at (-2,3) {\small $-\su_1+\su_2$};
	\end{tikzpicture}

	\caption{FTCY graph for the 3-dimensional formal relative geometry $(\hcY, \hcD)$.}
	\label{fig:ExDim3}
\end{figure}
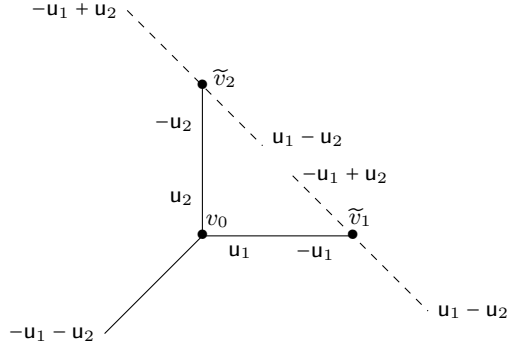

The FTCY graph $\vGa = (\hGa, \sfp, \sff_2, \sff_3)$ specifying the 3-dimensional formal relative geometry $(\hcY, \hcD)$ can be described as follows. The vertex set of the underlying graph is
$$
	V(\hGa) = \{v_0, \tv_1, \tv_2\}
$$
where $v_0$ corresponds to the 3-cone $\sigma_{\{123\}}$ and is the only trivalent vertex. In addition, we have
$$
	\sff_2(\tv_1) = \sff_2(\tv_2) = -\su_1 + \su_2, \qquad \sff_3(\tv_1) = \sff_3(\tv_2) = \su_1 - \su_2.
$$
See Figure \ref{fig:ExDim3} for an illustration.

The formal geometry may also be realized via a non-formal construction as follows. Let $\cY$ be the toric 3-fold specified by the fan in $N_{\bR}$ whose rays are spanned by
$$
	b_1 = (1,0,1), \quad b_2 = (0,1,1), \quad b_3 = (0,0,1), \quad b_4^{(0)} = (-1,-1,0)
$$
respectively and whose 3-cones are $\sigma_{\{123\}}, \sigma_{\{134\}}, \sigma_{\{234\}}$. Let $\cD$ be the toric divisor corresponding to the ray spanned by $b_4$. The the pair $(\cY, \cD)$ is log Calabi-Yau, and $(\hcY, \hcD)$ is the formal completion of $\cY$ along the closed subset of the toric 1-skeleton which is the union of the $T$-invariant lines corresponding to the edges in $\hGa$. The relative Gromov-Witten invariants of $(\cY, \cD)$ can be computed by localization (cf. \cite{Dolfen21}) and with our choice of insertions on $\cD$, by a vanishing argument similar to that in Section \ref{sect:Vanishing}, they are identified with the formal relative invariants of $(\hcY, \hcD)$.

\subsection{4-dimensional intermediate geometries}
\begin{figure}[h]
	\begin{tikzpicture}
		\draw (-2.3,-1.3) -- (-1,0);
		\draw (-1,0) -- (-1,2);
		\draw[dashed] (-2.6,2.3) -- (-1,2);
		\draw (-2.6, 0.3) -- (-1,0) -- (1,0) -- (2.3, 0.3);
		\draw (0.2,0.8) -- (2,-1);
		\draw[dashed] (-2,3) -- (-0.2,1.2);
		\node at (-1,0) {$\bullet$};
		\node at (1,0) {$\bullet$};
		\node at (-1,2) {$\bullet$};
		\node at (-0.8,0.2) {\small $v_0$};
		\node at (1.1,0.25) {\small $\tv_1$};
		\node at (-0.7,2.1) {\small $\tv_2$};
		\node[below] at (-0.5, 0) {\small $\su_1$};
		\node[below] at (0.5, 0) {\small $-\su_1$};
		\node[left] at (-1,1.5) {\small $-\su_2$};
		\node[left] at (-1,0.5) {\small $\su_2$};
		\node[left] at (-2.3,-1.3) {\small $-\su_1-\su_2-\su_4$};
		\node[right] at (0.2,0.8) {\small $-\su_1+\su_2$};
		\node[right] at (2,-1) {\small $\su_1-\su_2-\su_4$};
		\node[right] at (-0.2,1.3) {\small $\su_1-\su_2$};
		\node[left] at (-2,3) {\small $-\su_1+\su_2-\su_4$};
		\node[left] at (-2.6, 0.3) {\small $\su_4$};
		\node[left] at (-2.6, 2.3) {\small $\su_4$};
		\node[right] at (2.3,0.3) {\small $\su_1+\su_4$};
	\end{tikzpicture}

	\caption{FTCY graph for the 4-dimensional formal relative geometry $(\hcY^{(1)}, \hcD^{(1)})$.}
	\label{fig:ExDim4}
\end{figure}
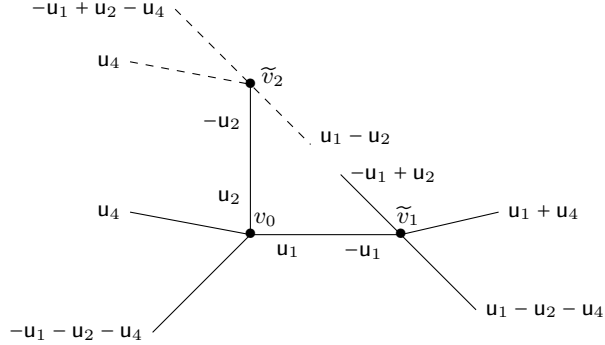

The FTCY graph $\vGa^{(1)}$ specifying the 4-dimensional formal intermediate geometry $(\hcY^{(1)}, \hcD^{(1)})$ has the same underlying vertex set $V(\hGa)$, except that $v_0$ and $\tv_1$ are 4-valent and $\tv_2$ is univalent. See Figure \ref{fig:ExDim4} for an illustration of the graph and the weights.

The formal geometry may also be realized via a non-formal construction as follows. Let $\cY^{(1)}$ be the toric 4-fold specified by the fan in $N_{\bR} \oplus \bR$ whose rays are spanned by
$$
	b_1^{(1)} = (1,0,1,0), \quad b_2^{(1)} = (0,1,1,0), \quad b_3^{(1)} = (0,0,1,0),
$$
$$
	b_4^{(1)} = (-1,-1,1,1), \quad b_5^{(1)} = (0,0,1,1), \quad b_6^{(1)} = (-1,-1,0,0)
$$
respectively and whose 4-cones are $\sigma_{\{1235\}}, \sigma_{\{2345\}}, \sigma_{\{1356\}}$. Let $\cD^{(1)}$ be the toric divisor corresponding to the ray spanned by $b_6^{(1)}$. The the pair $(\cY^{(1)}, \cD^{(1)})$ is log Calabi-Yau, and $(\hcY^{(1)}, \hcD^{(1)})$ is the formal completion of $\cY^{(1)}$ along the closed subset of the toric 1-skeleton which is the union of the $T_{(1)}$-invariant lines corresponding to the edges in $\hGa^{(1)}$. Similar to before, the relative Gromov-Witten invariants of $(\cY^{(1)}, \cD^{(1)})$ can be computed by localization and identified with the formal relative invariants of $(\hcY^{(1)}, \hcD^{(1)})$.

\subsection{5-dimensional closed geometries}
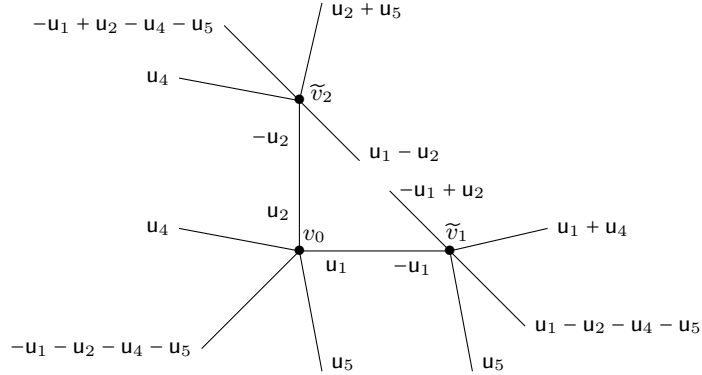
\begin{figure}[h]
	\begin{tikzpicture}
		\draw (-2.3,-1.3) -- (-1,0);
		\draw (-0.7,-1.6) -- (-1,0) -- (-1,2) -- (-0.7,3.3);
		\draw (-2.6,2.3) -- (-1,2);
		\draw (-2.6, 0.3) -- (-1,0) -- (1,0) -- (2.3, 0.3);
		\draw (1,0) -- (1.3,-1.6);
		\draw (0.2,0.8) -- (2,-1);
		\draw (-2,3) -- (-0.2,1.2);
		\node at (-1,0) {$\bullet$};
		\node at (1,0) {$\bullet$};
		\node at (-1,2) {$\bullet$};
		\node at (-0.8,0.2) {\small $v_0$};
		\node at (1.1,0.25) {\small $\tv_1$};
		\node at (-0.7,2.1) {\small $\tv_2$};
		\node[below] at (-0.5, 0) {\small $\su_1$};
		\node[below] at (0.5, 0) {\small $-\su_1$};
		\node[left] at (-1,1.5) {\small $-\su_2$};
		\node[left] at (-1,0.5) {\small $\su_2$};
		\node[left] at (-2.3,-1.3) {\small $-\su_1-\su_2-\su_4-\su_5$};
		\node[right] at (0.2,0.8) {\small $-\su_1+\su_2$};
		\node[right] at (2,-1) {\small $\su_1-\su_2-\su_4-\su_5$};
		\node[right] at (-0.2,1.3) {\small $\su_1-\su_2$};
		\node[left] at (-2,3) {\small $-\su_1+\su_2-\su_4-\su_5$};
		\node[left] at (-2.6, 0.3) {\small $\su_4$};
		\node[left] at (-2.6, 2.3) {\small $\su_4$};
		\node[right] at (2.3,0.3) {\small $\su_1+\su_4$};
		\node[right] at (-0.7, -1.5) {\small $\su_5$};
		\node[right] at (1.3, -1.5) {\small $\su_5$};
		\node[right] at (-0.7, 3.2) {\small $\su_2+\su_5$};
	\end{tikzpicture}

	\caption{FTCY graph for the 5-dimensional formal geometry $\hcY^{(2)}$.}
	\label{fig:ExDim5}
\end{figure}

In the fan $\tSi \subset \tN_{\bR} \cong N_{\bR} \oplus \bR^2$ of the 5-dimensional closed geometry $\tcX$, the rays are spanned by 
$$
	\tb_1 = (1,0,1,0,0), \quad \tb_2 = (0,1,1,0,0), \quad \tb_3 = (0,0,1,0,0),
$$
$$
	\tb_4 = (-1,-1,1,1,0), \quad \tb_5 = (0,0,1,1,0), \quad \tb_6 = (-1,-1,1,0,1), \quad \tb_7 = (0,0,1,0,1)
$$
respectively and the 5-cones are $\sigma_{\{12357\}}, \sigma_{\{23457\}}, \sigma_{\{13567\}}$. See Figure \ref{fig:ExDim5} for an illustration of the FTCY graph $\vGa^{(2)}$ specifying the formal geometry $\hcY^{(2)}$ which is the formal completion of $\tcX$ along $\tcX^1$. In this case, $\tcX$ is a smooth toric Calabi-Yau 5-fold and the discussion on BPS integrality in Section \ref{sect:BPS} can be applied.

\begin{remark}\rm{
In addition, we may work with the natural semi-projective partial compactification $\tcX'$ of $\tcX$ whose fan $\tSi'$ contains $\tSi$ as a subfan and has support the convex hull of the support of $\tSi$. By a vanishing argument similar to that in Section \ref{sect:Vanishing}, the closed invariants of $\tcX'$ with our choice of insertions are identified with the closed invariants of $\tcX$. Therefore, they fit into the correspondences in Theorem \ref{thm:OpenClosedStatement} and have the integrality property described in Corollary \ref{cor:ClosedBPS}. 
}\end{remark}

\begin{remark}\label{rem:OrbifoldExample}\rm{
If at the beginning we include the third Aganagic-Vafa brane $\cL_3$ with $\tau_3 = \tau_{12}$ into the boundary condition, the corresponding parallel framing would be $f_3 = -\frac{1}{2}$ and the resulting 6-dimensional closed geometry is an orbifold. For any general $f \in \bQ$, we have $f_1f_2f_3 = 1$ and it is not possible to obtain a smooth closed geometry.
}\end{remark}

 \nocite{*}
\bibliographystyle{plain}

\end{document}